\documentclass[10.5pt]{article}
\usepackage{amsmath, amssymb, amscd, amsthm, amsfonts}
\usepackage{graphicx}
\usepackage[hyperfootnotes=false]{hyperref}
\usepackage{enumitem}
\usepackage{esint}
\usepackage{booktabs}
\usepackage{tocloft} 
\usepackage{etoolbox}
\usepackage{xcolor}
\usepackage{marvosym}
\usepackage[title]{appendix}
\hypersetup{
    colorlinks=true, 
    linkcolor=blue,      
    citecolor=red!40!purple ,
    urlcolor=blue
}

\let\oldthebibliography\thebibliography
\let\endoldthebibliography\endthebibliography
\renewenvironment{thebibliography}[1]{%
  \oldthebibliography{#1}%
  \setlength{\itemsep}{4.3pt}%
  \setlength{\parsep}{0pt}%
  \setlength{\parskip}{0pt}%
}{%
  \endoldthebibliography
}

\makeatletter
\let\oldcite\cite
\renewcommand{\cite}[1]{{\fontsize{9}{12}\selectfont\oldcite{#1}}}
\makeatother

\title{Singular Global Compactness for Logarithmically Coupled \(p\)-Laplacian Systems}
\author{Wenzheng Hu\textsuperscript{1}}
\date{}

\newtheorem{theorem}{Theorem}[section]
\newtheorem{lemma}{Lemma}[section]

\newtheorem{remark}{Remark}[section]
\newtheorem{definition}{Definition}[section]
\newtheorem{proposition}{Proposition}[section]
\newtheorem{corollary}{Corollary}[section]
\newtheorem{example}{Example}[section]

\allowdisplaybreaks[1]
\usepackage{amsmath}
\numberwithin{equation}{section}
\begin{document}

\begingroup
\let\thefootnote\relax
\footnotetext{%
  \fontsize{10}{12}\selectfont
  \parbox[t]{\linewidth}{%
    \begin{tabular}{@{}l@{\hspace{0.5em}}p{0.85\textwidth}@{}}
      \raisebox{-0.1em}{\large\Letter} & Wenzheng Hu \\[4pt]
      & wzhu25@m.fudan.edu.cn \\[12pt]
      \textsuperscript{1} & School of Mathematical Sciences, Fudan University, Shanghai 200433, China.
    \end{tabular}%
  }%
}
\endgroup
\maketitle	
	
	\begin{center}
\Large\bfseries Abstract
\end{center}

We establish global compactness for logarithmically coupled \(p\)-Laplacian systems on \(\mathbb{R}^N\), where \(p=N>4\) and bounded positive potentials converge to positive constants at infinity. The finite-energy domain need not be closed under subtraction, obstructing direct profile subtraction. We develop a scalar localisation method using slowly expanding cutoffs and residual commutator estimates uniform in the regularisation parameter. As the regularisation vanishes, residual Palais–Smale sequences split into a solution of the original system and finitely many separating autonomous profiles, with a strongly vanishing remainder and splitting of the logarithmic densities. A converse construction characterises all levels of compactness failure. Applications include convergence of ground-state levels, a compactness–escape alternative for minimising sequences, and a directional response formula for the ground-state energy.



\vspace{0.3cm}
\tableofcontents

\section{Introduction}\label{sec:introduction}

Equations with logarithmic nonlinearities have attracted sustained attention over the past few decades. Following the logarithmic Schr\"odinger model of nonlinear wave mechanics introduced by Bia{\l}ynicki-Birula and Mycielski~\cite{BM76}, early analytical developments included Cazenave's study of stable solutions~\cite{Caz83}. Subsequent work addressed existence, multiplicity, and qualitative properties of stationary solutions; in particular, d'Avenia, Montefusco, and Squassina~\cite{DMS14} developed a nonsmooth variational approach to the logarithmic Schr\"odinger equation. Recent advances include the existence and concentration of normalised solutions for logarithmic \(p\)-Laplacian equations~\cite{SS25}, as well as convergence from power-type equations and log-concavity of positive solutions on convex domains~\cite{GMS26}. These developments highlight the continuing importance of variational methods and compactness arguments in the analysis of logarithmic problems. For further studies of PDEs with logarithmic nonlinearities, we refer the reader to \cite{APYZ26,CD25,CF24,CG18,Fer21} and the references therein.

Compactness plays a central role in variational approaches to nonlinear elliptic equations, allowing minimising and Palais--Smale sequences to yield solutions. On unbounded domains, however, bounded sequences may lose compactness through translation to infinity, even when local compactness is available. Lions' concentration--compactness principle~\cite{Lio84a,Lio84b} provides a systematic framework for analysing this loss through the alternatives of compactness, vanishing, and dichotomy. A complementary description is furnished by global compactness theorems, which resolve suitable Palais--Smale sequences into a weak limit and a finite collection of nontrivial profiles, together with a strongly vanishing remainder and an energy decomposition. Struwe established a foundational result of this type for elliptic problems with critical Sobolev growth~\cite{BC88, Bec93, BN83, Str84}, and Mercuri and Willem subsequently developed a corresponding theory for critical \(p\)-Laplacian problems~\cite{MW10}. These theories make it possible to identify the mechanisms of compactness failure and recover compactness at suitable energy levels. For further developments on compactness, concentration phenomena,
and quantitative profile decomposition in elliptic PDEs, we refer
the reader to \cite{DMMT22,DSW25,DT20,FG20,Han22,IMN15,JLX14}
and the references therein.

This is the second series in our study of logarithmically coupled \(p\)-Laplacian systems and a sequel work to \cite{Hu26}. We study the logarithmically coupled \(p\)-Laplacian system
\begin{equation}\label{eq:system}
\left\{
\begin{aligned}
 -\Delta_p u+a(x)|u|^{p-2}u
 &=\frac{p-2}{p}|u|^{p-4}u v^2\log(v^2+\varepsilon)
   +\frac{2}{p}|v|^{p-2}u\log(u^2+\varepsilon)
   +\frac{2|v|^{p-2}u^3}{p(u^2+\varepsilon)},
 &&\text{in }\mathbb R^N,\\[8pt]
 -\Delta_p v+b(x)|v|^{p-2}v
 &=\frac{p-2}{p}|v|^{p-4}v u^2\log(u^2+\varepsilon)
   +\frac{2}{p}|u|^{p-2}v\log(v^2+\varepsilon)
   +\frac{2|u|^{p-2}v^3}{p(v^2+\varepsilon)},
 &&\text{in }\mathbb R^N,
\end{aligned}
\right.
\end{equation}
where \(p=N>4\) and \(0\le\varepsilon\le\tfrac12\). The potentials \(a\) and \(b\)
are bounded, uniformly positive, and converge at infinity to
positive constants \(\bar a\) and \(\bar b\), respectively.
The case \(\varepsilon=0\) is the original system, with the
nonlinear terms interpreted by continuity at zero; positive
values of \(\varepsilon\) give its regularization.

For convenience, we introduce the interaction density
\[
 \ell_\varepsilon(s,t)
 =|s|^{p-2}t^2\log(t^2+\varepsilon)
  +|t|^{p-2}s^2\log(s^2+\varepsilon).
\]
Then \eqref{eq:system} takes the form
\[
 \begin{cases}
 -\Delta_p u+a(x)|u|^{p-2}u
   =\dfrac1p\,\partial_s\ell_\varepsilon(u,v),\\[4pt]
 -\Delta_p v+b(x)|v|^{p-2}v
   =\dfrac1p\,\partial_t\ell_\varepsilon(u,v).
 \end{cases}
\]
At \(\varepsilon=0\), the logarithmic terms are understood by continuous extension at zero. Here \(\Delta_p u=\operatorname{div}(|\nabla u|^{p-2}\nabla u)\), and the potential pair \(q=(a,b)\) satisfies
\[
a,b\in L^\infty(\mathbb{R}^N),\qquad
\operatorname*{ess\,inf}_{\mathbb{R}^N}a>0,\qquad
\operatorname*{ess\,inf}_{\mathbb{R}^N}b>0,
\]
together with
\[
\operatorname*{ess\,sup}_{|x|>R}
\bigl(|a(x)-\bar a|+|b(x)-\bar b|\bigr)
\longrightarrow0
\qquad\text{as }R\to\infty,
\]
for positive constants \(\bar a,\bar b\). Unlike the confining continuous setting considered in \cite{Hu26}, the present setting permits profiles to escape to infinity, where they are governed by the autonomous system with potential pair \(\bar q=(\bar a,\bar b)\).

The system also admits a phenomenological interpretation as a
stationary model for two interacting fields in an inhomogeneous
medium. The \(p\)-Laplacian describes nonlinear diffusion, while
the potentials represent spatially varying energetic costs.
The mixed logarithmic terms introduce an amplitude-dependent
interaction between the two components. Localised ground states
represent energetically preferred stationary configurations,
whereas escape to infinity describes the failure of spatial
inhomogeneity to bind a minimising sequence. The binding--escape
alternative therefore distinguishes localisation from loss of
localisation in a precise variational sense. This interpretation
suggests potential relevance to nonlinear diffusion and
two-component field models.

The logarithmic interaction in \eqref{eq:system} presents two intertwined obstacles. First, the coupling is nonseparable: each logarithmic term is weighted by a power of the opposite component. At \(\varepsilon=0\), the negative logarithmic part may fail to be integrable on the ambient Sobolev space. The energy is therefore extended-valued, and its finite-energy domain need not be closed under subtraction, obstructing the usual procedure of subtracting profiles. Second, asymptotically constant potentials permit translation to infinity, so local compactness does not prevent mass from escaping. The central difficulty is thus to describe this loss of compactness while passing to the singular limit \(\varepsilon\to 0\), where the original energy has no global Fr\'echet derivative and estimates for fixed positive regularisation may deteriorate. A standard Palais--Smale argument for a fixed \(C^1\) functional does not resolve these simultaneous difficulties. The system \eqref{eq:system} and the initial exponent-calibration argument were introduced in \cite{Hu26}. Here we extend the analysis to the noncompact setting and study the potential dependence of the ground-state energy through radial normalisation.

Our analysis builds on the concentration--compactness and global compactness theories discussed above. More recent developments include profile and integral-functional decompositions in Sobolev spaces~\cite{Oku22} and min-max analysis for parametrised elliptic problems on unbounded domains~\cite{DS22}; a fundamental quasilinear counterpart is provided by Mercuri and Willem~\cite{MW10}. These results do not apply directly to the present singular limit: the endpoint energy is extended-valued, its finite-energy domain need not be closed under subtraction, and the regularisation parameter varies along the sequence. Moreover, almost-everywhere convergence alone does not justify splitting the logarithmic terms. Br\'ezis--Lieb splitting~\cite{BL83} and its extensions to superposition operators~\cite{Ack18} provide a useful point of comparison. Here we formulate the endpoint equation through a distributional residual and establish localisation estimates uniform in the regularisation parameter. The Nehari constraint~\cite{Neh60} and Ekeland's principle~\cite{Eke74}
allow us to approximate arbitrary Nehari minimising sequences by states
with small residuals, uniformly in the regularization parameter.
We then determine when these sequences converge to a ground state and
when they escape to infinity, and how escape affects the directional
derivatives of the ground-state energy at the threshold.

We now describe the main results, which fall into three connected parts: global compactness at fixed and vanishing regularisation, and its applications to ground-state convergence and the characterisation of compactness thresholds.

\begin{theorem}[Singular global compactness]
\label{thm:singular-global-decomposition}
Assume (C) holds. Let \(0\le\varepsilon_j\le1/2\),
\(\varepsilon_j\to0\), and let \(z_j\in\mathcal X\) be admissible for
\(\mathfrak R_{q,\varepsilon_j}\). Suppose that
\begin{equation}\label{eq:singular-ps-hypothesis}
 J_{q,\varepsilon_j}(z_j)\longrightarrow c\in\mathbb R,
 \qquad
 \|\mathfrak R_{q,\varepsilon_j}(z_j)\|_{\mathcal X^*}
 \longrightarrow0.
\end{equation}
Then, after passage to a subsequence, there exist an integer \(k\ge0\),
a possibly zero finite-energy solution \(z^0\) for \(q\), nonzero
finite-energy solutions \(Z^1,\ldots,Z^k\) for \(\bar q\), and centres
\(y_j^1,\ldots,y_j^k\in\mathbb R^N\) satisfying
\[
 |y_j^i|\longrightarrow\infty
 \quad (1\le i\le k),
 \qquad
 |y_j^i-y_j^h|\longrightarrow\infty
 \quad (1\le i<h\le k),
\]
such that
\begin{equation}\label{eq:singular-profile-decomposition}
 z_j=z^0+\sum_{i=1}^k\mathsf U_{y_j^i}Z^i+e_j,
 \qquad
 \|e_j\|_q\longrightarrow0.
\end{equation}
Set \(Z^0=z^0\) and \(y_j^0=0\). Then
\begin{align}
 S_q(z_j)
 &=S_q(z^0)+\sum_{i=1}^k S_{\bar q}(Z^i)+o(1),
 \label{eq:singular-norm-splitting}\\
 c
 &=J_q(z^0)+\sum_{i=1}^k J_{\bar q}(Z^i),
 \label{eq:singular-energy-splitting}\\
 \left\|\ell_{\varepsilon_j}^{\pm}(z_j)
 -\sum_{i=0}^k\mathsf U_{y_j^i}\ell_0^{\pm}(Z^i)\right\|_1
 &\longrightarrow0,
 \label{eq:singular-entropy-splitting}\\
 \left\|r_{\varepsilon_j}(z_j)
 -\sum_{i=0}^k\mathsf U_{y_j^i}m(Z^i)\right\|_1
 +\left\|m(z_j)
 -\sum_{i=0}^k\mathsf U_{y_j^i}m(Z^i)\right\|_1
 &\longrightarrow0.
 \label{eq:singular-mass-splitting}
\end{align}
The nonlinear sources, interpreted as continuous extensions of their
distributions, satisfy
\begin{equation}\label{eq:singular-source-splitting}
 \left\|D_\xi\ell_{\varepsilon_j}(z_j)
 -\sum_{i=0}^k\mathsf U_{y_j^i}D_\xi\ell_0(Z^i)
 \right\|_{\mathcal X^*}
 \longrightarrow0.
\end{equation}
Every nonzero limiting state belongs to its unregularised Nehari set.
In particular,
\begin{equation}\label{eq:singular-profile-count}
 c\ge0,
 \qquad
 k d_{\bar q,0}\le c-J_q(z^0),
 \qquad
 J_q(z^0)\ge d_{q,0}\quad\text{if }z^0\ne0.
\end{equation}
\end{theorem}

\begin{remark}
The main feature of this theorem is that profile decomposition remains
valid as the regularisation vanishes, including for sequences with
\(\varepsilon_j=0\). At this singular endpoint, the energy is
extended-valued and its finite-energy domain need not be closed under
subtraction, so the usual iterative subtraction of profiles cannot be
applied directly. We overcome this obstruction by scalar localisation
with slowly expanding cutoffs and residual commutator estimates uniform
in \(\varepsilon\). This yields not only strong convergence of the
remainder and energy splitting, but also separate \(L^1\)-splitting of
the positive and negative logarithmic densities and convergence of the
nonlinear sources in \(\mathcal X^*\). The conclusion retains
information that strong Sobolev convergence alone does not provide.
\end{remark}

\begin{theorem}[The binding--escape alternative]
\label{thm:binding-minimisers}
Let \(0\le\varepsilon_j\le1/2\), \(\varepsilon_j\to0\), and
\(\widehat z_j\in\mathcal N_{q,\varepsilon_j}\) satisfy
\begin{equation}\label{eq:binding-minimising-hypothesis}
J_{q,\varepsilon_j}(\widehat z_j)-d_{q,\varepsilon_j}
\longrightarrow0.
\end{equation}
After passage to a subsequence, exactly one of the following
alternatives holds:
\begin{enumerate}
\item \(\widehat z_j\to z\) strongly in \(\mathcal X\), where
\(z\) is a finite-energy solution for \(q\) and \(J_q(z)=d_{q,0}\);
\item \(d_{q,0}=d_{\bar q,0}\), and there exist centres
\(|y_j|\to\infty\) and an autonomous ground state \(Z\) such that
\(\widehat z_j-\mathsf U_{y_j}Z\to0\) strongly in \(\mathcal X\).
\end{enumerate}
Every minimiser of \(J_q\) on \(\mathcal N_{q,0}\) is a
finite-energy solution, although \(J_q\) need not be
differentiable on \(\mathcal X\). In particular,
\(d_{q,0}<d_{\bar q,0}\) implies attainment by a componentwise
nonnegative solution. Moreover,
\begin{equation}\label{eq:binding-compactness-equivalence}
d_{q,0}<d_{\bar q,0}
\quad\Longleftrightarrow\quad
\text{every minimising sequence in \(\mathcal N_{q,0}\)
has a strongly convergent subsequence.}
\end{equation}
\end{theorem}

\begin{remark}
The preceding theorems describe residual Palais--Smale sequences
and their compactness thresholds. The present theorem extends
this analysis to arbitrary Nehari minimising sequences, for which
small residuals are not assumed. It yields a binding--escape
alternative and an exact compactness criterion for strict binding,
while identifying constrained minimisers as solutions despite
the lack of global differentiability.
\end{remark}

\begin{theorem}[Density response including escape]
\label{thm:binding-density-response}
Let \(q\in\Omega_\infty\).
\begin{enumerate}
\item[(i)] For every \(h\in\mathcal P_0\), the right directional
derivative exists in the Hadamard sense and satisfies
\begin{equation}\label{eq:binding-directional-response}
 F'(q;h)
 =\lim_{\substack{t\downarrow0\\ h'\to h}}
   \frac{F(q+th')-F(q)}{t}
 =\min_{\lambda\in\mathcal A(q)}\lambda[h].
\end{equation}

\item[(ii)] The superdifferential of \(F\) at \(q\) is
\begin{equation}\label{eq:binding-superdifferential}
 \partial^+F(q)
 =\overline{\operatorname{co}}^{\,w^*}\mathcal A(q),
\end{equation}
where the closure is taken in \(\sigma(\mathcal P_0^*,\mathcal P_0)\).

\item[(iii)] If \(d_{q,0}<d_{\bar q,0}\), then \(F\) is G\^ateaux
differentiable at \(q\) if and only if \(\mathcal W(q)\) is a
singleton. In that case, writing \(\mathcal W(q)=\{w\}\), \(F\)
is in fact Fr\'echet differentiable at \(q\), with \(DF(q)=\rho(w)\).

If \(d_{q,0}=d_{\bar q,0}\), then \(F\) is not Fr\'echet differentiable
at \(q\). It is G\^ateaux differentiable at \(q\) if and only if
\(\mathcal W(q)=\varnothing\), in which case its derivative is zero.
\end{enumerate}
\end{theorem}

\begin{remark}
The preceding compactness results yield either a ground-state
limit or an escaping profile. The present theorem translates
this alternative into an energy response formula, with escape
contributing the zero functional to the active set. This
additional contribution determines differentiability at the
threshold and, when the ground level is not attained, separates
G\^ateaux from Fr\'echet differentiability.
\end{remark}

The paper is organised as follows. Section~2 presents the
variational framework and preliminary estimates for the continuous
problem. Section~3 establishes the Palais--Smale decomposition
at fixed positive regularisation, while Section~4 develops
the corresponding decomposition as the regularisation vanishes.
Section~5 investigates mass escape and its implications for
compactness, ground states, and energy response. The appendices
provide the uniform passage from Nehari minimisation to small
residuals and the recovery of radial and logarithmic densities.

\section{Variational framework and preliminary estimates}\label{sec:framework}

This section establishes the radial structure and uniform estimates
for the problem on $\mathbb R^N$ with asymptotically constant
potentials. In particular, we prove uniform bounds for
Palais--Smale sequences and Nehari sublevel sets, without
assuming the existence of a ground state.

All functions are real-valued, and \(p=N>4\).
We write \(q=(a,b)\) for the potential pair. Throughout the paper, we assume
\begin{equation}\tag{C}\label{hyp:continuous}
\begin{gathered}
a,b\in L^\infty_{\mathrm{loc}}(\mathbb R^N),
\qquad a,b\ge V_0>0,\\
\lim_{R\to\infty}
\operatorname*{ess\,sup}_{|x|>R}
\bigl(|a(x)-\bar a|+|b(x)-\bar b|\bigr)=0,
\qquad \bar a,\bar b>0.
\end{gathered}
\end{equation}
In particular, \(a,b\in L^\infty(\mathbb R^N)\).
We denote the limiting constant pair by
\(\bar q=(\bar a,\bar b)\) and set
\[
 \mathcal X=W^{1,p}(\mathbb R^N)\times W^{1,p}(\mathbb R^N).
\]
For \(z=(u,v)\in \mathcal X\), define
\begin{equation}\label{eq:principal-energy}
S_q(z)=\int_{\mathbb R^N}
\bigl(|\nabla u|^p+a(x)|u|^p
      +|\nabla v|^p+b(x)|v|^p\bigr)\,dx,
\qquad
\|z\|_q=S_q(z)^{1/p}.
\end{equation}
The norm \(\|\cdot\|_q\) is equivalent to the usual
\(p\)-product Sobolev norm, and its dual norm is denoted by
\(\|\cdot\|_{q,*}\). Equipped with this norm, \(\mathcal X\) is
reflexive and uniformly convex, since it is isometric to a
closed subspace of a product of \(L^p\) spaces, with Euclidean
fibres for the gradients.

\begin{lemma}[Embeddings]\label{lem:embeddings}
Under \eqref{hyp:continuous}, each component space embeds
continuously into \(L^r(\mathbb R^N)\) for \(p\le r<\infty\),
and compactly into \(L^r_{\mathrm{loc}}(\mathbb R^N)\)
for every \(1\le r<\infty\). Moreover,
\begin{equation}\label{eq:common-embedding}
\|u\|_{L^r}+\|v\|_{L^r}
\le C_r\|z\|_q,
\qquad p\le r<\infty.
\end{equation}
\end{lemma}

\begin{proof}
Since \(p=N\), assumption \eqref{hyp:continuous} makes the
weighted component norms equivalent to the usual \(W^{1,N}\)
norm. For \(N<r<\infty\), the Gagliardo--Nirenberg inequality gives
\[
\|u\|_{L^r}
\le C_{N,r}\|\nabla u\|_{L^N}^{1-N/r}
                 \|u\|_{L^N}^{N/r}.
\]
Together with the corresponding estimate for \(v\) and the
positive lower bound on the potentials, this proves
\eqref{eq:common-embedding}; the case \(r=N\) follows directly
from the definition of \(\|\cdot\|_q\).
The local compact embeddings follow from the
Rellich--Kondrachov theorem on bounded balls and a diagonal
subsequence argument.
\end{proof}

\begin{definition}\label{def:logarithmic-densities}
For \(0\le\varepsilon\le1/2\) and \((s,t)\in\mathbb R^2\), define
\begin{align}
\ell_\varepsilon(s,t)
&=|s|^{p-2}t^2\log(t^2+\varepsilon)
 +|t|^{p-2}s^2\log(s^2+\varepsilon),
\label{eq:log-density}\\
m(s,t)
&=|s|^{p-2}t^2+|t|^{p-2}s^2,
\label{eq:mass-density}\\
r_\varepsilon(s,t)
&=|s|^{p-2}\frac{t^4}{t^2+\varepsilon}
 +|t|^{p-2}\frac{s^4}{s^2+\varepsilon},
\qquad r_0=m.
\label{eq:radial-density}
\end{align}
Set
\[
\ell_\varepsilon^\pm(s,t)
=|s|^{p-2}t^2[\log(t^2+\varepsilon)]_\pm
 +|t|^{p-2}s^2[\log(s^2+\varepsilon)]_\pm,
\]
where \([\tau]_\pm=\max\{\pm\tau,0\}\); all expressions at
zero components are defined by continuity. For \(z=(u,v)\in \mathcal X\),
write
\[
L_\varepsilon^\pm(z)=\int_{\mathbb R^N}\ell_\varepsilon^\pm(z)\,dx,
\qquad
M(z)=\int_{\mathbb R^N}m(z)\,dx,
\qquad
R_\varepsilon(z)=\int_{\mathbb R^N}r_\varepsilon(z)\,dx.
\]
\end{definition}

\begin{lemma}[Uniform logarithmic bounds]\label{lem:log-bounds}
For every \(0<\sigma<1\), there exists \(C_\sigma>0\),
independent of \(0\le\varepsilon\le1/2\), such that
\begin{align}
\ell_\varepsilon^+(\xi)
&\le C_\sigma|\xi|^{p+\sigma},
\label{eq:positive-growth}\\
\ell_\varepsilon^+(\xi)+\ell_\varepsilon^-(\xi)
&\le C_\sigma\bigl(|\xi|^{p-\sigma}+|\xi|^{p+\sigma}\bigr),
\label{eq:absolute-growth}\\
|D_\xi\ell_\varepsilon(\xi)|
&\le C_\sigma\bigl(
|\xi|^{p-1-\sigma}+|\xi|^{p-1+\sigma}\bigr).
\label{eq:derivative-growth}
\end{align}
Moreover, \(\ell_\varepsilon\in C^1(\mathbb R^2)\), including
\(\varepsilon=0\), and
\begin{equation}\label{eq:mass-control}
0\le R_\varepsilon(z)\le M(z)
\le\int_{\mathbb R^N}(|u|^p+|v|^p)\,dx
\le V_0^{-1}S_q(z).
\end{equation}
\end{lemma}

\begin{proof}
For \(h_\varepsilon(s)=s^2\log(s^2+\varepsilon)\), elementary
estimates give, uniformly in \(0\le\varepsilon\le1/2\),
\[
[h_\varepsilon(s)]_+\le C_\sigma|s|^{2+\sigma},
\qquad
|h_\varepsilon(s)|
\le C_\sigma\bigl(|s|^{2-\sigma}+|s|^{2+\sigma}\bigr),
\]
and
\[
|h_\varepsilon'(s)|
\le C_\sigma\bigl(|s|^{1-\sigma}+|s|^{1+\sigma}\bigr),
\]
with continuous extensions at zero.
Differentiation and weighted Young inequalities yield
\eqref{eq:positive-growth}--\eqref{eq:derivative-growth};
the same formulas show that \(\ell_\varepsilon\in C^1(\mathbb R^2)\).
Finally, \eqref{eq:mass-control} follows from
\(0\le r_\varepsilon\le m\), Young's inequality, and \(a,b\ge V_0\).
\end{proof}
Set
\[
L_\varepsilon(z)=L_\varepsilon^+(z)-L_\varepsilon^-(z).
\]
For $0<\varepsilon\le1/2$, the energy functional associated
with \eqref{eq:system} is
\begin{equation}
J_{q,\varepsilon}(z)
=\frac1p\bigl(S_q(z)-L_\varepsilon(z)\bigr),
\qquad z\in\mathcal X.
\label{eq:regularised-energy}
\end{equation}
At $\varepsilon=0$, the original system is associated with
the extended-valued energy
\begin{equation}
\begin{gathered}
J_q(z)=J_{q,0}(z)
=\frac1p\bigl(S_q(z)-L_0^+(z)+L_0^-(z)\bigr),
\qquad z\in\mathcal X,\\
\mathcal D_q
=\{z\in\mathcal X:L_0^-(z)<\infty\},
\end{gathered}
\label{eq:effective-domain}
\end{equation}
where \(\mathcal D_q\) is the domain of the oringinal energy functional.
\begin{proposition}[Regularity and radial identity]
\label{prop:energy-regularity}
Assume \eqref{hyp:continuous}. For every
\(0<\varepsilon\le1/2\),
\(J_{q,\varepsilon}\in C^1(\mathcal X,\mathbb R)\), and
\begin{equation}
 G_{q,\varepsilon}(z)
 :=\langle J_{q,\varepsilon}'(z),z\rangle
 =S_q(z)-L_\varepsilon(z)-\frac2pR_\varepsilon(z).
 \label{eq:nehari-functional}
\end{equation}
Moreover,
\(G_{q,\varepsilon}\in C^1(\mathcal X,\mathbb R)\), and
\begin{equation}
 J_{q,\varepsilon}(z)-\frac1pG_{q,\varepsilon}(z)
 =\frac2{p^2}R_\varepsilon(z).
 \label{eq:radial-energy-identity}
\end{equation}
The extended-valued functional \(J_q\) is strongly lower
semicontinuous on \(\mathcal X\), and
\begin{equation}
 J_{q,\varepsilon}(z)\uparrow J_q(z)
 \qquad\text{as }\varepsilon\downarrow0
 \quad\text{for every }z\in \mathcal X.
 \label{eq:monotone-regularisation}
\end{equation}
\end{proposition}

\begin{proof}
For fixed \(\varepsilon>0\) and \(0<\sigma<1\),
\[
 |D_\xi\ell_\varepsilon(\xi)|
 \le C_{\varepsilon,\sigma}
 \bigl(|\xi|^{p-1}+|\xi|^{p-1+\sigma}\bigr).
\]
The embeddings in Lemma~\ref{lem:embeddings} and the standard
differentiability theorem for integral functionals
\cite{Wil96} therefore yield
\(J_{q,\varepsilon}\in C^1(\mathcal X,\mathbb R)\).
Likewise, \(R_\varepsilon\in C^1(\mathcal X,\mathbb R)\), since
\[
 |D_\xi r_\varepsilon(\xi)|\le C|\xi|^{p-1}.
\]
The radial identities
\eqref{eq:nehari-functional}--\eqref{eq:radial-energy-identity}
were established in \cite{Hu26}; in particular,
\(G_{q,\varepsilon}\in C^1(\mathcal X,\mathbb R)\).

If \(z_j\to z\) in \(\mathcal X\), then
\(S_q(z_j)\to S_q(z)\) and
\(L_0^+(z_j)\to L_0^+(z)\), whereas Fatou's lemma gives
\[
 L_0^-(z)\le\liminf_{j\to\infty}L_0^-(z_j).
\]
Thus \(J_q=(S_q-L_0^++L_0^-)/p\) is strongly lower
semicontinuous.

Finally, for fixed \(0<\varepsilon_0\le1/2\), the densities
\(\ell_{\varepsilon_0}-\ell_\varepsilon\) increase nonnegatively
to \(\ell_{\varepsilon_0}-\ell_0\) as
\(\varepsilon\downarrow0\).
The monotone convergence theorem proves
\eqref{eq:monotone-regularisation}, including when
\(J_q(z)=+\infty\).
\end{proof}

\begin{definition}[Nehari set and level]
Following \cite{Hu26}, for \(\varepsilon=0\) define
\[
 G_{q,0}(z)=S_q(z)-L_0(z)-\frac2pM(z),
 \qquad z\in\mathcal D_q.
\]
This is the radial derivative
\(\left.\frac{d}{dt}J_q(tz)\right|_{t=1}\).
Indeed, \(\mathcal D_q\) is invariant under positive scaling, since
\[
 L_0^-(tz)
 \le t^p\bigl(L_0^-(z)+2|\log t|M(z)\bigr),
 \qquad t>0.
\]
Set \(\mathcal D_{q,\varepsilon}=\mathcal X\) for \(\varepsilon>0\)
and \(\mathcal D_{q,0}=\mathcal D_q\). Define
\begin{equation}
 \mathcal N_{q,\varepsilon}
 =\{z\in\mathcal D_{q,\varepsilon}:
       M(z)>0,\ G_{q,\varepsilon}(z)=0\},
 \qquad
 d_{q,\varepsilon}
 =\inf_{z\in\mathcal N_{q,\varepsilon}}J_{q,\varepsilon}(z).
 \label{eq:nehari-level}
\end{equation}
We also write \(d(q)=d_{q,0}\), without assuming attainment.

For \(\varepsilon>0\), put
\[
 Q_\varepsilon(z)
 =\int_{\mathbb R^N}
 \left(
 \varepsilon|u|^{p-2}\frac{v^4}{(v^2+\varepsilon)^2}
 +\varepsilon|v|^{p-2}\frac{u^4}{(u^2+\varepsilon)^2}
 \right)\,dx,
 \qquad Q_0(z)=0.
\]
Then \(0\le Q_\varepsilon(z)\le R_\varepsilon(z)\).
\end{definition}

\begin{proposition}[Radial projection]
\label{prop:radial-projection}
Let \(0\le\varepsilon\le1/2\) and
\(z\in\mathcal D_{q,\varepsilon}\) with \(M(z)>0\).
There exists a unique \(t_{q,\varepsilon}(z)>0\) such that
\(t_{q,\varepsilon}(z)z\in\mathcal N_{q,\varepsilon}\).
This scale is the unique global maximiser of
\(t\mapsto J_{q,\varepsilon}(tz)\) on \((0,\infty)\).
For \(\varepsilon=0\),
\begin{align}
 \log t_{q,0}(z)
 &=\frac{S_q(z)-L_0(z)}{2M(z)}-\frac1p,
 \label{eq:exact-ray-scale}\\
 \max_{t>0}J_q(tz)
 &=\frac{2M(z)}{p^2}
   \exp\!\left(
   \frac{p(S_q(z)-L_0(z))}{2M(z)}-1
   \right).
 \label{eq:exact-ray-energy}
\end{align}
For \(\varepsilon>0\), the set
\(\mathcal N_{q,\varepsilon}\) is a \(C^1\) hypersurface in
\(\mathcal X\), and
\begin{equation}
 DG_{q,\varepsilon}(z)[z]
 =-2R_\varepsilon(z)-\frac4pQ_\varepsilon(z)<0,
 \qquad z\in\mathcal N_{q,\varepsilon}.
 \label{eq:radial-transversality}
\end{equation}
Every constrained critical point of \(J_{q,\varepsilon}\)
on \(\mathcal N_{q,\varepsilon}\) is a critical point in \(\mathcal X\).
\end{proposition}

\begin{proof}
For \(\varepsilon>0\), the radial projection and maximum
properties follow by the same argument as in Lemma~3.5
of \cite{Hu26}, which does not require compact embeddings.

For \(\varepsilon=0\), the scaling identity
\begin{equation}
 L_0(tz)=t^p\bigl(L_0(z)+2M(z)\log t\bigr)
 \label{eq:log-scaling}
\end{equation}
gives
\[
 J_q(tz)=\frac{t^p}{p}
 \bigl(S_q(z)-L_0(z)-2M(z)\log t\bigr).
\]
Differentiation yields
\[
 \frac{d}{dt}J_q(tz)
 =t^{p-1}
 \left(
 S_q(z)-L_0(z)-2M(z)\log t-\frac2pM(z)
 \right).
\]
The expression in parentheses decreases strictly from
\(+\infty\) to \(-\infty\), proving uniqueness and the maximum
property. Solving for its zero and evaluating the energy give
\eqref{eq:exact-ray-scale}--\eqref{eq:exact-ray-energy}.

For \(\varepsilon>0\), direct differentiation at
\(z\in\mathcal N_{q,\varepsilon}\) gives
\[
 DG_{q,\varepsilon}(z)[z]
 =-2R_\varepsilon(z)-\frac4pQ_\varepsilon(z)<0.
\]
The implicit function theorem proves the hypersurface assertion.
Finally, if
\(J_{q,\varepsilon}'(z)=\lambda DG_{q,\varepsilon}(z)\),
pairing with \(z\) gives \(\lambda=0\).
\end{proof}

\begin{corollary}[Uniform positive levels]
\label{cor:uniform-levels}
For fixed \(q\), there exist constants \(d_*,d^*>0\),
independent of \(0\le\varepsilon\le1/2\), such that
\begin{equation}
 0<d_*\le d_{q,\varepsilon}\le d^*<\infty.
 \label{eq:uniform-levels}
\end{equation}
Moreover,
\begin{equation}
 S_q(z)\ge C(V_0,p,d_*)>0,
 \qquad z\in\mathcal N_{q,\varepsilon}.
 \label{eq:nehari-gaps}
\end{equation}
The infimum defining \(d_{q,\varepsilon}\) is unchanged
if restricted to componentwise nonnegative pairs.
\end{corollary}

\begin{proof}
The energy bounds follow by the same small-sphere and
fixed-ray comparison arguments as in \cite{Hu26}.
Indeed, \eqref{eq:positive-growth} and
\eqref{eq:common-embedding} give
\[
 J_{q,\varepsilon}(z)
 \ge \frac1p\|z\|_q^p-C\|z\|_q^{p+\sigma},
\]
uniformly in \(\varepsilon\), while a fixed nonnegative pair
\(\zeta=(\chi,\chi)\), with
\(0\ne\chi\in C_c^\infty(\mathbb R^N)\), gives
\[
 d_{q,\varepsilon}
 \le\max_{t>0}J_{q,\varepsilon}(t\zeta)
 \le\max_{t>0}J_q(t\zeta)<\infty.
\]
Proposition~\ref{prop:radial-projection} therefore yields
\eqref{eq:uniform-levels}.
On \(\mathcal N_{q,\varepsilon}\), the radial identity and
\eqref{eq:mass-control} imply
\[
 d_*\le J_{q,\varepsilon}(z)
 =\frac2{p^2}R_\varepsilon(z)
 \le\frac{2}{p^2V_0}S_q(z),
\]
which proves \eqref{eq:nehari-gaps}.
Finally, taking componentwise absolute values preserves the
nonlinear terms and does not increase the principal energy;
comparison of radial maxima proves the last assertion.
\end{proof}

The following lemma gives a localised version of the
exponent-calibration argument developed in \cite{Hu26}.
The estimate is uniform in \(\varepsilon\) and will be used
in the concentration analysis and singular profile decomposition.

\begin{lemma}[Calibrated positive-part estimate]
\label{lem:calibration}
Let \(0<\theta<1\), \(\sigma>0\), and
\[
 r=p+\frac{\sigma}{1-\theta}.
\]
For every measurable \(A\subset\mathbb R^N\), every
\(z\in \mathcal X\), and every \(0\le\varepsilon\le1/2\),
\begin{equation}
 \int_A\ell_\varepsilon^+(z)\,dx
 \le C_{\theta,\sigma}
 \left(\int_A r_\varepsilon(z)\,dx\right)^\theta
 S_q(z)^{\,1-\theta+\sigma/p}.
 \label{eq:calibration}
\end{equation}
The constant is independent of \(A\) and \(\varepsilon\).
\end{lemma}

\begin{proof}
Consider the first logarithmic term and set
\(A_v=A\cap\{v^2+\varepsilon\ge1\}\).
On \(A_v\), \(v^2\ge1/2\ge\varepsilon\), so
\[
 |u|^{p-2}v^2
 \le2|u|^{p-2}\frac{v^4}{v^2+\varepsilon},
 \qquad
 \log(v^2+\varepsilon)\le C_\sigma|v|^\sigma.
\]
H\"older's inequality gives
\begin{align*}
 &\int_{A_v}|u|^{p-2}v^2
       \log(v^2+\varepsilon)\,dx\\
 &\quad\le C_\sigma
 \left(\int_{A_v}|u|^{p-2}v^2\,dx\right)^\theta
 \left(
 \int_{A_v}|u|^{p-2}
 |v|^{2+\sigma/(1-\theta)}\,dx
 \right)^{1-\theta}\\
 &\quad\le C_{\theta,\sigma}
 \left(\int_A r_\varepsilon(z)\,dx\right)^\theta
 \|u\|_{L^r}^{(p-2)(1-\theta)}
 \|v\|_{L^r}^{2(1-\theta)+\sigma}.
\end{align*}
The second application of H\"older's inequality is valid because
\[
 (p-2)+2+\frac{\sigma}{1-\theta}=r.
\]
By \eqref{eq:common-embedding}, the product of norms is bounded by
\[
 C_{\theta,\sigma}\|z\|_q^{p(1-\theta)+\sigma}
 =C_{\theta,\sigma}S_q(z)^{1-\theta+\sigma/p}.
\]
Interchanging \(u\) and \(v\) gives the same bound for the
second logarithmic term. Adding the two estimates proves
\eqref{eq:calibration}.
\end{proof}

\begin{lemma}[Uniform bounds without compact embeddings]
\label{thm:uniform-ps-bound}
Assume \eqref{hyp:continuous}.
For every \(E_0,K_1\ge0\), there exists \(C>0\),
independent of \(0<\varepsilon\le1/2\), such that every
\(z\in \mathcal X\) satisfying
\[
 |J_{q,\varepsilon}(z)|\le E_0,
 \qquad
 \|J_{q,\varepsilon}'(z)\|_{q,*}\le K_1
\]
obeys
\begin{equation}
 \|z\|_q\le C,\qquad
 L_\varepsilon^+(z)\le C,\qquad
 L_\varepsilon^-(z)\le C.
 \label{eq:uniform-ps-bound}
\end{equation}
In particular, these bounds apply to regularised
Palais--Smale sequences uniformly in the regularisation
parameter.

For every \(E_0\ge0\), the same bounds hold for all
\(0\le\varepsilon\le1/2\) and
\(z\in\mathcal N_{q,\varepsilon}\) satisfying
\(J_{q,\varepsilon}(z)\le E_0\), without any assumption
on the full derivative.
\end{lemma}

\begin{proof}
Write \(Y=S_q(z)\) and
\(\eta=\|J_{q,\varepsilon}'(z)\|_{q,*}\le K_1\).
Since
\[
 |G_{q,\varepsilon}(z)|\le\eta Y^{1/p},
\]
the radial energy identity gives
\begin{equation}
 R_\varepsilon(z)
 \le\frac{p^2}{2}E_0+\frac p2\eta Y^{1/p}
 \le C(1+Y^{1/p}).
 \label{eq:radial-residual-control}
\end{equation}
Choose \(0<\theta<1\) and
\(0<\sigma<\min\{1,\theta(p-1)\}\), and set
\[
 \beta=1-\theta+\frac{\sigma}{p},
 \qquad
 \gamma=\beta+\frac{\theta}{p}<1.
\]
By \eqref{eq:nehari-functional}, \eqref{eq:calibration},
and \eqref{eq:radial-residual-control},
\begin{align}
 Y
 &\le L_\varepsilon^+(z)
      +\frac2pR_\varepsilon(z)
      +|G_{q,\varepsilon}(z)|\notag\\
 &\le C R_\varepsilon(z)^\theta Y^\beta
      +C(1+Y^{1/p})\notag\\
 &\le C\bigl(Y^\beta+Y^\gamma+1+Y^{1/p}\bigr).
 \label{eq:sublinear-absorption}
\end{align}
All powers on the last line are strictly less than one,
so \(Y\), and hence \(\|z\|_q\), is uniformly bounded.
The positive-part estimate and the continuous embeddings
then give \(L_\varepsilon^+(z)\le C\).
Finally,
\[
 L_\varepsilon^-(z)
 =pJ_{q,\varepsilon}(z)-S_q(z)+L_\varepsilon^+(z)
 \le pE_0+L_\varepsilon^+(z)\le C.
\]

For the Nehari sublevels, including \(\varepsilon=0\),
the radial identity gives
\[
 G_{q,\varepsilon}(z)=0,
 \qquad
 R_\varepsilon(z)=\frac{p^2}{2}J_{q,\varepsilon}(z)
 \le\frac{p^2}{2}E_0.
\]
Applying \eqref{eq:calibration} with
\(0<\sigma<\min\{1,p\theta\}\), we obtain
\[
 Y\le L_\varepsilon^+(z)+\frac2pR_\varepsilon(z)
 \le CY^{1-\theta+\sigma/p}+C.
\]
Since \(1-\theta+\sigma/p<1\), this bounds \(Y\).
The bounds for \(L_\varepsilon^\pm(z)\) follow as above.
\end{proof}

\begin{remark}\label{rem:ps-not-cerami}
The proof of Lemma~\ref{thm:uniform-ps-bound} does not assume
 $\langle J_{q,\varepsilon_j}'(z_j),z_j\rangle=o(1)$
before boundedness is established. Instead, the residual term
in \eqref{eq:radial-residual-control} is retained and absorbed.
Hence the lemma applies to Palais--Smale sequences without
requiring the stronger Cerami condition.
\end{remark}

\begin{lemma}[Local convergence of the nonlinearities]
\label{lem:local-nonlinear-convergence}
Suppose that \(0<\varepsilon_j\le1/2\),
\(\varepsilon_j\to0\), and \(z_j\rightharpoonup z\) in \(\mathcal X\).
After passage to a subsequence, for every compact
\(K\subset\mathbb R^N\),
\begin{align}
 D_\xi\ell_{\varepsilon_j}(z_j)
 &\longrightarrow D_\xi\ell_0(z)
 &&\text{in }L^{p'}(K;\mathbb R^2),
 \label{eq:local-source-convergence}\\
 \ell_{\varepsilon_j}^\pm(z_j)
 &\longrightarrow\ell_0^\pm(z)
 &&\text{in }L^1(K),
 \label{eq:local-log-convergence}\\
 r_{\varepsilon_j}(z_j)
 &\longrightarrow m(z)
 &&\text{in }L^1(K).
 \label{eq:local-radial-convergence}
\end{align}
\end{lemma}

\begin{proof}
These conclusions follow by the local compactness and
uniform integrability arguments used in the proof of
Theorem~1.8 of \cite{Hu26}. We omit the details.
\end{proof}

\section{Global compactness at fixed regularization}

Throughout this section, we assume \eqref{hyp:continuous}
and fix \(0<\varepsilon\le1/2\).
We retain the notation of Section~\ref{sec:framework} and write
\[
 \bar q=(\bar a,\bar b).
\]
We equip \(\mathcal X\) with the norm
\(\|z\|_{\bar q}=S_{\bar q}(z)^{1/p}\)
and denote its dual norm by \(\|\cdot\|_{\mathcal X^*}\).
Since the potentials are bounded and uniformly positive,
this norm is equivalent to \(\|\cdot\|_q\).

For \(y\in\mathbb R^N\), define
\[
 (\mathsf U_yz)(x)=z(x-y).
\]
Translations are isometries of
\((\mathcal X,\|\cdot\|_{\bar q})\).
A Palais--Smale sequence at level \(c\) is a sequence
\((z_j)\subset\mathcal X\) satisfying
\begin{equation}
 J_{q,\varepsilon}(z_j)\longrightarrow c,
 \qquad
 \|J_{q,\varepsilon}'(z_j)\|_{\mathcal X^*}
 \longrightarrow0.
 \label{eq:fixed-ps-definition}
\end{equation}
Lemma~\ref{thm:uniform-ps-bound} supplies boundedness directly
from \eqref{eq:fixed-ps-definition}, with the order of the
estimates explained in Remark~\ref{rem:ps-not-cerami}.
The principal result is Theorem 6.1.
It gives strong convergence of the final remainder and exact
splitting of the energy and both logarithmic parts. Splitting
of the derivative is proved before each profile is subtracted.
The reverse construction identifies all levels at which
compactness fails.
\subsection{Splitting and local recovery}
\label{subsec:fixed-splitting}

In what follows, $L^{p'}+L^{r'}$ denotes the sum of the two
Lebesgue spaces, with norm
\[
 \|f\|_{L^{p'}+L^{r'}}
 =\inf_{f=f_1+f_2}
       \bigl(\|f_1\|_{p'}+\|f_2\|_{r'}\bigr).
\]
Vector-valued Lebesgue spaces use the Euclidean norm.
As in Section~\ref{sec:framework}, product Sobolev spaces
carry the $p$-product norm. We use $p'=p/(p-1)$ and
$r'=r/(r-1)$. By
Lemma~\ref{lem:embeddings}, this space acts continuously on
$\mathcal X$ whenever $p<r<\infty$.

\begin{lemma}[Splitting at a positive regularisation parameter]
\label{lem:fixed-nonlinear-splitting}
Fix \(0<\varepsilon\le1/2\) and \(0<\sigma<1\), and put
\(r=p+\sigma\).
Suppose that \(\zeta_j\to\zeta\) almost everywhere in
\(\mathbb R^N\) and that \((\zeta_j)\) is bounded in
\[
 L^p(\mathbb R^N;\mathbb R^2)
 \cap L^r(\mathbb R^N;\mathbb R^2).
\]
Set \(\omega_j=\zeta_j-\zeta\).
For every
\begin{equation}
 H\in
 \{\ell_\varepsilon,\ell_\varepsilon^+,
   \ell_\varepsilon^-,m,r_\varepsilon\},
 \label{eq:fixed-density-family}
\end{equation}
one has
\begin{equation}
 \|H(\zeta_j)-H(\omega_j)-H(\zeta)\|_{L^1}
 \longrightarrow0.
 \label{eq:fixed-density-splitting}
\end{equation}
Moreover,
\begin{equation}
 \bigl\|
 D_\xi\ell_\varepsilon(\zeta_j)
 -D_\xi\ell_\varepsilon(\omega_j)
 -D_\xi\ell_\varepsilon(\zeta)
 \bigr\|_{L^{p'}+L^{r'}}
 \longrightarrow0.
 \label{eq:fixed-source-splitting}
\end{equation}
If $U_j\to U$ almost everywhere and $(U_j)$ is bounded in
$L^p$, then
\begin{align}
 \big\||U_j|^p-|U_j-U|^p-|U|^p\big\|_1&\longrightarrow0,
 \label{eq:fixed-power-splitting}\\
 \big\||U_j|^{p-2}U_j-|U_j-U|^{p-2}(U_j-U)
                     -|U|^{p-2}U\big\|_{p'}&\longrightarrow0.
 \label{eq:fixed-flux-splitting}
\end{align}
\end{lemma}
\begin{proof}
   Since $p>4$ and $\varepsilon>0$, the density $\ell_\varepsilon$ is $C^2$.
Together with the formula for $h_\varepsilon'$ in the proof of
Lemma~\ref{lem:log-bounds}, we use
\[
 h_\varepsilon''(s)
 =2\log(s^2+\varepsilon)+\frac{10s^2}{s^2+\varepsilon}
                         -\frac{4s^4}{(s^2+\varepsilon)^2}.
\]
The second derivatives are
\begin{align*}
 \partial_{ss}\ell_\varepsilon(s,t)
  &=(p-2)(p-3)|s|^{p-4}h_\varepsilon(t)
                         +|t|^{p-2}h_\varepsilon''(s),\\
 \partial_{st}\ell_\varepsilon(s,t)
  &=(p-2)|s|^{p-4}s\,h_\varepsilon'(t)
                         +(p-2)|t|^{p-4}t\,h_\varepsilon'(s),
\end{align*}
with the symmetric expression for $\partial_{tt}\ell_\varepsilon$.
The powers at zero are well defined because $p>4$.
The bounds on $h_\varepsilon,h_\varepsilon',h_\varepsilon''$ and weighted Young
inequalities give
\begin{equation}
 |D_\xi^2\ell_\varepsilon(\xi)|
       \le C_{\varepsilon,\sigma}(|\xi|^{p-2}+|\xi|^{r-2}).
 \label{eq:fixed-hessian-growth}
\end{equation}
Each density in \eqref{eq:fixed-density-family} is locally
Lipschitz, vanishes at zero, and satisfies
$|DH(\xi)|\le C_{\varepsilon,\sigma}
(|\xi|^{p-1}+|\xi|^{r-1})$ wherever its derivative exists. Consequently, for every
$\eta>0$, the mean value estimate and Young's inequality yield
\begin{align}
 |H(\alpha+\beta)-H(\alpha)|
 &\le\eta(|\alpha|^p+|\alpha|^r)
      +C_\eta(|\beta|^p+|\beta|^r),
 \label{eq:fixed-value-increment}\\
 |D_\xi\ell_\varepsilon(\alpha+\beta)
   -D_\xi\ell_\varepsilon(\alpha)|
 &\le\eta(|\alpha|^{p-1}+|\alpha|^{r-1})
      +C_\eta(|\beta|^{p-1}+|\beta|^{r-1}).
 \label{eq:fixed-source-increment}
\end{align}
Here $C_\eta$ may also depend on $\varepsilon,p,\sigma$. Let $K_j=H(\zeta_j)-H(\omega_j)-H(\zeta)$.
Then $K_j\to0$ almost everywhere. By
\eqref{eq:fixed-value-increment},
\[
 \bigl(|K_j|-\eta(|\omega_j|^p+|\omega_j|^r)\bigr)_+
       \le C_\eta'(|\zeta|^p+|\zeta|^r).
\]
The left-hand side tends to zero almost everywhere, so dominated
convergence and boundedness of $(\omega_j)$ give
$\limsup_j\|K_j\|_1\le C\eta$. Letting $\eta\downarrow0$
proves \eqref{eq:fixed-density-splitting}.
For the derivative, apply the same argument to the norm of
the vector defect in
\eqref{eq:fixed-source-splitting}, using
\eqref{eq:fixed-source-increment}. After subtraction of
$\eta(|\omega_j|^{p-1}+|\omega_j|^{r-1})$, the positive
remainder is bounded by
$C_\eta'(|\zeta|^{p-1}+|\zeta|^{r-1})$.
On $\{|\zeta|\le1\}$ this bound lies in $L^{p'}$, and on
$\{|\zeta|>1\}$ it lies in $L^{r'}$.
Dominated convergence on these two sets, followed by
$\eta\downarrow0$, proves \eqref{eq:fixed-source-splitting}.
Applying the identical argument to $|U|^p$ and
$|U|^{p-2}U$ gives
\eqref{eq:fixed-power-splitting}--\eqref{eq:fixed-flux-splitting}.
For the continuity assertion, the mean value bounds before
Young's inequality give local Lipschitz estimates for the
superposition maps by H\"older's inequality. The Lipschitz
constants are bounded on bounded subsets of the stated
domain spaces.
\end{proof}

\begin{lemma}[Local strong convergence and the potential at infinity]
\label{lem:fixed-local-recovery}
The following statements hold.
\begin{enumerate}
\item Let \(q_j=(a_j,b_j)\) satisfy
\[
 0<\underline V\le a_j,b_j\le\overline V<\infty
 \quad\text{almost everywhere, uniformly in }j,
\]
and suppose that
\(q_j\to q_*\) in \(L^\infty_{\mathrm{loc}}(\mathbb R^N)^2\).
If
\[
 z_j\rightharpoonup z
 \quad\text{in }\mathcal X,
 \qquad
 \|J_{q_j,\varepsilon}'(z_j)\|_{\mathcal X^*}
 \longrightarrow0,
\]
then
\[
 z_j\to z
 \quad\text{in }W^{1,p}_{\mathrm{loc}}(\mathbb R^N)^2,
 \qquad
 J_{q_*,\varepsilon}'(z)=0
 \quad\text{in }\mathcal X^*.
\]

\item If \(w_j\rightharpoonup0\) in \(\mathcal X\), then
\begin{align}
 S_q(w_j)-S_{\bar q}(w_j)
 &\longrightarrow0,
 \label{eq:compact-potential-energy}\\
 \|J_{q,\varepsilon}'(w_j)
   -J_{\bar q,\varepsilon}'(w_j)\|_{\mathcal X^*}
 &\longrightarrow0.
 \label{eq:compact-potential-derivative}
\end{align}
In particular, these conclusions apply to
\(w_j=\mathsf U_{y_j}W\) whenever \(W\in\mathcal X\)
and \(|y_j|\to\infty\).
\end{enumerate}
\end{lemma}
\begin{proof}
    Write $z_j=(u_j,v_j)$ and $z=(u,v)$.
The local compact embedding in Lemma~\ref{lem:embeddings}
gives strong convergence in every finite local Lebesgue space.
The fixed-$\varepsilon$ growth estimate in
Proposition~\ref{prop:energy-regularity} consequently gives
\[
 D_\xi\ell_\varepsilon(z_j)\longrightarrow D_\xi\ell_\varepsilon(z)
                 \quad\hbox{in }L^{p'}_{\mathrm{loc}}.
\]
This follows, for example, by choosing a local Lebesgue exponent
larger than $(p-1+\sigma)p'$ and applying uniform integrability.
Let $\chi\in C_c^\infty(\mathbb R^N)$ satisfy $0\le\chi\le1$.
Testing the derivative equation with $\chi(z_j-z)$, the potential
and source terms tend to zero. The terms containing $\nabla\chi$
also tend to zero, by the boundedness of the gradient fluxes
in $L^{p'}$ and local strong $L^p$ convergence. Subtracting
the pairings with the weakly convergent gradients yields
\begin{align*}
 &\int\chi
 \bigl(|\nabla u_j|^{p-2}\nabla u_j-|\nabla u|^{p-2}\nabla u\bigr)
                     \cdot(\nabla u_j-\nabla u)\,dx\\
 &\quad+
 \int\chi
 \bigl(|\nabla v_j|^{p-2}\nabla v_j-|\nabla v|^{p-2}\nabla v\bigr)
                     \cdot(\nabla v_j-\nabla v)\,dx=o(1).
\end{align*}
The inequality
\[
 (|A|^{p-2}A-|B|^{p-2}B)\cdot(A-B)
                        \ge2^{2-p}|A-B|^p
\]
proves the asserted local strong convergence. We can now pass
to the limit in the first-variation formula against compactly
supported smooth tests. The resulting derivative is continuous
on $\mathcal X$, and density extends the identity to every test
in $\mathcal X$. In particular, gradient convergence precedes identification
of the limiting $p$-Laplacians.

For the second assertion, local strong $L^p$ convergence and
\eqref{hyp:continuous} give, componentwise,
\begin{align*}
 \int |a-\bar a|\,|(w_j)_1|^p\,dx&\longrightarrow0,\\
 \|(a-\bar a)|(w_j)_1|^{p-2}(w_j)_1\|_{p'}&\longrightarrow0.
\end{align*}
Indeed, split each integral into $B_R$ and its complement;
first let $j\to\infty$, then $R\to\infty$.
The same estimates hold for the second component.
They imply \eqref{eq:compact-potential-energy} and
\eqref{eq:compact-potential-derivative}.
Finally, translates of a fixed $W\in\mathcal X$ converge
strongly to zero in $W^{1,p}$ on bounded sets when their centres
escape to infinity. Their boundedness implies weak convergence
to zero in $\mathcal X$, proving the final assertion.
\end{proof}

\subsection{Uniform vanishing and the energy of a profile}
\label{subsec:uniform-vanishing}
The next estimate is uniform as the regularization parameter
tends to zero. The equation residual controls the global
\(L^p\) norm, which local vanishing alone does not control.

\begin{lemma}[Vanishing with a small equation residual]
\label{lem:uniform-ps-vanishing}
Let \(0<\varepsilon_j\le1/2\), and let
\(q_j=(a_j,b_j)\) satisfy
\[
 0<\underline V\le a_j,b_j\le\overline V<\infty
 \quad\text{almost everywhere, uniformly in }j.
\]
Suppose that \((z_j)\) is bounded in \(\mathcal X\) and
\begin{equation}
 \|J_{q_j,\varepsilon_j}'(z_j)\|_{\mathcal X^*}
 \longrightarrow0,
 \qquad
 \sup_{y\in\mathbb R^N}
 \int_{B_1(y)}|z_j|^p\,dx\longrightarrow0.
 \label{eq:uniform-vanishing-hypotheses}
\end{equation}
Then
\begin{equation}
 \|z_j\|_{\bar q}\longrightarrow0,
 \qquad
 J_{q_j,\varepsilon_j}(z_j)\longrightarrow0.
 \label{eq:uniform-vanishing-conclusion}
\end{equation}
\end{lemma}
\begin{proof}
   Firstly, we partition $\mathbb{R}^N$ into unit cubes $Q$.
The local embedding $W^{1,p}(Q)\hookrightarrow L^{2p}(Q)$,
with a constant independent of the cube, gives
\[
 \int_Q|z_j|^{3p/2} dx
 \le\left(\int_Q|z_j|^p dx\right)^{1/2}
       \|z_j\|_{L^{2p}(Q)}^p
 \le C\left(\int_Q|z_j|^p dx\right)^{1/2}
       \|z_j\|_{W^{1,p}(Q)^2}^p.
\]
Since every cube \(Q\) is covered by at most \(C_N\) unit balls, where \(C_N\) is a finite number, and is therefore,
\[\int_Q |z_j|^p dx\le C_N\delta_j,\]
which follows that 
\[
 \int_Q|z_j|^{3p/2} dx
 \le\left(\int_Q|z_j|^p dx\right)^{1/2}
       \|z_j\|_{L^{2p}(Q)}^p
 \le C\left(\int_Q|z_j|^p dx\right)^{1/2}
       \|z_j\|_{W^{1,p}(Q)^2}^p.
\]
Summing over the cubes, we obtain
\[ \begin{aligned} \int_{\mathbb R^N}|z_j|^{3p/2}\,dx &\le C\delta_j^{1/2} \sum_Q\|z_j\|_{W^{1,p}(Q;\mathbb R^2)}^p\\ &=C\delta_j^{1/2} \|z_j\|_{W^{1,p}(\mathbb R^N;\mathbb R^2)}^p\\ &\le C\delta_j^{1/2}\longrightarrow0. \end{aligned} \]
Where we use the fact that \((z_j)\) is uniformly bounded in \(\mathcal{X}\).
For \(p<r<3p/2\), interpolation between the bounded
\(L^p\) norm and the vanishing \(L^{3p/2}\) norm gives
\(\|z_j\|_{L^r}\to0\).
For \(r>3p/2\), choose a finite \(s>r\) and interpolate
between \(L^{3p/2}\) and \(L^s\), using the uniform
\(L^s\) bound supplied by Lemma~\ref{lem:embeddings}.
Thus
\begin{equation}
 z_j\longrightarrow0
 \quad\text{in }L^r(\mathbb R^N;\mathbb R^2)
 \quad\text{for every }p<r<\infty.
 \label{eq:uniform-vanishing-lebesgue}
\end{equation}
On the other hand, for \(\varepsilon>0\), \[ \ell_\varepsilon(\xi)+\frac2p r_\varepsilon(\xi)=|s|^{p-2}t^2
 \left(\log(t^2+\varepsilon)+\frac2p\frac{t^2}{t^2+\varepsilon}\right).
\]
If the expression in parentheses is positive, then
\(t^2+\varepsilon>e^{-2/p}\), and hence
\[
 |t|^2>e^{-2/p}-1/2>0.
\]
Its positive part is therefore bounded by
\(C_\sigma|t|^\sigma\), uniformly in
\(0<\varepsilon\le1/2\).
The symmetric term satisfies the same estimate.
Weighted Young inequalities then yield
\begin{equation}
 \ell_\varepsilon(\xi)+\frac2p r_\varepsilon(\xi)
 \le C_\sigma|\xi|^{p+\sigma},
 \qquad
 \xi\in\mathbb R^2,\quad 0\le\varepsilon\le1/2.
 \label{eq:uniform-radial-source-bound}
\end{equation}
The case \(\varepsilon=0\) follows by passage to the limit. Therefore,
\[
 S_{q_j}(z_j)
 \le C_\sigma\|z_j\|_{p+\sigma}^{p+\sigma}+o(1)
 \longrightarrow0,
\]
which means \(\|z_j\|_{\bar{q}}\rightarrow 0\), it also follows that \(R_{\varepsilon_j}(z_j)\to 0\); \eqref{eq:radial-energy-identity} yields convergence of the
energy to zero, and so we complete the proof.
\end{proof}
\begin{proposition}[Uniform lower bounds for autonomous profiles]
\label{prop:uniform-profile-cost}
There exists \(\delta>0\), depending only on
\(p,\bar a,\bar b\), such that every nonzero critical point
\(Z\) of \(J_{\bar q,\varepsilon}\), with
\(0<\varepsilon\le1/2\), satisfies
\begin{equation}
 Z\in\mathcal N_{\bar q,\varepsilon},
 \qquad
 J_{\bar q,\varepsilon}(Z)\ge d_{\bar q,\varepsilon}\ge\delta,
 \qquad
 S_{\bar q}(Z)\ge C(V_0,p,\delta)>0.
 \label{eq:uniform-profile-cost}
\end{equation}
For every \(E_0>0\), there exists \(\eta_{E_0}>0\),
independent of \(\varepsilon\), such that every such critical
point with \(J_{\bar q,\varepsilon}(Z)\le E_0\) satisfies
\begin{equation}
 \sup_{y\in\mathbb R^N}
 \int_{B_1(y)}|Z|^p\,dx\ge\eta_{E_0}.
 \label{eq:uniform-profile-local-mass}
\end{equation}
\end{proposition}

\begin{proof}
If \(M(Z)=0\), then
\(L_\varepsilon(Z)=R_\varepsilon(Z)=0\).
Since \(Z\) is a critical point,
\eqref{eq:nehari-functional} gives \(S_{\bar q}(Z)=0\),
contrary to \(Z\ne0\).
Thus \(Z\in\mathcal N_{\bar q,\varepsilon}\), and
Corollary~\ref{cor:uniform-levels} proves
\eqref{eq:uniform-profile-cost}.

If \eqref{eq:uniform-profile-local-mass} were false,
there would exist \(0<\varepsilon_j\le1/2\) and nonzero
critical points \(Z_j\) of \(J_{\bar q,\varepsilon_j}\) such that
\[
 J_{\bar q,\varepsilon_j}(Z_j)\le E_0,
 \qquad
 \sup_{y\in\mathbb R^N}
 \int_{B_1(y)}|Z_j|^p\,dx\longrightarrow0.
\]
Lemma~\ref{thm:uniform-ps-bound} gives a uniform bound in
\(\mathcal X\).
Lemma~\ref{lem:uniform-ps-vanishing} then yields
\(\|Z_j\|_{\bar q}\to0\), contradicting
\eqref{eq:uniform-profile-cost}.
\end{proof}

\subsection{The finite profile decomposition}
\label{subsec:fixed-profile-decomposition}
\begin{theorem}[Global Palais--Smale decomposition]
\label{thm:fixed-global-decomposition}
Assume \eqref{hyp:continuous} and fix \(0<\varepsilon\le1/2\).
Every sequence satisfying \eqref{eq:fixed-ps-definition}
has, after passage to a subsequence, a decomposition
\begin{equation}
 z_j=z^0+\sum_{i=1}^k\mathsf U_{y_j^i}Z^i+e_j,
 \qquad \|e_j\|_q\longrightarrow0,
 \label{eq:fixed-profile-decomposition}
\end{equation}
where \(k\ge0\) is an integer,
\begin{equation}
 J_{q,\varepsilon}'(z^0)=0,\qquad
 Z^i\ne0,\quad J_{\bar q,\varepsilon}'(Z^i)=0
 \quad(1\le i\le k),
 \label{eq:fixed-profile-equations}
\end{equation}
and the centres satisfy
\begin{equation}
 |y_j^i|\longrightarrow\infty,\qquad
 |y_j^i-y_j^h|\longrightarrow\infty
 \quad(1\le i\ne h\le k).
 \label{eq:fixed-profile-separation}
\end{equation}
The norm and energy split as
\begin{align}
 S_q(z_j)
 &=S_q(z^0)+\sum_{i=1}^kS_{\bar q}(Z^i)+o(1),
 \label{eq:fixed-profile-norm}\\
 c
 &=J_{q,\varepsilon}(z^0)
   +\sum_{i=1}^kJ_{\bar q,\varepsilon}(Z^i).
 \label{eq:fixed-profile-energy}
\end{align}
In particular, \(c\ge0\), and
\begin{equation}
 k d_{\bar q,\varepsilon}
 \le c-J_{q,\varepsilon}(z^0),
 \qquad k\delta\le c,
 \label{eq:fixed-profile-count}
\end{equation}
where \(\delta\) is given by
Proposition~\ref{prop:uniform-profile-cost}.
\end{theorem}
\begin{proof}
   \emph{Step 1.}
Lemma~\ref{thm:uniform-ps-bound} bounds \((z_j)\) in
\(\mathcal X\). Hence, after passage to a subsequence,
\(z_j\rightharpoonup z^0\).
Lemma~\ref{lem:fixed-local-recovery} gives
\(J_{q,\varepsilon}'(z^0)=0\) and local strong
\(W^{1,p}\) convergence.
Taking a further diagonal subsequence gives almost everywhere
convergence of both components and their gradients.
Set \(w_j^0=z_j-z^0\).
By Lemma~\ref{lem:fixed-local-recovery}, the convergence is
strong in \(W^{1,p}_{\mathrm{loc}}\).
Passing to a further diagonal subsequence, both the components
and their gradients converge almost everywhere.
Lemma~\ref{lem:fixed-nonlinear-splitting}, together with the
power and flux splitting formulas applied to the components
and gradients, therefore gives
\begin{align}
 S_q(z_j)
 &=S_q(z^0)+S_q(w_j^0)+o(1),
 \label{eq:fixed-initial-norm-splitting}\\
 J_{q,\varepsilon}(z_j)
 &=J_{q,\varepsilon}(z^0)
   +J_{q,\varepsilon}(w_j^0)+o(1),
 \label{eq:fixed-initial-energy-splitting}\\
 \bigl\|
 J_{q,\varepsilon}'(z_j)
 -J_{q,\varepsilon}'(z^0)
 -J_{q,\varepsilon}'(w_j^0)
 \bigr\|_{\mathcal X^*}
 &\longrightarrow0.
 \label{eq:fixed-initial-derivative-splitting}
\end{align}
In fact, \eqref{eq:fixed-initial-norm-splitting} and \eqref{eq:fixed-initial-energy-splitting} are trivial. To justify \eqref{eq:fixed-initial-derivative-splitting},
the power-flux splitting formula gives
\[
 \bigl\|
 |\nabla u_j|^{p-2}\nabla u_j
 -|\nabla u^0|^{p-2}\nabla u^0
 -|\nabla(w_j^0)_1|^{p-2}\nabla(w_j^0)_1
 \bigr\|_{L^{p'}}
 \longrightarrow0,
\]
and the analogous convergence for the scalar fluxes
\(|u_j|^{p-2}u_j\), as well as for the second component.
Since \(a,b\in L^\infty(\mathbb R^N)\), H\"older's inequality
therefore gives convergence of the principal derivative
defects in \(\mathcal X^*\).
For the nonlinear term,
\eqref{eq:fixed-source-splitting} gives convergence in
\(L^{p'}+L^{r'}\), which embeds continuously into
\(\mathcal X^*\) by Lemma~\ref{lem:embeddings}.
Combining these estimates proves
\eqref{eq:fixed-initial-derivative-splitting}.
Since $w_j^0\rightharpoonup0$, it is easy to verify that
\[|J_{q.\varepsilon}(w_j^0)-J_{\bar{q}.\varepsilon}(w_j^0)|\to 0,\quad \|J_{q,\varepsilon}'(w_j^0)-J_{\bar{q}.\varepsilon}(w_j^0)\|_{\mathcal{X}^*}\to 0,\]
which immediately implies that
\begin{equation}
 J_{\bar q,\varepsilon}(w_j^0)\longrightarrow
 c-J_{q,\varepsilon}(z^0),
 \qquad
 J_{\bar q,\varepsilon}'(w_j^0)\longrightarrow0
 \quad\text{in }\mathcal X^*.
 \label{eq:fixed-autonomous-remainder}
\end{equation}
Every bounded Palais--Smale sequence has nonnegative limiting
energy, since its radial residual tends to zero and
\eqref{eq:radial-energy-identity} applies.
Every nonzero critical point likewise belongs to its Nehari
set and has positive energy by Corollary~\ref{cor:uniform-levels}.

\emph{Step 2.}
Suppose first that
\[
 \sup_y\int_{B_1(y)}|w_j^0|^p\,dx\longrightarrow0.
\]
Lemma~\ref{lem:uniform-ps-vanishing} gives \(w_j^0\to0\) in
\(\mathcal X\), and the conclusion holds with \(k=0\).
Otherwise, after passage to a subsequence, there exist
\(\eta>0\) and centres \(y_j^1\) such that
\[
 \int_{B_1(y_j^1)}|w_j^0|^p\,dx\ge\eta.
\]
Local strong convergence of \(w_j^0\) to zero forces
\(|y_j^1|\to\infty\).
The translated sequence \(\mathsf U_{-y_j^1}w_j^0\)
is bounded and is Palais--Smale for
\(J_{\bar q,\varepsilon}\) by translation invariance.
After passage to a subsequence, it converges weakly to some
\(Z^1\).
Lemma~\ref{lem:fixed-local-recovery} gives strong local
\(W^{1,p}\) convergence and
\(J_{\bar q,\varepsilon}'(Z^1)=0\).
The local mass inequality implies \(Z^1\ne0\).
Define
\[
 w_j^1=w_j^0-\mathsf U_{y_j^1}Z^1.
\]
As before, after passing to a further subsequence, the
components and gradients converge almost everywhere in the
translated coordinates. Applying the same splitting formulas yields
\begin{align}
 S_{\bar q}(w_j^0)
 &=S_{\bar q}(Z^1)+S_{\bar q}(w_j^1)+o(1),
 \label{eq:fixed-one-step-norm}\\
 J_{\bar q,\varepsilon}(w_j^0)
 &=J_{\bar q,\varepsilon}(Z^1)
   +J_{\bar q,\varepsilon}(w_j^1)+o(1),
 \label{eq:fixed-one-step-energy}\\
 \bigl\|
 J_{\bar q,\varepsilon}'(w_j^0)
 -\mathsf U_{y_j^1}J_{\bar q,\varepsilon}'(Z^1)
 -J_{\bar q,\varepsilon}'(w_j^1)
 \bigr\|_{\mathcal X^*}
 &\longrightarrow0.
 \label{eq:fixed-one-step-derivative}
\end{align}
The density splitting holds in \(L^1\) as well.
Consequently, \(w_j^1\) is a bounded autonomous Palais--Smale
sequence. It tends strongly to zero locally in the original
coordinates and in the coordinates centred at \(y_j^1\).

\emph{Step 3.}
Assume that \(i\) profiles have been extracted, and write
\[
 w_j^i=z_j-z^0-\sum_{h=1}^i\mathsf U_{y_j^h}Z^h.
\]
By induction, \((w_j^i)\) is a bounded Palais--Smale sequence
for \(J_{\bar q,\varepsilon}\). With \(y_j^0=0\) denoting
the original frame, we have
\[
 \mathsf U_{-y_j^h}w_j^i\to0
 \quad\text{in }W^{1,p}_{\mathrm{loc}}(\mathbb R^N;\mathbb R^2),
 \qquad 0\le h\le i.
\]
If the vanishing condition fails, then, after passage to a
subsequence, there exist \(\eta>0\) and centres \(y_j^{i+1}\)
such that
\[
 \int_{B_1(y_j^{i+1})}|w_j^i|^p\,dx\ge\eta.
\]
We claim that \(|y_j^{i+1}-y_j^h|\to\infty\) for every
\(0\le h\le i\). Otherwise, for some \(h\) and along a
subsequence, \(|y_j^{i+1}-y_j^h|\le C\). Hence
\[
 \eta
 \le\int_{B_1(y_j^{i+1}-y_j^h)}
       |\mathsf U_{-y_j^h}w_j^i|^p\,dx
 \le\int_{B_{C+1}(0)}
       |\mathsf U_{-y_j^h}w_j^i|^p\,dx
 \longrightarrow0,
\]
a contradiction. In particular, \(h=0\) gives
\(|y_j^{i+1}|\to\infty\).

By translation invariance, \((\mathsf U_{-y_j^{i+1}}w_j^i)\)
is again a bounded Palais--Smale sequence for
\(J_{\bar q,\varepsilon}\). As in Step~2, after passage to a
subsequence, Lemma~\ref{lem:fixed-local-recovery} gives
\[
 \mathsf U_{-y_j^{i+1}}w_j^i\to Z^{i+1}
 \quad\text{in }W^{1,p}_{\mathrm{loc}},
 \qquad J_{\bar q,\varepsilon}'(Z^{i+1})=0.
\]
The local mass bound passes to the limit on \(B_1(0)\), giving
\(\int_{B_1(0)}|Z^{i+1}|^p\,dx\ge\eta\). Thus \(Z^{i+1}\ne0\).

Set \(w_j^{i+1}=w_j^i-\mathsf U_{y_j^{i+1}}Z^{i+1}\).
In the new frame,
\[
 \mathsf U_{-y_j^{i+1}}w_j^{i+1}
 =\mathsf U_{-y_j^{i+1}}w_j^i-Z^{i+1}
 \to0
 \quad\text{in }W^{1,p}_{\mathrm{loc}}.
\]
In each old frame \(0\le h\le i\),
\[
 \mathsf U_{-y_j^h}w_j^{i+1}
 =\mathsf U_{-y_j^h}w_j^i
  -\mathsf U_{y_j^{i+1}-y_j^h}Z^{i+1}
 \to0
 \quad\text{in }W^{1,p}_{\mathrm{loc}}.
\]
Indeed, the first term tends to zero by induction, and the
second by separation of the centres, since \(Z^{i+1}\) and its
gradient belong to \(L^p\). The new remainder is bounded in
\(\mathcal X\), and the energy and derivative splittings from
Step~2 show that it is again a Palais--Smale sequence for
\(J_{\bar q,\varepsilon}\). This closes the induction.

At every step, the counterparts of
\eqref{eq:fixed-one-step-norm}--\eqref{eq:fixed-one-step-derivative}
hold. In particular, the remainders are bounded, and their
limiting energies satisfy
\[
 0\le\lim_{j\to\infty}J_{\bar q,\varepsilon}(w_j^i)
 =c-J_{q,\varepsilon}(z^0)
   -\sum_{h=1}^iJ_{\bar q,\varepsilon}(Z^h).
\]
Here nonnegativity follows from the radial energy identity
and the Palais--Smale property of the remainders.
By Proposition~\ref{prop:uniform-profile-cost}, each extraction
lowers the limiting energy by at least \(\delta>0\).
The procedure must therefore stop after finitely many steps:
otherwise, choosing
\[
 i>\frac{c-J_{q,\varepsilon}(z^0)}{\delta}
\]
would give a remainder with negative limiting energy,
a contradiction.

\emph{Step 4.}
For the terminal index \(k\),
Lemma~\ref{lem:uniform-ps-vanishing} gives
\[
 w_j^k\longrightarrow0\quad\text{in }\mathcal X,
 \qquad
 J_{\bar q,\varepsilon}(w_j^k)\longrightarrow0.
\]
Set \(e_j=w_j^k\).
Summing the initial and successive splittings, together with
\eqref{eq:compact-potential-energy}, yields
\eqref{eq:fixed-profile-norm} and
\eqref{eq:fixed-profile-energy}.

For completeness, the density splittings also give, for every
\(H\) in \eqref{eq:fixed-density-family},
\begin{equation}
 \left\|
 H(z_j)-H(z^0)
 -\sum_{i=1}^k\mathsf U_{y_j^i}H(Z^i)
 \right\|_{L^1}\longrightarrow0.
 \label{eq:fixed-profile-densities}
\end{equation}
Indeed, the growth bounds at fixed \(\varepsilon>0\)
and \(w_j^k\to0\) in \(\mathcal X\) imply
\(H(w_j^k)\to0\) in \(L^1\).

Finally, \eqref{eq:uniform-profile-cost} and the energy
splitting imply \eqref{eq:fixed-profile-count} and we complete the proof.
\end{proof}

\begin{remark}\label{rem:profile-energy-meaning}
The autonomous profiles in
Theorem~\ref{thm:fixed-global-decomposition} are nonzero
critical points, not necessarily ground states.
Their energies are at least \(d_{\bar q,\varepsilon}\),
but need not equal it.
Thus \eqref{eq:fixed-profile-energy} expresses the total
energy as a sum of critical values, rather than integer
multiples of the ground-state level.
Neither the decomposition nor the profile count requires
a sign assumption on the original sequence.
\end{remark}

\subsection{Reverse construction and the exact compactness threshold}
\label{subsec:fixed-sharp-threshold}

\begin{theorem}[Autonomous ground states]
\label{thm:fixed-autonomous-ground-state}
For every \(0<\varepsilon\le1/2\), the level
\(d_{\bar q,\varepsilon}\) is attained by a nonzero,
componentwise nonnegative critical point of
\(J_{\bar q,\varepsilon}\).
Every Palais--Smale sequence for \(J_{\bar q,\varepsilon}\)
at level \(d_{\bar q,\varepsilon}\) admits a subsequence
converging strongly in \(\mathcal X\), after suitable
translations, to a ground state.
The same compactness conclusion holds for every minimising
sequence of \(J_{\bar q,\varepsilon}\) on
\(\mathcal N_{\bar q,\varepsilon}\).
\end{theorem}
\begin{proof}
Corollary~\ref{cor:uniform-levels} gives
\(0<d_{\bar q,\varepsilon}<\infty\).
The Nehari set is closed, hence complete in \(\mathcal X\).
Indeed, if \(z_j\in\mathcal N_{\bar q,\varepsilon}\)
and \(z_j\to z\) in \(\mathcal X\), continuity and
\eqref{eq:nehari-gaps} give
\(G_{\bar q,\varepsilon}(z)=0\) and \(z\ne0\).
Moreover, \(M(z)=0\) would imply
\(G_{\bar q,\varepsilon}(z)=S_{\bar q}(z)>0\),
a contradiction.
Thus \(z\in\mathcal N_{\bar q,\varepsilon}\).

Let $(\widehat z_j)$ be any Nehari minimising sequence.
Choose $\lambda_j>0$ tending to zero so that
$J_{\bar q,\varepsilon}(\widehat z_j)\le d_{\bar q,\varepsilon}+\lambda_j^2$.
Ekeland's variational principle on the complete Nehari set
provides $\zeta_j\in\mathcal N_{\bar q,\varepsilon}$ satisfying
\[
 \|\zeta_j-\widehat z_j\|_{\bar q}\le\lambda_j,
 \qquad J_{\bar q,\varepsilon}(\zeta_j)
                  \le J_{\bar q,\varepsilon}(\widehat z_j),
\]
and
\[
 J_{\bar q,\varepsilon}(w)
 \ge J_{\bar q,\varepsilon}(\zeta_j)
                   -\lambda_j\|w-\zeta_j\|_{\bar q}
 \qquad(w\in\mathcal N_{\bar q,\varepsilon}).
\]
By Proposition~\ref{prop:radial-projection}, the Nehari set is
a \(C^1\) hypersurface. Applying the last inequality to
\(C^1\) curves through \(\zeta_j\) gives
\begin{equation}
 |\langle J_{\bar q,\varepsilon}'(\zeta_j),h\rangle|
 \le\lambda_j\|h\|_{\bar q}
 \quad\text{if }DG_{\bar q,\varepsilon}(\zeta_j)[h]=0.
 \label{eq:fixed-tangential-ekeland}
\end{equation}
Indeed, for such \(h\), choose a \(C^1\) curve
\(\gamma:(-\delta,\delta)\to\mathcal N_{\bar q,\varepsilon}\)
with \(\gamma(0)=\zeta_j\) and \(\gamma'(0)=h\).
The preceding variational inequality yields
\[
 J_{\bar q,\varepsilon}(\gamma(t))
 -J_{\bar q,\varepsilon}(\zeta_j)
 \ge-\lambda_j\|\gamma(t)-\zeta_j\|_{\bar q}.
\]
Dividing by \(t>0\) and letting \(t\downarrow0\) gives
\(\langle J_{\bar q,\varepsilon}'(\zeta_j),h\rangle
\ge-\lambda_j\|h\|_{\bar q}\).
Replacing \(h\) by \(-h\) proves
\eqref{eq:fixed-tangential-ekeland}.
The Nehari bound in Lemma~\ref{thm:uniform-ps-bound} gives
boundedness of \((\zeta_j)\) in \(\mathcal X\).
By \eqref{eq:radial-transversality} and the Nehari energy identity,
\[
 |DG_{\bar q,\varepsilon}(\zeta_j)[\zeta_j]|
 \ge2R_\varepsilon(\zeta_j)
 =p^2J_{\bar q,\varepsilon}(\zeta_j)
 \ge p^2 d_{\bar q,\varepsilon}
 \ge p^2\delta.
\]
For fixed \(\varepsilon>0\),
\(DG_{\bar q,\varepsilon}\) is bounded on bounded subsets
of \(\mathcal X\), by \eqref{eq:nehari-functional},
the fixed-parameter source bound, and
 $|D_\xi r_\varepsilon(\xi)|\le C|\xi|^{p-1}.$
 Hence the projections
\[
 P_jh=h-
 \frac{DG_{\bar q,\varepsilon}(\zeta_j)[h]}
      {DG_{\bar q,\varepsilon}(\zeta_j)[\zeta_j]}
 \,\zeta_j
\]
onto \(T_{\zeta_j}\mathcal N_{\bar q,\varepsilon}\)
satisfy
\[
 \|P_jh\|_{\bar q}
 \le
 \left(
 1+
 \frac{\|DG_{\bar q,\varepsilon}(\zeta_j)\|_{\mathcal X^*}
       \|\zeta_j\|_{\bar q}}
      {|DG_{\bar q,\varepsilon}(\zeta_j)[\zeta_j]|}
 \right)\|h\|_{\bar q}
 \le C\|h\|_{\bar q}.
\]
Since
\(\langle J_{\bar q,\varepsilon}'(\zeta_j),\zeta_j\rangle=0\),
applying \eqref{eq:fixed-tangential-ekeland} to \(P_jh\) gives
\[
 \bigl|\langle J_{\bar q,\varepsilon}'(\zeta_j),h\rangle\bigr|
 =
 \bigl|\langle J_{\bar q,\varepsilon}'(\zeta_j),P_jh\rangle\bigr|
 \le\lambda_j\|P_jh\|_{\bar q}
 \le C\lambda_j\|h\|_{\bar q}.
\]
Thus
\(\|J_{\bar q,\varepsilon}'(\zeta_j)\|_{\mathcal X^*}\to0\),
and \((\zeta_j)\) is a Palais--Smale sequence in \(\mathcal X\)
at level \(d_{\bar q,\varepsilon}\).

Apply Theorem~\ref{thm:fixed-global-decomposition} with
\(q=\bar q\). The energy splitting gives
\[
 d_{\bar q,\varepsilon}
 =J_{\bar q,\varepsilon}(z^0)
  +\sum_{i=1}^kJ_{\bar q,\varepsilon}(Z^i).
\]
Every nonzero state in this decomposition belongs to
\(\mathcal N_{\bar q,\varepsilon}\) and therefore has energy
at least \(d_{\bar q,\varepsilon}>0\).
Hence exactly one nonzero state occurs, and its energy equals
\(d_{\bar q,\varepsilon}\).
If this state is \(z^0\), then \(\zeta_j\to z^0\) strongly.
Otherwise,
\[
 \zeta_j=\mathsf U_{y_j^1}Z^1+o_{\mathcal X}(1),
 \qquad
 \mathsf U_{-y_j^1}\zeta_j\to Z^1
 \quad\text{strongly in }\mathcal X.
\]
Thus, in either case, there exist translations \(y_j\)
and a ground state \(Z\) such that
\(\mathsf U_{-y_j}\zeta_j\to Z\) strongly in \(\mathcal X\).

Since translations are isometries,
\[
 \|\mathsf U_{-y_j}\widehat z_j-Z\|_{\bar q}
 \le \|\widehat z_j-\zeta_j\|_{\bar q}
    +\|\mathsf U_{-y_j}\zeta_j-Z\|_{\bar q}
 \longrightarrow0.
\]
This proves attainment and compactness modulo translations
for every Nehari minimising sequence.
The same energy-counting argument applies directly to any
Palais--Smale sequence at level \(d_{\bar q,\varepsilon}\),
without the Ekeland approximation.

Finally, replacing an attained minimiser $(u,v)$ by
$(|u|,|v|)$ preserves the principal energy and every mixed
density in the continuous setting. The new pair is again a
Nehari minimiser and hence a free critical point by
Proposition~\ref{prop:radial-projection}.
\end{proof}

\begin{proposition}[Separated critical points generate Palais--Smale sequences]
\label{prop:fixed-profile-gluing}
Fix \(0<\varepsilon\le1/2\).
Let \(z^0\in\mathcal X\) satisfy
\(J_{q,\varepsilon}'(z^0)=0\), and let
\(Z^1,\ldots,Z^k\) be nonzero critical points of
\(J_{\bar q,\varepsilon}\), where \(k\ge1\).
For any centres satisfying \eqref{eq:fixed-profile-separation},
the sequence
\begin{equation}
 z_j=z^0+\sum_{i=1}^k\mathsf U_{y_j^i}Z^i
 \label{eq:fixed-glued-sequence}
\end{equation}
is a Palais--Smale sequence for \(J_{q,\varepsilon}\) at level
\[
 J_{q,\varepsilon}(z^0)
 +\sum_{i=1}^kJ_{\bar q,\varepsilon}(Z^i).
\]
Moreover, \eqref{eq:fixed-profile-norm} and
\eqref{eq:fixed-profile-densities} hold, and
\(z_j\rightharpoonup z^0\) in \(\mathcal X\),
but \((z_j)\) has no strongly convergent subsequence.
\end{proposition}

\begin{proof}
We divide the proof into three steps.

\emph{Step 1.}
Throughout the proof, we use the reference norm
\(\|\cdot\|_{\bar q}\), which is equivalent to
\(\|\cdot\|_q\).

Fix \(0<\delta<1\). Choose
\(\phi^0,\phi^1,\ldots,\phi^k
\in C_c^\infty(\mathbb R^N;\mathbb R^2)\) such that
\[
 \|z^0-\phi^0\|_{\bar q}
 +\sum_{i=1}^k\|Z^i-\phi^i\|_{\bar q}<\delta,
\]
and set
\[
 \widetilde z_j
 =\phi^0+\sum_{i=1}^k\mathsf U_{y_j^i}\phi^i.
\]
Since translations preserve the reference norm,
\[
 \|z_j-\widetilde z_j\|_{\bar q}<\delta.
\]

For these fixed approximations, the separation of the centres
ensures that the supports of
\(\phi^0,\mathsf U_{y_j^1}\phi^1,\ldots,
\mathsf U_{y_j^k}\phi^k\)
are pairwise disjoint for all sufficiently large \(j\).
Locality therefore gives
\begin{align*}
 S_q(\widetilde z_j)
 &=S_q(\phi^0)
   +\sum_{i=1}^kS_q(\mathsf U_{y_j^i}\phi^i),\\
 J_{q,\varepsilon}(\widetilde z_j)
 &=J_{q,\varepsilon}(\phi^0)
   +\sum_{i=1}^k
     J_{q,\varepsilon}(\mathsf U_{y_j^i}\phi^i),\\
 J_{q,\varepsilon}'(\widetilde z_j)
 &=J_{q,\varepsilon}'(\phi^0)
   +\sum_{i=1}^k
     J_{q,\varepsilon}'(\mathsf U_{y_j^i}\phi^i).
\end{align*}
Likewise, for every \(H\) in
\eqref{eq:fixed-density-family},
\[
 H(\widetilde z_j)
 =H(\phi^0)
  +\sum_{i=1}^k\mathsf U_{y_j^i}H(\phi^i).
\]

All profiles, their approximations, and their translated sums
remain in a common bounded subset of \(\mathcal X\).
The continuity estimates underlying
Lemma~\ref{lem:fixed-nonlinear-splitting}, together with
the Sobolev embeddings and the boundedness of \(a,b\),
show that the approximation errors tend to zero uniformly
in \(j\) as \(\delta\downarrow0\).
First letting \(j\to\infty\) and then
\(\delta\downarrow0\), we obtain
\begin{align*}
 S_q(z_j)
 &=S_q(z^0)
   +\sum_{i=1}^kS_q(\mathsf U_{y_j^i}Z^i)+o(1),\\
 J_{q,\varepsilon}(z_j)
 &=J_{q,\varepsilon}(z^0)
   +\sum_{i=1}^k
     J_{q,\varepsilon}(\mathsf U_{y_j^i}Z^i)+o(1),
\end{align*}
and
\[
 \left\|
 J_{q,\varepsilon}'(z_j)
 -J_{q,\varepsilon}'(z^0)
 -\sum_{i=1}^k
   J_{q,\varepsilon}'(\mathsf U_{y_j^i}Z^i)
 \right\|_{\mathcal X^*}
 \longrightarrow0.
\]
The same approximation argument proves
\eqref{eq:fixed-profile-densities}.

\emph{Step 2.}
Since \(|y_j^i|\to\infty\),
\[
 \mathsf U_{y_j^i}Z^i\rightharpoonup0
 \quad\text{in }\mathcal X.
\]
Lemma~\ref{lem:fixed-local-recovery} and translation
invariance of the autonomous functional give
\begin{align*}
 S_q(\mathsf U_{y_j^i}Z^i)
 &=S_{\bar q}(Z^i)+o(1),\\
 J_{q,\varepsilon}(\mathsf U_{y_j^i}Z^i)
 &=J_{\bar q,\varepsilon}(Z^i)+o(1).
\end{align*}
Substituting these identities into the splittings from Step~1
proves \eqref{eq:fixed-profile-norm} and
\[
 J_{q,\varepsilon}(z_j)
 \longrightarrow
 J_{q,\varepsilon}(z^0)
 +\sum_{i=1}^kJ_{\bar q,\varepsilon}(Z^i).
\]

For the derivatives,
Lemma~\ref{lem:fixed-local-recovery} yields
\[
 \left\|
 J_{q,\varepsilon}'(\mathsf U_{y_j^i}Z^i)
 -J_{\bar q,\varepsilon}'(\mathsf U_{y_j^i}Z^i)
 \right\|_{\mathcal X^*}\longrightarrow0.
\]
By translation invariance and
\(J_{\bar q,\varepsilon}'(Z^i)=0\),
\[
 J_{\bar q,\varepsilon}'(\mathsf U_{y_j^i}Z^i)=0.
\]
Together with the derivative splitting from Step~1 and
\(J_{q,\varepsilon}'(z^0)=0\), this gives
\[
 \|J_{q,\varepsilon}'(z_j)\|_{\mathcal X^*}
 \longrightarrow0.
\]
Thus \((z_j)\) is a Palais--Smale sequence at the stated level.

\emph{Step 3.}
The weak convergence of the translated profiles gives
\[
 z_j\rightharpoonup z^0
 \quad\text{in }\mathcal X.
\]
Applying the norm-splitting argument from Step~1 to
\(z_j-z^0\), with the constant potentials \(\bar q\), yields
\[
 \|z_j-z^0\|_{\bar q}^p
 =\sum_{i=1}^k\|Z^i\|_{\bar q}^p+o(1).
\]
Since \(k\ge1\) and every \(Z^i\ne0\), the right-hand side
has a strictly positive limit.
Any strongly convergent subsequence would necessarily
converge to the weak limit \(z^0\), contradicting this identity.
Hence \((z_j)\) has no strongly convergent subsequence
in \(\mathcal X\).
\end{proof}

\begin{corollary}[Failure levels and the exact first threshold]
\label{cor:fixed-ps-failure-levels}
Fix \(0<\varepsilon\le1/2\), and define
\[
 \mathfrak C_{q,\varepsilon}
 =\{J_{q,\varepsilon}(z):
   z\in\mathcal X,\ J_{q,\varepsilon}'(z)=0\}.
\]
In particular, \(0\in\mathfrak C_{q,\varepsilon}\).
We say that \(J_{q,\varepsilon}\) satisfies the
Palais--Smale condition at level \(c\) if every sequence
satisfying \eqref{eq:fixed-ps-definition} admits a strongly
convergent subsequence.
The levels at which this condition fails are precisely
\begin{equation}
 \left\{
 c_0+\sum_{i=1}^k c_i:
 k\ge1,\quad
 c_0\in\mathfrak C_{q,\varepsilon},\quad
 c_i\in\mathfrak C_{\bar q,\varepsilon}\setminus\{0\}
 \ (1\le i\le k)
 \right\}.
 \label{eq:fixed-ps-failure-set}
\end{equation}
Consequently, \(J_{q,\varepsilon}\) satisfies the
Palais--Smale condition at every \(c<d_{\bar q,\varepsilon}\)
and fails it at \(c=d_{\bar q,\varepsilon}\).
Thus \(d_{\bar q,\varepsilon}\) is the exact first
loss-of-compactness level.
\end{corollary}

\begin{proof}
Suppose that the Palais--Smale condition fails at level \(c\).
Choose a Palais--Smale sequence \((z_j)\) at this level
with no strongly convergent subsequence.
Theorem~\ref{thm:fixed-global-decomposition} gives
\[
 z_j=z^0+\sum_{i=1}^k\mathsf U_{y_j^i}Z^i
      +o_{\mathcal X}(1).
\]
If \(k=0\), then \(z_j\to z^0\) strongly, a contradiction.
Thus \(k\ge1\), and
\[
 c=J_{q,\varepsilon}(z^0)
   +\sum_{i=1}^kJ_{\bar q,\varepsilon}(Z^i),
\]
where
\[
 J_{q,\varepsilon}(z^0)\in\mathfrak C_{q,\varepsilon},
 \qquad
 J_{\bar q,\varepsilon}(Z^i)
 \in\mathfrak C_{\bar q,\varepsilon}\setminus\{0\}.
\]
Hence \(c\) belongs to \eqref{eq:fixed-ps-failure-set}.

Conversely, let
\[
 c=c_0+\sum_{i=1}^k c_i
\]
be any level in \eqref{eq:fixed-ps-failure-set}.
Choose critical points \(z^0,Z^1,\ldots,Z^k\) representing
these values, and set
\[
 y_j^i=ij\mathbf e,
 \qquad
 z_j=z^0+\sum_{i=1}^k\mathsf U_{y_j^i}Z^i,
\]
where \(\mathbf e\in\mathbb R^N\) is a fixed unit vector.
Then
\[
 |y_j^i|\to\infty,
 \qquad
 |y_j^i-y_j^h|=|i-h|j\to\infty
 \quad(i\ne h).
\]
Proposition~\ref{prop:fixed-profile-gluing} gives
\[
 J_{q,\varepsilon}(z_j)\to c,
 \qquad
 \|J_{q,\varepsilon}'(z_j)\|_{\mathcal X^*}\to0,
\]
but \((z_j)\) has no strongly convergent subsequence.
Thus the Palais--Smale condition fails at \(c\).

Finally, every critical point satisfies
\[
 J_{q,\varepsilon}(z)
 =\frac2{p^2}R_\varepsilon(z)\ge0,
\]
and every nonzero autonomous critical point satisfies
\[
 J_{\bar q,\varepsilon}(Z)\ge d_{\bar q,\varepsilon}.
\]
Consequently, every failure level obeys
\[
 c=c_0+\sum_{i=1}^k c_i
 \ge k d_{\bar q,\varepsilon}\ge d_{\bar q,\varepsilon}.
\]
The Palais--Smale condition therefore holds below
\(d_{\bar q,\varepsilon}\).

By Theorem~\ref{thm:fixed-autonomous-ground-state},
there exists a nonzero critical point \(Z\) with
\(J_{\bar q,\varepsilon}(Z)=d_{\bar q,\varepsilon}\).
Taking \(z^0=0\), \(k=1\), and \(z_j=\mathsf U_{j\mathbf e}Z\)
in the gluing proposition proves failure at
\(c=d_{\bar q,\varepsilon}\).
\end{proof}

\begin{corollary}[Parameters with a positive limit]
\label{cor:positive-regularisation-limit}
Let \(\varepsilon_j\to\varepsilon_*>0\), with
\(0<\varepsilon_j\le1/2\), and suppose that
\[
 J_{q,\varepsilon_j}(z_j)\longrightarrow c,
 \qquad
 J_{q,\varepsilon_j}'(z_j)\longrightarrow0
 \quad\text{in }\mathcal X^*.
\]
Then the conclusions of
Theorem~\ref{thm:fixed-global-decomposition} hold with
\(\varepsilon=\varepsilon_*\).
In \eqref{eq:fixed-profile-densities}, the density applied
to \(z_j\) may also be evaluated at \(\varepsilon_j\).
\end{corollary}

\begin{proof}
Lemma~\ref{thm:uniform-ps-bound} gives boundedness of
\((z_j)\) in \(\mathcal X\).
Choose \(0<\delta<\varepsilon_*\).
For all sufficiently large \(j\), both
\(\varepsilon_j\) and \(\varepsilon_*\) lie in
\([\delta,1/2]\).

For \(\varepsilon\in[\delta,1/2]\), we have
\[
 |\partial_\varepsilon h_\varepsilon(s)|
 =\frac{s^2}{s^2+\varepsilon}
 \le\delta^{-1}|s|^2,
 \qquad
 |\partial_\varepsilon h_\varepsilon'(s)|
 =\frac{2|s|\varepsilon}{(s^2+\varepsilon)^2}
 \le2\delta^{-1}|s|.
\]
The mean value theorem and weighted Young inequalities
therefore give
\begin{align*}
 |\ell_{\varepsilon_j}(\xi)-\ell_{\varepsilon_*}(\xi)|
 &\le C_\delta|\varepsilon_j-\varepsilon_*||\xi|^p,\\
 |D_\xi\ell_{\varepsilon_j}(\xi)
   -D_\xi\ell_{\varepsilon_*}(\xi)|
 &\le C_\delta|\varepsilon_j-\varepsilon_*||\xi|^{p-1}.
\end{align*}
Since the principal part is independent of \(\varepsilon\),
integration and H\"older's inequality yield
\begin{align*}
 |J_{q,\varepsilon_j}(z_j)-J_{q,\varepsilon_*}(z_j)|
 &\le C_\delta|\varepsilon_j-\varepsilon_*|
       \|z_j\|_{L^p}^p\longrightarrow0,\\
 \|J_{q,\varepsilon_j}'(z_j)
   -J_{q,\varepsilon_*}'(z_j)\|_{\mathcal X^*}
 &\le C_\delta|\varepsilon_j-\varepsilon_*|
       \|z_j\|_{L^p}^{p-1}\longrightarrow0.
\end{align*}
Together with the hypotheses, these estimates imply
\[
 J_{q,\varepsilon_*}(z_j)\longrightarrow c,
 \qquad
 J_{q,\varepsilon_*}'(z_j)\longrightarrow0
 \quad\text{in }\mathcal X^*.
\]
Thus \((z_j)\) is a Palais--Smale sequence for the fixed
functional \(J_{q,\varepsilon_*}\), and
Theorem~\ref{thm:fixed-global-decomposition} applies.

Finally, the Lipschitz continuity of the positive and negative
parts gives the same parameter estimate for
\(\ell_\varepsilon^\pm\).
For the radial density,
\[
 |\partial_\varepsilon r_\varepsilon(\xi)|
 \le\delta^{-1}m(\xi)\le C_\delta|\xi|^p,
\]
whereas \(m\) is independent of \(\varepsilon\).
Consequently, for every parameter-dependent density
\(H_\varepsilon\) in \eqref{eq:fixed-density-family},
\[
 \|H_{\varepsilon_j}(z_j)-H_{\varepsilon_*}(z_j)\|_{L^1}
 \le C_\delta|\varepsilon_j-\varepsilon_*|
       \|z_j\|_{L^p}^p\longrightarrow0.
\]
Hence \(H_{\varepsilon_*}(z_j)\) may be replaced by
\(H_{\varepsilon_j}(z_j)\) in
\eqref{eq:fixed-profile-densities}, while the limiting
profiles retain the parameter \(\varepsilon_*\).
\end{proof}

\begin{remark}[Uniform estimates and the singular limit]
\label{rem:fixed-versus-singular-compactness}
The bounds in Lemma~\ref{thm:uniform-ps-bound},
Proposition~\ref{prop:uniform-profile-cost}, and
Lemma~\ref{lem:uniform-ps-vanishing} are uniform as
\(\varepsilon\downarrow0\).
In particular, a fixed-parameter decomposition at level
\(c\le E_0\) contains at most \(E_0/\delta\) escaping profiles.
However, the splitting and parameter estimates above
are not uniform near zero, so they do not directly yield
a decomposition when \(\varepsilon_j\to0\).
The singular limit requires separate control of the negative
logarithmic energy and the residual, while respecting the
finite-energy domain \(\mathcal D_q\).
\end{remark}

\section{Singular global compactness and removal of regularisation}
\label{sec:singular-global-compactness}

Lions' concentration--compactness principle \cite{Lio84a}
provides a framework for studying loss of compactness on
unbounded domains. For semilinear elliptic equations on
\(\mathbb R^N\) with asymptotically constant coefficients,
Kuzin and Pohozaev \cite{KP97} described the loss of compactness
of Palais--Smale sequences through translated solutions of the
limiting equation. Hua, Li and Wang obtained a discrete
counterpart on Cayley graphs of groups of polynomial growth
(Theorem~1.3 of \cite{HLW23}). In the present system, the
profiles at infinity likewise solve an autonomous equation.
The additional difficulty is that the regularization parameter
varies along the sequence, and subtracting a profile need not
preserve the finite-energy domain at the singular limit.

Throughout this section, we assume \eqref{hyp:continuous},
with \(p=N>4\), retain the notation of the preceding section,
and set \(d_{\bar q,0}=d_{\bar q,0}\).
We prove a finite profile decomposition for
\(\varepsilon_j\to0\), allowing \(\varepsilon_j=0\).
A uniform commutator estimate controls the equation residual
under slowly expanding cutoffs. The resulting remainder
converges strongly to zero in \(\mathcal X\), and both its
positive and negative logarithmic energies vanish.

\subsection{The equation residual and local recovery}
\label{subsec:singular-residual-framework}

For \(z\in\mathcal X\) and \(0\le\varepsilon\le1/2\),
define the equation residual by
\begin{equation}
 \langle\mathfrak R_{q,\varepsilon}(z),h\rangle
 =\frac1p S_q'(z)[h]
 -\frac1p\int_{\mathbb R^N}D_\xi\ell_\varepsilon(z)\cdot h\,dx,
 \quad
 h\in C_c^\infty(\mathbb R^N;\mathbb R^2).
 \label{eq:singular-residual-definition}
\end{equation}
By \eqref{eq:derivative-growth} and
Lemma~\ref{lem:embeddings},
\(D_\xi\ell_\varepsilon(z)\in L^{p'}_{\mathrm{loc}}\).
Thus, by density, the same formula holds for compactly
supported \(W^{1,p}\) test pairs.
For \(\varepsilon>0\), the residual coincides with
\(J_{q,\varepsilon}'(z)\in\mathcal X^*\).
At \(\varepsilon=0\), we call \(z\in\mathcal D_q\) admissible if
\(\mathfrak R_{q,0}(z)\) extends continuously to \(\mathcal X\).
The extension is unique by density and is denoted by
\(\mathfrak R_{q,0}(z)\).
A finite-energy solution for $q$ means an element of
;$\mathcal D_q$ satisfying \eqref{eq:system} with
$\varepsilon=0$ in $\mathcal D'(\mathbb R^N)$ it is admissible and has zero residual.
At \(\varepsilon=0\), Palais--Smale terminology refers to
the energy-and-residual condition
\eqref{eq:singular-ps-hypothesis}, without assuming
Fr\'echet differentiability of \(J_q\) or invoking a
weak-slope formulation.

The same definitions apply to auxiliary potential pairs
satisfying the uniform bounds in
Lemma~\ref{lem:fixed-local-recovery}.
Whenever \(\varepsilon=0\), the state is assumed to belong
to the corresponding residual domain.

\begin{lemma}[Radial identities and uniform bounds]
\label{lem:singular-radial-testing}
Let \(0\le\varepsilon\le1/2\).
For every \(z\in\mathcal X\) if \(\varepsilon>0\), and
every admissible \(z\) if \(\varepsilon=0\),
\begin{align}
 \langle\mathfrak R_{q,\varepsilon}(z),z\rangle
 &=S_q(z)-L_\varepsilon(z)-\frac2pR_\varepsilon(z),
 \label{eq:singular-radial-residual}\\
 J_{q,\varepsilon}(z)
 -\frac1p\langle\mathfrak R_{q,\varepsilon}(z),z\rangle
 &=\frac2{p^2}R_\varepsilon(z).
 \label{eq:singular-radial-energy}
\end{align}
Here \(L_0=L\), \(R_0=M\), and \(J_{q,0}=J_q\).
The bounds in Lemma~\ref{thm:uniform-ps-bound} remain valid
with \(J_{q,\varepsilon}'\) replaced by
\(\mathfrak R_{q,\varepsilon}\), including at
\(\varepsilon=0\).
The constants may be chosen uniformly for potentials satisfying
$0<\underline V\le a,b\le\overline V<\infty$.
\end{lemma}

\begin{proof}
For \(\varepsilon>0\), the identities follow from
Proposition~\ref{prop:energy-regularity}.
It remains to justify radial testing at \(\varepsilon=0\).

Choose \(\chi_R\in C_c^\infty(\mathbb R^N)\) such that
\[
 0\le\chi_R\le1,\qquad
 \chi_R=1\text{ on }B_R,\qquad
 \chi_R=0\text{ outside }B_{2R},\qquad
 |\nabla\chi_R|\le C/R.
\]
Since \(D_\xi\ell_0(z)\in L^{p'}_{\mathrm{loc}}\), density
permits testing \eqref{eq:singular-residual-definition}
with \(\chi_Rz\). Writing \(z=(u,v)\), we obtain
\begin{align*}
 \langle\mathfrak R_{q,0}(z),\chi_Rz\rangle
 ={}&
 \int_{\mathbb R^N}\chi_R
 \bigl(
 |\nabla u|^p+a|u|^p+|\nabla v|^p+b|v|^p
 \bigr)\,dx\\
 &+\int_{\mathbb R^N}
 \bigl(
 u|\nabla u|^{p-2}\nabla u
 +v|\nabla v|^{p-2}\nabla v
 \bigr)\cdot\nabla\chi_R\,dx\\
 &-\frac1p\int_{\mathbb R^N}
 \chi_R D_\xi\ell_0(z)\cdot z\,dx.
\end{align*}
The boundary terms vanish as \(R\to\infty\); for example,
\[
 \left|
 \int_{\mathbb R^N}
 u|\nabla u|^{p-2}\nabla u\cdot\nabla\chi_R\,dx
 \right|
 \le\frac CR
 \|\nabla u\|_{L^p(B_{2R}\setminus B_R)}^{p-1}
 \|u\|_{L^p(B_{2R}\setminus B_R)}
 \longrightarrow0.
\]
Moreover,
\[
 D_\xi\ell_0(z)\cdot z=p\ell_0(z)+2m(z),
 \qquad
 |D_\xi\ell_0(z)\cdot z|
 \le p\bigl(\ell_0^+(z)+\ell_0^-(z)\bigr)+2m(z)
 \in L^1.
\]
Since \(\chi_Rz\to z\) in \(\mathcal X\), continuity of
the residual and dominated convergence give
\[
 \langle\mathfrak R_{q,0}(z),z\rangle
 =S_q(z)-L_0(z)-\frac2pM(z).
\]
Subtracting \(1/p\) times this identity from
\(J_q(z)=(S_q(z)-L_0(z))/p\) proves
\eqref{eq:singular-radial-energy} at \(\varepsilon=0\).

With the radial identities established, the proof of
Lemma~\ref{thm:uniform-ps-bound} applies verbatim, replacing
\(J_{q,\varepsilon}'\) by \(\mathfrak R_{q,\varepsilon}\),
and now includes \(\varepsilon=0\).
The constants are uniform under the stated potential bounds,
since the Sobolev embedding constants are uniform.
\end{proof}

\begin{lemma}[Recovery in a moving frame]
\label{lem:singular-moving-recovery}
Let \(q_j=(a_j,b_j)\) satisfy
\[
 0<\underline V\le a_j,b_j\le\overline V<\infty
 \quad\text{almost everywhere, uniformly in }j,
\]
and suppose that
\(q_j\to q_*\) in \(L^\infty_{\mathrm{loc}}(\mathbb R^N)^2\).
Let \(0\le\varepsilon_j\le1/2\), \(\varepsilon_j\to0\),
and let \(z_j\) be admissible for
\(\mathfrak R_{q_j,\varepsilon_j}\), with
\[
 z_j\rightharpoonup Z\quad\text{in }\mathcal X,
 \qquad
 \sup_j L_{\varepsilon_j}^-(z_j)<\infty,
 \qquad
 \|\mathfrak R_{q_j,\varepsilon_j}(z_j)\|_{\mathcal X^*}
 \longrightarrow0.
\]
Then, after passage to a subsequence,
\begin{equation}
 z_j\to Z
 \quad\text{in }W^{1,p}_{\mathrm{loc}}(\mathbb R^N;\mathbb R^2),
 \qquad
 Z\in\mathcal D_{q_*},
 \qquad
 \mathfrak R_{q_*,0}(Z)=0
 \quad\text{in }\mathcal X^*.
 \label{eq:singular-moving-recovery}
\end{equation}
Moreover, 
\begin{equation}
 S_{q_*}(Z)=L_0(Z)+\frac2pM(Z),
 \qquad
 J_{q_*}(Z)=\frac2{p^2}M(Z).
 \label{eq:singular-profile-radial-identity}
\end{equation}
\end{lemma}

\begin{proof}
By local compactness and a diagonal argument, after passage
to a subsequence,
\[
 z_j\to Z\quad\text{almost everywhere and in }
 L^r_{\mathrm{loc}}(\mathbb R^N;\mathbb R^2)
 \quad(1\le r<\infty).
\]
The proof of Lemma~\ref{lem:local-nonlinear-convergence}
also allows zero parameters, since its estimates are uniform
for \(0\le\varepsilon\le1/2\). Hence
\[
 D_\xi\ell_{\varepsilon_j}(z_j)\to D_\xi\ell_0(Z)
 \quad\text{in }L^{p'}_{\mathrm{loc}},
\]
and
\[
 \ell_{\varepsilon_j}^\pm(z_j)\to\ell_0^\pm(Z),
 \qquad
 r_{\varepsilon_j}(z_j)\to m(Z)
 \quad\text{in }L^1_{\mathrm{loc}}.
\]
Fatou's lemma gives
\[
 L_0^-(Z)
 \le\liminf_{j\to\infty}L_{\varepsilon_j}^-(z_j)<\infty,
\]
so \(Z\in\mathcal D_{q_*}\).

Write \(z_j=(u_j,v_j)\) and \(Z=(Z_u,Z_v)\).
Fix \(0\le\chi\le1\) in \(C_c^\infty(\mathbb R^N)\).
The test pairs \(h_j=(\chi(u_j-Z_u),0)\) are bounded
in \(\mathcal X\), and therefore
\[
 \bigl|\langle
 \mathfrak R_{q_j,\varepsilon_j}(z_j),h_j
 \rangle\bigr|
 \le
 \|\mathfrak R_{q_j,\varepsilon_j}(z_j)\|_{\mathcal X^*}
 \|h_j\|_{\bar q}
 \longrightarrow0.
\]
Let \(K=\operatorname{supp}\chi\).
The lower-order and cutoff terms satisfy
\begin{align*}
 \left|\int a_j|u_j|^{p-2}u_j\chi(u_j-Z_u)\,dx\right|
 &\le C\|u_j\|_{L^p(K)}^{p-1}
       \|u_j-Z_u\|_{L^p(K)}\to0,\\
 \left|\int
 \partial_s\ell_{\varepsilon_j}(z_j)\chi(u_j-Z_u)\,dx\right|
 &\le C\|D_\xi\ell_{\varepsilon_j}(z_j)\|_{L^{p'}(K)}
       \|u_j-Z_u\|_{L^p(K)}\to0,\\
 \left|\int
 |\nabla u_j|^{p-2}\nabla u_j\cdot
 (u_j-Z_u)\nabla\chi\,dx\right|
 &\le C\|\nabla u_j\|_{L^p(K)}^{p-1}
       \|u_j-Z_u\|_{L^p(K)}\to0.
\end{align*}
Expanding the residual pairing thus gives
\[
 \int\chi|\nabla u_j|^{p-2}\nabla u_j
 \cdot(\nabla u_j-\nabla Z_u)\,dx=o(1).
\]
Weak convergence also yields
\[
 \int\chi|\nabla Z_u|^{p-2}\nabla Z_u
 \cdot(\nabla u_j-\nabla Z_u)\,dx=o(1).
\]
Subtracting and using monotonicity, we obtain
\begin{align*}
 2^{2-p}\int\chi|\nabla u_j-\nabla Z_u|^p\,dx
 &\le
 \int\chi
 \bigl(|\nabla u_j|^{p-2}\nabla u_j
       -|\nabla Z_u|^{p-2}\nabla Z_u\bigr)
 \cdot(\nabla u_j-\nabla Z_u)\,dx\\
 &=o(1).
\end{align*}
Choosing \(\chi=1\) on any prescribed compact set proves
local strong convergence of \(\nabla u_j\).
The same argument applies to \(v_j\), so
\[
 z_j\to Z
 \quad\text{in }W^{1,p}_{\mathrm{loc}}(\mathbb R^N;\mathbb R^2).
\]

For every \(h\in C_c^\infty(\mathbb R^N;\mathbb R^2)\),
the local convergences above and
\(q_j\to q_*\) in \(L^\infty_{\mathrm{loc}}\) give
\[
 \langle\mathfrak R_{q_*,0}(Z),h\rangle
 =\lim_{j\to\infty}
   \langle\mathfrak R_{q_j,\varepsilon_j}(z_j),h\rangle
 =0.
\]
Thus the limiting residual extends as zero to \(\mathcal X\).
Finally, Lemma~\ref{lem:singular-radial-testing} gives
\[
 S_{q_*}(Z)=L_0(Z)+\frac2pM(Z),
 \qquad
 J_{q_*}(Z)=\frac2{p^2}M(Z).
\]
\end{proof}

\subsection{Uniform localisation of the logarithmic equation}
\label{subsec:singular-cutoff-calculus}
Multiplication of both components by the same scalar admits
a derivative estimate uniform down to zero regularisation.
This estimate controls the slow cutoffs used in the main theorem,
including at zero components and for arbitrarily unequal
component sizes.

\begin{lemma}[Scalar modulation of the mixed logarithms]
\label{lem:singular-scalar-modulation}
There exists \(C_p>0\), depending only on \(p\), such that
the following statements hold for every \(\xi\in\mathbb R^2\)
and \(0\le\varepsilon\le1/2\).
\begin{enumerate}
\item[\textup{(i)}]
For \(0<t\le1\),
\begin{equation}
 |D_\xi\ell_\varepsilon(t\xi)
   -t^{p-1}D_\xi\ell_\varepsilon(\xi)|
 \le C_p t^{p-1}
 \bigl(|\log t|+1-t^2\bigr)|\xi|^{p-1}.
 \label{eq:singular-modulation-source}
\end{equation}
The defect vanishes at \(t=0\), by continuous extension,
and at \(t=1\).

\item[\textup{(ii)}]
For \(0<t\le1\),
\begin{align}
 \ell_\varepsilon^+(t\xi)
 &\le t^p\ell_\varepsilon^+(\xi),
 \label{eq:singular-modulation-positive}\\
 \ell_\varepsilon^-(t\xi)
 &\le t^p\ell_\varepsilon^-(\xi)
      +2t^p|\log t|m(\xi),
 \label{eq:singular-modulation-negative}\\
 r_\varepsilon(t\xi)
 &\le t^p r_\varepsilon(\xi),
 \qquad m(t\xi)=t^pm(\xi).
 \label{eq:singular-modulation-mass}
\end{align}
These estimates also hold at \(t=0\) by continuous extension.
\end{enumerate}
\end{lemma}

\begin{proof}
Assume first that \(\varepsilon>0\) and \(0<t\le1\).
For every \(s\in\mathbb R\),
\[
 t^2(s^2+\varepsilon)\le t^2s^2+\varepsilon
 \le s^2+\varepsilon,
\]
and hence
\begin{equation}
 2\log t
 \le \log(t^2s^2+\varepsilon)-\log(s^2+\varepsilon)
 \le0.
 \label{eq:singular-log-ratio}
\end{equation}

\emph{(i)}
In addition to \eqref{eq:singular-log-ratio}, we have
\[
 \left|
 \frac{t^2s^2}{t^2s^2+\varepsilon}
 -\frac{s^2}{s^2+\varepsilon}
 \right|
 =
 \frac{(1-t^2)\varepsilon s^2}
 {(t^2s^2+\varepsilon)(s^2+\varepsilon)}
 \le1-t^2.
\]
Using \(h_\varepsilon(s)=s^2\log(s^2+\varepsilon)\)
and its derivative, we obtain
\begin{align*}
 |h_\varepsilon(ts)-t^2h_\varepsilon(s)|
 &\le 2t^2|s|^2|\log t|,\\
 |h_\varepsilon'(ts)-t h_\varepsilon'(s)|
 &\le 2t|s|\bigl(2|\log t|+1-t^2\bigr).
\end{align*}
Applying these bounds to the two mixed terms in
\(\ell_\varepsilon\), together with their symmetric counterparts, proves
\eqref{eq:singular-modulation-source}.

\emph{(ii)}
By \eqref{eq:singular-log-ratio}, scaling can only decrease
the positive part of each logarithm and increases its negative
part by at most \(2|\log t|\).
Since the mixed weights scale by \(t^p\), this gives
\eqref{eq:singular-modulation-positive} and
\eqref{eq:singular-modulation-negative}.
The inequality \(t^2(s^2+\varepsilon)\le t^2s^2+\varepsilon\)
gives \(r_\varepsilon(t\xi)\le t^p r_\varepsilon(\xi)\)
term by term, while \(m(t\xi)=t^pm(\xi)\) by homogeneity.

Passing to the limit as \(\varepsilon\downarrow0\) gives
the same estimates at \(\varepsilon=0\).
The assertions at \(t=0\) follow by continuity, and the
source defect vanishes identically at \(t=1\).
\end{proof}

\begin{proposition}[Uniform residual commutator estimate]
\label{prop:singular-residual-commutator}
Fix \(B>0\) and \(0\le\varepsilon\le1/2\).
Let \(z\) be admissible for \(\mathfrak R_{q,\varepsilon}\),
with \(\|z\|_{\bar q}\le B\), and let
\(\vartheta\in W^{1,\infty}(\mathbb R^N)\) satisfy
\(0\le\vartheta\le1\) and \(\|\nabla\vartheta\|_\infty\le1\).
Then \(\vartheta z\) is also admissible, and
\begin{equation}
 \begin{aligned}
 &\bigl\|
 \mathfrak R_{q,\varepsilon}(\vartheta z)
 -\vartheta^{p-1}\mathfrak R_{q,\varepsilon}(z)
 \bigr\|_{\mathcal X^*}\\
 &\qquad\le C_B\left(
 \|\nabla\vartheta\|_\infty
 +\|z\|_{L^p(\{0<\vartheta<1\})}^{p-1}
 \right).
 \end{aligned}
 \label{eq:singular-residual-commutator}
\end{equation}
Here \(\vartheta^{p-1}f\) is defined by duality:
\(\langle\vartheta^{p-1}f,h\rangle
=\langle f,\vartheta^{p-1}h\rangle\).
The constants are uniform in \(\varepsilon\) and \(\vartheta\),
and for potentials satisfying fixed positive lower and upper bounds.
\end{proposition}
\begin{proof}
At \(\varepsilon=0\),
Lemma~\ref{lem:singular-scalar-modulation} gives
\[
 L_0^-(\vartheta z)
 \le L_0^-(z)+2\sup_{0<t\le1}t^p|\log t|\,M(z)<\infty.
\]
Since \(\vartheta z\in\mathcal X\), it follows that
\(\vartheta z\in\mathcal D_q\).

We first estimate the residual on compactly supported smooth
test pairs. For the scalar component \(u\) and test function
\(\phi\), the potential terms cancel, and the principal
difference is
\begin{equation}
\label{eq:principal-commutator-pairing}
\begin{aligned}
 &\bigl\langle
   -\Delta_p(\vartheta u)-\vartheta^{p-1}(-\Delta_p u),\phi
 \bigr\rangle\\
 &\quad =
 \int_{\mathbb R^N}
 \Bigl(
   |\nabla(\vartheta u)|^{p-2}\nabla(\vartheta u)
   -\vartheta^{p-1}|\nabla u|^{p-2}\nabla u
 \Bigr)\cdot\nabla\phi\,dx\\
 &\qquad {}-(p-1)\int_{\mathbb R^N}
 \vartheta^{p-2}\phi\,|\nabla u|^{p-2}
 \nabla u\cdot\nabla\vartheta\,dx.
\end{aligned}
\end{equation}
Applying the elementary estimate
\[
 \bigl\||A|^{p-2}A-|B|^{p-2}B\bigr\|_{L^{p'}}
 \le C_p\bigl(\|A\|_{L^p}+\|B\|_{L^p}\bigr)^{p-2}
 \|A-B\|_{L^p}
\]
with \(A=\nabla(\vartheta u)\) and \(B=\vartheta\nabla u\),
we have \(A-B=u\nabla\vartheta\). Thus the first integral
is bounded in absolute value by
\(C_B\|\nabla\vartheta\|_\infty\|\nabla\phi\|_{L^p}\). On the other hand,
H\"older's inequality bounds the second integral of \eqref{eq:principal-commutator-pairing} by
\(C_B\|\nabla\vartheta\|_\infty\|\phi\|_{L^p}\).
The same estimates hold for the second component \(v\).

By \eqref{eq:singular-modulation-source}, the source defect
vanishes wherever \(\vartheta=0\) or \(\vartheta=1\), and
\[
 \bigl\|D_\xi\ell_\varepsilon(\vartheta z)
 -\vartheta^{p-1}D_\xi\ell_\varepsilon(z)\bigr\|_{L^{p'}}
 \le C_p\|z\|_{L^p(\{0<\vartheta<1\})}^{p-1}.
\]
This global \(L^{p'}\) bound is needed only for the difference.
Combining it with the principal estimates gives
\eqref{eq:singular-residual-commutator} on smooth test pairs.

Finally, the product rule and the cutoff bounds give
\(\|\vartheta^{p-1}h\|_{\bar q}\le C\|h\|_{\bar q}\).
Thus \(\vartheta^{p-1}\mathfrak R_{q,\varepsilon}(z)\)
belongs to \(\mathcal X^*\), and the commutator extends to
\(\mathcal X^*\) by density. Their sum gives the required
extension of \(\mathfrak R_{q,\varepsilon}(\vartheta z)\),
proving admissibility also at \(\varepsilon=0\), which completes the proof.
\end{proof}

\begin{lemma}[Vanishing in the residual formulation]
\label{lem:singular-residual-vanishing}
Let \(q_j=(a_j,b_j)\) satisfy
\(0<\underline V\le a_j,b_j\le\overline V<\infty\)
almost everywhere, uniformly in \(j\).
Let \(0\le\varepsilon_j\le1/2\), and suppose that
\((w_j)\) is bounded in \(\mathcal X\), with each \(w_j\)
admissible for \(\mathfrak R_{q_j,\varepsilon_j}\).
If
\[
 \|\mathfrak R_{q_j,\varepsilon_j}(w_j)\|_{\mathcal X^*}\to0,
 \qquad
 \sup_{y\in\mathbb R^N}\int_{B_1(y)}|w_j|^p\,dx\to0,
\]
then
\begin{equation}
 \|w_j\|_{\bar q}\to0,\qquad
 L_{\varepsilon_j}^+(w_j)+L_{\varepsilon_j}^-(w_j)\to0,
 \qquad
 J_{q_j,\varepsilon_j}(w_j)\to0.
 \label{eq:singular-vanishing-conclusion}
\end{equation}
No convergence of \((\varepsilon_j)\) is required.
\end{lemma}

\begin{proof}
The cube argument in Lemma~\ref{lem:uniform-ps-vanishing}
gives \(w_j\to0\) in \(L^r(\mathbb R^N;\mathbb R^2)\)
for every \(p<r<\infty\).
Since \((w_j)\) is bounded in \(\mathcal X\),
\[
 \langle\mathfrak R_{q_j,\varepsilon_j}(w_j),w_j\rangle=o(1).
\]
Fix \(0<\sigma<1\). By \eqref{eq:singular-radial-residual}
and \eqref{eq:uniform-radial-source-bound},
\[
 S_{q_j}(w_j)
 =\int_{\mathbb R^N}
 \left(\ell_{\varepsilon_j}(w_j)
       +\frac2p r_{\varepsilon_j}(w_j)\right)\,dx+o(1)
 \le C_\sigma\|w_j\|_{L^{p+\sigma}}^{p+\sigma}+o(1)
 \longrightarrow0.
\]
Uniform equivalence of the potential norms now gives
\(\|w_j\|_{\bar q}\to0\).
By \eqref{eq:mass-control}, \(R_{\varepsilon_j}(w_j)\to0\);
hence \eqref{eq:singular-radial-energy} yields
\(J_{q_j,\varepsilon_j}(w_j)\to0\).

Finally, \eqref{eq:positive-growth} gives
\(L_{\varepsilon_j}^+(w_j)
\le C_\sigma\|w_j\|_{L^{p+\sigma}}^{p+\sigma}\to0\).
The energy identity then gives
\[
 L_{\varepsilon_j}^-(w_j)
 =pJ_{q_j,\varepsilon_j}(w_j)-S_{q_j}(w_j)
  +L_{\varepsilon_j}^+(w_j)
 \longrightarrow0.
\]
\end{proof}

\subsection{Proof of Theorem~\ref{thm:singular-global-decomposition}}
\label{subsec:singular-profile-decomposition}
\begin{proof} The process is divided into several steps.

   \emph{Step 1.}
Lemma~\ref{lem:singular-radial-testing} gives uniform bounds
on \(\|z_j\|_{\bar q}\) and \(L_{\varepsilon_j}^\pm(z_j)\).
After passage to a subsequence, \(z_j\rightharpoonup z^0\)
in \(\mathcal X\). Applying
Lemma~\ref{lem:singular-moving-recovery} with \(q_j=q\),
we find that \(z^0\) is a finite-energy solution for \(q\) and
\[
 z_j\to z^0
 \quad\text{in }W^{1,p}_{\mathrm{loc}}(\mathbb R^N;\mathbb R^2).
\]
Since \((z_j)\) is bounded and
\(\mathfrak R_{q,\varepsilon_j}(z_j)\to0\) in \(\mathcal X^*\),
we have
\(\langle\mathfrak R_{q,\varepsilon_j}(z_j),z_j\rangle=o(1)\).
The radial identity \eqref{eq:singular-radial-energy} therefore gives
\begin{equation}
 R_{\varepsilon_j}(z_j)\longrightarrow\frac{p^2}{2}c,
 \qquad c\ge0.
 \label{eq:singular-total-radial-mass}
\end{equation}
Set \(e_j^0=z_j-z^0\). Then \(e_j^0\to0\) locally in \(W^{1,p}\).
In the profile extraction below, we use the residual equation
for the original sequence \(z_j\).

\emph{Step 2.}
Suppose that \(i\) profiles have been extracted at pairwise
separating centres, and write
\[
 e_j^i=z_j-z^0-\sum_{h=1}^i\mathsf U_{y_j^h}Z^h.
\]
By induction, \((e_j^i)\) is bounded in \(\mathcal X\) and
converges strongly to zero locally in \(W^{1,p}\), both in
the original frame and in every previously extracted frame.
If
\begin{equation}
 \sup_{y\in\mathbb R^N}
 \int_{B_1(y)}|e_j^i|^p\,dx\longrightarrow0,
 \label{eq:singular-geometric-vanishing}
\end{equation}
we stop. Otherwise, after passage to a subsequence, there exist
centres \(y_j^{i+1}\) and \(\eta>0\) such that
\[
 \int_{B_1(y_j^{i+1})}|e_j^i|^p\,dx\ge\eta.
\]
As in the proof of Theorem~\ref{thm:fixed-global-decomposition},
local convergence to zero in the original and previously
extracted frames implies
\[
 |y_j^{i+1}|\to\infty,
 \qquad
 |y_j^{i+1}-y_j^h|\to\infty
 \quad(1\le h\le i).
\]
Consider the original sequence in the new frame,
\(\mathsf U_{-y_j^{i+1}}z_j\).
Since \(|y_j^{i+1}|\to\infty\), the translated potentials
\(q(\cdot+y_j^{i+1})\) converge to \(\bar q\) in
\(L^\infty_{\mathrm{loc}}\) by \eqref{hyp:continuous}.
A change of variables gives
\[
 \mathfrak R_{q(\cdot+y_j^{i+1}),\varepsilon_j}
 \bigl(\mathsf U_{-y_j^{i+1}}z_j\bigr)
 =\mathsf U_{-y_j^{i+1}}\mathfrak R_{q,\varepsilon_j}(z_j)
 \longrightarrow0
 \quad\text{in }\mathcal X^*,
\]
where the convergence follows from translation invariance
of the dual norm. Translations also preserve the reference
norm and both logarithmic energies. Hence, after passage
to a subsequence, the translated sequence converges weakly
in \(\mathcal X\) to some \(Z^{i+1}\).
Lemma~\ref{lem:singular-moving-recovery} shows that
\(Z^{i+1}\) is a finite-energy solution for \(\bar q\) and that
\[
 \mathsf U_{-y_j^{i+1}}z_j\to Z^{i+1}
 \quad\text{strongly in }W^{1,p}_{\mathrm{loc}}.
\]
By separation of the centres, \(z^0\) and all previously
extracted profiles vanish locally in the new frame.
Thus \(\mathsf U_{-y_j^{i+1}}e_j^i\to Z^{i+1}\) in \(L^p(B_1)\),
and the local mass bound gives
\[
 \int_{B_1}|Z^{i+1}|^p\,dx
 =\lim_{j\to\infty}\int_{B_1(y_j^{i+1})}|e_j^i|^p\,dx
 \ge\eta.
\]
In particular, \(Z^{i+1}\ne0\). 

As in the proof of Theorem~\ref{thm:fixed-global-decomposition},
the new remainder
\(e_j^{i+1}=e_j^i-\mathsf U_{y_j^{i+1}}Z^{i+1}\)
converges to zero in \(W^{1,p}_{\mathrm{loc}}\) in the original
frame and in all frames chosen so far. This closes the induction.
For each finite step, boundedness follows from
\[
\sup_j \|e_j^i\|_{\bar{q}}
\le \sup_j \|z_j\|_{\bar{q}} + \|z^0\|_{\bar{q}}
+ \sum_{h=1}^i \|Z^h\|_{\bar{q}}.
\]
This provides all the boundedness required during extraction. 

For any finite collection of extracted profiles and any
fixed $R$, the balls $B_R(y_j^h)$, including $y_j^0=0$,
are disjoint for large $j$. Local convergence of the radial
density and its nonnegativity give
\[
 \sum_{h=0}^i\int_{B_R}m(Z^h) dx
          \le\lim_jR_{\epsilon_j}(z_j)=\frac{p^2}{2}c.
\]
Letting $R\to\infty$ and using
\eqref{eq:singular-profile-radial-identity}, we obtain
\begin{equation}
 J_q(z^0)+\sum_{h=1}^iJ_{\bar q}(Z^h)\le c.
 \label{eq:singular-extraction-budget}
\end{equation}
A nonzero profile has positive mixed mass, since otherwise
its radial identity would give $S_{\bar q}=0$.
It therefore belongs to $\mathcal N_{\bar q,0}$ and has
energy at least $d_{\bar q,0}>0$, by
Corollary~\ref{cor:uniform-levels}. Thus the procedure must
terminate after a finite number $k$ of profiles, and \(e_j=e_j^k\) is used in Step~4.

\emph{Step 3.}
For each $0\le i\le k$, the original sequence in the frame
$y_j^i$ converges locally strongly to $Z^i$ in $W^{1,p}$;
its logarithmic densities converge locally in $L^1$.
Since there are finitely many frames, a diagonal choice
gives $R_j\to\infty$ such that
\begin{equation}
4R_j<\min_{i\ne h}\left|y_j^i-y_j^h\right|,
\quad\text{when } k\ge1,
\label{eq:singular-cutoff-separation}
\end{equation}
We first fix the radius.
For each integer \(n\ge1\), the local convergence from
Lemma~\ref{lem:singular-moving-recovery} gives
\[
\begin{aligned}
 &\|\mathsf U_{-y_j^i}z_j-Z^i\|_{W^{1,p}(B_{2n};\mathbb R^2)}
 +\sum_{\pm}
 \bigl\|\ell_{\varepsilon_j}^\pm(\mathsf U_{-y_j^i}z_j)
       -\ell_0^\pm(Z^i)\bigr\|_{L^1(B_{2n})}
 \longrightarrow0
\end{aligned}
\]
in each frame \(i=0,\ldots,k\). Since there are only finitely
many frames, we can choose \(j_n>j_{n-1}\) (with \(j_0=0\))
so that, whenever \(j\ge j_n\), the error on \(B_{2n}\) is at
most \(1/n\) in every frame. By separation of the centres, we
may also require that \(|y_j^i-y_j^h|>8n\) for all \(i\ne h\)
whenever \(j\ge j_n\).

We now let the radius increase in stages. Set
\[
 R_j=n \qquad\text{for }j_n\le j<j_{n+1}.
\]
On this entire range, \(B_{2R_j}=B_{2n}\), so the error on
\(B_{2R_j}\) is at most \(1/n\). As \(j\to\infty\), the
corresponding \(n\) tends to infinity. Thus \(R_j\to\infty\)
and, for every \(0\le i\le k\),
\begin{equation}
 \begin{aligned}
 &\|\mathsf U_{-y_j^i}z_j-Z^i\|_{W^{1,p}(B_{2R_j};\mathbb R^2)}\\
 &\quad+
 \sum_{\pm}
 \bigl\|\ell_{\varepsilon_j}^\pm(\mathsf U_{-y_j^i}z_j)
       -\ell_0^\pm(Z^i)\bigr\|_{L^1(B_{2R_j})}
 \longrightarrow0.
 \end{aligned}
 \label{eq:singular-expanding-ball-convergence}
\end{equation}
The same choice gives \(|y_j^i-y_j^h|>8R_j\) for \(i\ne h\),
so the balls \(B_{2R_j}(y_j^i)\) are pairwise disjoint.

Let
\[
 D_j=B_{2R_j}\setminus B_{R_j},
 \qquad
 A_j=\bigcup_{i=0}^k(y_j^i+D_j).
\]
Fix \(0\le i\le k\). Since \(D_j\) lies inside \(B_{2R_j}\)
and outside \(B_{R_j}\), we have
\[
 \|\mathsf U_{-y_j^i}z_j\|_{W^{1,p}(D_j)}
 \le
 \|\mathsf U_{-y_j^i}z_j-Z^i\|_{W^{1,p}(B_{2R_j})}
 +\|Z^i\|_{W^{1,p}(\mathbb R^N\setminus B_{R_j})}
 \longrightarrow0.
\]
The first term tends to zero by
\eqref{eq:singular-expanding-ball-convergence};
the second tends to zero since \(Z^i\in\mathcal X\) and
\(R_j\to\infty\). Likewise, for each sign,
\[
\begin{aligned}
 \int_{D_j}\ell_{\varepsilon_j}^{\pm}
       (\mathsf U_{-y_j^i}z_j)\,dx
 &\le
 \bigl\|\ell_{\varepsilon_j}^{\pm}(\mathsf U_{-y_j^i}z_j)
       -\ell_0^{\pm}(Z^i)\bigr\|_{L^1(B_{2R_j})}\\
 &\quad+
 \int_{\mathbb R^N\setminus B_{R_j}}\ell_0^{\pm}(Z^i)\,dx
 \longrightarrow0
\end{aligned}
\]
by \eqref{eq:singular-expanding-ball-convergence} and the
integrability of \(\ell_0^{\pm}(Z^i)\).
Translating back and summing over the finitely many annuli gives
\begin{equation}
 \int_{A_j}\Bigl(
 |\nabla u_j|^p+|\nabla v_j|^p+|z_j|^p
 +\ell_{\varepsilon_j}^+(z_j)+\ell_{\varepsilon_j}^-(z_j)
 \Bigr)\,dx\longrightarrow0.
 \label{eq:singular-small-annuli}
\end{equation}
Let \(\chi\in C_c^\infty(\mathbb R^N)\) satisfy
\(0\le\chi\le1\), \(\chi=1\) on \(B_1\), and \(\chi=0\)
outside \(B_2\). Set \(\chi_R(x)=\chi(x/R)\) and define
\begin{equation}
 \widehat Z_j^i
 :=\chi_{R_j}\mathsf U_{-y_j^i}z_j
 \quad(0\le i\le k),
 \qquad
 w_j:=\left(1-\sum_{i=0}^k
       \mathsf U_{y_j^i}\chi_{R_j}\right)z_j.
 \label{eq:singular-localised-pieces}
\end{equation}
By the product rule and
\(\|\nabla\chi_{R_j}\|_\infty\le C/R_j\),
\[
 \|\widehat Z_j^i-Z^i\|_{\bar q}
 \le C\|\mathsf U_{-y_j^i}z_j-Z^i\|_{W^{1,p}(B_{2R_j})}
    +\|\chi_{R_j}Z^i-Z^i\|_{\bar q}.
\]
The first term tends to zero by
\eqref{eq:singular-expanding-ball-convergence}, and the second
tends to zero since \(Z^i\in\mathcal X\) and \(R_j\to\infty\).
Hence
\begin{equation}
 \widehat Z_j^i\longrightarrow Z^i
 \quad\text{in }\mathcal X
 \qquad(0\le i\le k).
 \label{eq:singular-cut-profile-strong}
\end{equation}
Fix \(0\le i\le k\). On \(B_{R_j}\), the cutoff equals one,
so \(\widehat Z_j^i=\mathsf U_{-y_j^i}z_j\). Hence
\eqref{eq:singular-expanding-ball-convergence} gives
\[
 \bigl\|\ell_{\varepsilon_j}^{\pm}(\widehat Z_j^i)
       -\ell_0^{\pm}(Z^i)\bigr\|_{L^1(B_{R_j})}
 \longrightarrow0.
\]
On the transition annulus \(D_j=B_{2R_j}\setminus B_{R_j}\),
the modulation estimates
\eqref{eq:singular-modulation-positive} and
\eqref{eq:singular-modulation-negative} give
\[
 \ell_{\varepsilon_j}^{\pm}(\widehat Z_j^i)
 \le C\Bigl(
 \ell_{\varepsilon_j}^{\pm}(\mathsf U_{-y_j^i}z_j)
 +m(\mathsf U_{-y_j^i}z_j)\Bigr),
\]
since \(t^p|\log t|\) is bounded for \(0<t\le1\).
As \(m(\xi)\le C|\xi|^p\) and \(y_j^i+D_j\subset A_j\),
a change of variables and \eqref{eq:singular-small-annuli}
show that
\[
 \int_{D_j}\ell_{\varepsilon_j}^{\pm}(\widehat Z_j^i)\,dx
 \longrightarrow0.
\]
Finally, \(\widehat Z_j^i=0\) outside \(B_{2R_j}\), while
\(\ell_0^\pm(Z^i)\in L^1(\mathbb R^N)\) implies
\[
 \int_{\mathbb R^N\setminus B_{R_j}}\ell_0^\pm(Z^i)\,dx
 \longrightarrow0.
\]
Combining the inner-ball estimate with these two tail bounds yields
\[
 \ell_{\varepsilon_j}^{\pm}(\widehat Z_j^i)
 \longrightarrow\ell_0^\pm(Z^i)
 \quad\text{in }L^1(\mathbb R^N).
\]
By \eqref{eq:singular-cut-profile-strong},
\(\widehat Z_j^i\to Z^i\) in \(L^p\).
Since \(|D_\xi r_\varepsilon(\xi)|\le C|\xi|^{p-1}\)
uniformly in \(\varepsilon\), H\"older's inequality gives
\[
 \|r_{\varepsilon_j}(\widehat Z_j^i)
   -r_{\varepsilon_j}(Z^i)\|_{L^1}
 \le C\bigl(\|\widehat Z_j^i\|_{L^p}+\|Z^i\|_{L^p}\bigr)^{p-1}
       \|\widehat Z_j^i-Z^i\|_{L^p}
 \longrightarrow0.
\]
Moreover, \(r_{\varepsilon_j}(Z^i)\to m(Z^i)\) pointwise and
\(0\le r_{\varepsilon_j}(Z^i)\le m(Z^i)\in L^1\).
Dominated convergence therefore yields
\(r_{\varepsilon_j}(\widehat Z_j^i)\to m(Z^i)\) in \(L^1\).
The same estimate with \(r_0=m\) gives
\(m(\widehat Z_j^i)\to m(Z^i)\) in \(L^1\).

Finally, strong convergence in \(\mathcal X\) gives
\(S_{\bar{q}}(\widehat Z_j^i)\to S_{\bar{q}}(Z^i)\).
Combining this with the \(L^1\) convergence of the positive
and negative logarithmic densities proved above, we obtain
\[
 J_{\bar{q},\varepsilon_j}(\widehat Z_j^i)
 \longrightarrow J_{\bar{q}}(Z^i).
\]

\emph{Step 4.}
Let \(\vartheta\) be either one of the cutoffs
\(\mathsf U_{y_j^i}\chi_{R_j}\), \(0\le i\le k\), or the
exterior cutoff
\[
 \vartheta_j=1-\sum_{i=0}^k\mathsf U_{y_j^i}\chi_{R_j}.
\]
Since the translated cutoffs have pairwise disjoint supports,
both choices satisfy
\[
 0\le\vartheta\le1,\qquad
 \{0<\vartheta<1\}\subset A_j,\qquad
 \|\nabla\vartheta\|_\infty\le C/R_j\le1
\]
for all sufficiently large \(j\).
By Step~1, \(\|z_j\|_{\bar q}\le B\) uniformly in \(j\).
Proposition~\ref{prop:singular-residual-commutator} therefore
gives
\[
\begin{aligned}
 &\bigl\|\mathfrak R_{q,\varepsilon_j}(\vartheta z_j)
 -\vartheta^{p-1}\mathfrak R_{q,\varepsilon_j}(z_j)
 \bigr\|_{\mathcal X^*}\\
 &\qquad\le C_B\left(
 R_j^{-1}
 +\left(\int_{A_j}|z_j|^p\,dx\right)^{(p-1)/p}
 \right)\longrightarrow0
\end{aligned}
\]
by \(R_j\to\infty\) and \eqref{eq:singular-small-annuli}.
Moreover, the uniform multiplier bound gives
\[
 \bigl\|\vartheta^{p-1}\mathfrak R_{q,\varepsilon_j}(z_j)
 \bigr\|_{\mathcal X^*}
 \le C\|\mathfrak R_{q,\varepsilon_j}(z_j)\|_{\mathcal X^*}
 \longrightarrow0.
\]
The triangle inequality now yields
\(\mathfrak R_{q,\varepsilon_j}(\vartheta z_j)\to0\)
in \(\mathcal X^*\). Taking the two choices of \(\vartheta\)
in turn, we obtain
\[
 \mathfrak R_{q,\varepsilon_j}
 \bigl(\mathsf U_{y_j^i}\widehat Z_j^i\bigr)\to0
 \quad(0\le i\le k),
 \qquad
 \mathfrak R_{q,\varepsilon_j}(w_j)\to0
 \quad\text{in }\mathcal X^*.
\]
For \(i\ge1\), the separation from the origin gives
\(|y_j^i|>8R_j\) for all sufficiently large \(j\).
Thus \(|y_j^i|-2R_j\to\infty\), and \eqref{hyp:continuous}
implies
\[
 \|q-\bar q\|_{L^\infty(B_{2R_j}(y_j^i))}\longrightarrow0.
\]
Write \(W_j=\mathsf U_{y_j^i}\widehat Z_j^i\), which is
supported in \(B_{2R_j}(y_j^i)\).
Only the potential terms differ between the residuals for
\(q\) and \(\bar q\). Hence, for every \(h\in\mathcal X\),
H\"older's inequality gives
\begin{align*}
 \bigl|\langle
 \mathfrak R_{q,\varepsilon_j}(W_j)
 -\mathfrak R_{\bar q,\varepsilon_j}(W_j),h
 \rangle\bigr|
 &\le C\|q-\bar q\|_{L^\infty(B_{2R_j}(y_j^i))}
       \int_{\mathbb R^N}|W_j|^{p-1}|h|\,dx\\
 &\le C\|q-\bar q\|_{L^\infty(B_{2R_j}(y_j^i))}
       \|W_j\|_{L^p}^{p-1}\|h\|_{\bar q}.
\end{align*}
Taking the supremum over \(\|h\|_{\bar q}\le1\), we obtain
\[
 \|\mathfrak R_{q,\varepsilon_j}(W_j)
   -\mathfrak R_{\bar q,\varepsilon_j}(W_j)\|_{\mathcal X^*}
 \le C\|q-\bar q\|_{L^\infty(B_{2R_j}(y_j^i))}
       \|W_j\|_{L^p}^{p-1}
 \longrightarrow0,
\]
since \(\|W_j\|_{L^p}=\|\widehat Z_j^i\|_{L^p}\) is bounded
by \eqref{eq:singular-cut-profile-strong}.

Combining this estimate with the residual bound for \(W_j\)
and using translation invariance for the constant potential
\(\bar q\), we obtain, for \(i\ge1\),
\[
 \|\mathfrak R_{\bar q,\varepsilon_j}(\widehat Z_j^i)\|_{\mathcal X^*}
 =
 \|\mathfrak R_{\bar q,\varepsilon_j}(W_j)\|_{\mathcal X^*}
 \longrightarrow0.
\]
For \(i=0\), the cutoff estimate already gives
\(\mathfrak R_{q,\varepsilon_j}(\widehat Z_j^0)\to0\)
in \(\mathcal X^*\). The potential comparison is uniform
in the regularization parameter.

By the exact partition, translation invariance of the reference
norm, and strong convergence of each \(\widehat Z_j^i\) to \(Z^i\),
\[
 \|w_j-e_j\|_{\bar q}
 =\left\|\sum_{i=0}^k
   \mathsf U_{y_j^i}(Z^i-\widehat Z_j^i)\right\|_{\bar q}
 \le \sum_{i=0}^k\|Z^i-\widehat Z_j^i\|_{\bar q}
 \longrightarrow 0.
\]
By the stopping criterion, the local \(L^p\) mass of \(e_j\)
vanishes uniformly in the centre. Hence
\[
 \sup_{y\in\mathbb R^N}\|w_j\|_{L^p(B_1(y))}
 \le \|w_j-e_j\|_{L^p(\mathbb R^N)}
   +\sup_{y\in\mathbb R^N}\|e_j\|_{L^p(B_1(y))}
 \longrightarrow 0.
\]
Moreover, \((w_j)\) is bounded in \(\mathcal X\), and the
cutoff estimate from Proposition~\ref{prop:singular-residual-commutator} gives
\[
 \|\mathfrak R_{q,\varepsilon_j}(w_j)\|_{\mathcal X^*}
 \longrightarrow 0.
\]
Applying the Lemma~\ref{lem:singular-residual-vanishing} for sequences with small residuals
to \(w_j\), we obtain
\[
 \|w_j\|_{\bar q}\to 0,
 \qquad
 L_{\varepsilon_j}^+(w_j)+L_{\varepsilon_j}^-(w_j)\to 0.
\]
Finally,
\[
 \|e_j\|_{\bar q}
 \le \|e_j-w_j\|_{\bar q}+\|w_j\|_{\bar q}
 \longrightarrow 0,
\]
which establishes \eqref{eq:singular-profile-decomposition}.

\emph{Step 5.}
Fix either sign. Outside \(A_j\), exactly one cutoff equals
one and the others vanish, so the logarithmic densities split
exactly. On \(A_j\), the transition annuli are pairwise disjoint.
Thus, if the cutoff on the annulus containing \(x\) has value
\(t\), the corresponding profile piece is \(tz_j(x)\) and the
remainder is \((1-t)z_j(x)\). Applying
Lemma~\ref{lem:singular-scalar-modulation}\textup{(ii)}
to both \(t\) and \(1-t\), we obtain
\[
 \bigl|\ell_{\varepsilon_j}^{\pm}(z_j)
       -\ell_{\varepsilon_j}^{\pm}(tz_j)
       -\ell_{\varepsilon_j}^{\pm}((1-t)z_j)\bigr|
 \le C\bigl(\ell_{\varepsilon_j}^{\pm}(z_j)+m(z_j)\bigr),
\]
since \(s^p|\log s|\) is bounded on \(0<s\le1\).
As \(m(z_j)\le C|z_j|^p\), integration gives
\[
\begin{aligned}
 &\left\|\ell_{\varepsilon_j}^{\pm}(z_j)
   -\sum_{i=0}^k\mathsf U_{y_j^i}
       \ell_{\varepsilon_j}^{\pm}(\widehat Z_j^i)
   -\ell_{\varepsilon_j}^{\pm}(w_j)\right\|_1\\
 &\qquad\le C\int_{A_j}
   \bigl(\ell_{\varepsilon_j}^{\pm}(z_j)+|z_j|^p\bigr)\,dx
   \longrightarrow0
\end{aligned}
\]
by \eqref{eq:singular-small-annuli}.

To replace the cutoff pieces by the limiting profiles, we add
and subtract their densities and the density of the remainder.
The triangle inequality and translation invariance of the
\(L^1\) norm yield
\[
\begin{aligned}
 &\left\|\ell_{\varepsilon_j}^{\pm}(z_j)
   -\sum_{i=0}^k\mathsf U_{y_j^i}\ell_0^\pm(Z^i)\right\|_1\\
 &\quad\le
 \left\|\ell_{\varepsilon_j}^{\pm}(z_j)
   -\sum_{i=0}^k\mathsf U_{y_j^i}
       \ell_{\varepsilon_j}^{\pm}(\widehat Z_j^i)
   -\ell_{\varepsilon_j}^{\pm}(w_j)\right\|_1\\
 &\qquad+
 \sum_{i=0}^k
 \left\|\ell_{\varepsilon_j}^{\pm}(\widehat Z_j^i)
       -\ell_0^\pm(Z^i)\right\|_1
 +L_{\varepsilon_j}^{\pm}(w_j).
\end{aligned}
\]
The first term tends to zero by the preceding estimate, the
finite sum by the \(L^1\) convergence proved in Step~3, and
the last term by Step~4. This proves
\eqref{eq:singular-entropy-splitting}. The same argument,
using \eqref{eq:singular-modulation-mass}, proves
\eqref{eq:singular-mass-splitting}.

We first establish the approximate splitting of the principal energy.
Outside the transition annuli $A_j$, exactly one cut-off coefficient
equals $1$, while all the others vanish. Hence the original energy
density agrees almost everywhere with the sum of the energy densities
of the truncated pieces. It therefore suffices to estimate the
discrepancy on $A_j$.

Fix one transition annulus and denote the corresponding cut-off
function by $\chi$. Since the transition annuli are pairwise
disjoint, only two pieces are involved there, namely
$\chi u_j$ and $(1-\chi)u_j$. By the product rule,
\[
\nabla(\chi u_j)
=\chi\nabla u_j+u_j\nabla\chi,
\qquad
\nabla((1-\chi)u_j)
=(1-\chi)\nabla u_j-u_j\nabla\chi.
\]
Using $0\leq\chi\leq1$ and
$|\nabla\chi|\leq C/R_j$, we obtain
\[
\begin{aligned}
\Big|
|\nabla u_j|^p
&-|\nabla(\chi u_j)|^p
-|\nabla((1-\chi)u_j)|^p
\Big|\\
&\leq C\bigl(
|\nabla u_j|^p+R_j^{-p}|u_j|^p
\bigr).
\end{aligned}
\]
The discrepancy in the potential term is bounded by
\[
a(x)|u_j|^p
\bigl|1-\chi^p-(1-\chi)^p\bigr|
\leq \|a\|_{L^\infty}|u_j|^p.
\]
The same estimates hold for the $v_j$-component. Summing over
all transition annuli and integrating, we thus obtain
\[
\begin{aligned}
&\Big|
S_q(z_j)
-\sum_{i=0}^k
S_q\bigl(\mathsf U_{y_j^i}\widehat Z_j^i\bigr)
-S_q(w_j)
\Big|\\
&\qquad\leq
C\int_{A_j}
\bigl(
|\nabla u_j|^p+
|\nabla v_j|^p+
|z_j|^p
\bigr)\,dx
+CR_j^{-p}\|z_j\|_{L^p(\mathbb R^N)^2}^p
\longrightarrow0.
\end{aligned}
\]
Indeed, the first term tends to zero by the previously
established vanishing of the energy on the transition annuli,
whereas the second tends to zero because $R_j\to\infty$ and
$(z_j)$ is bounded in $L^p(\mathbb R^N)^2$. Consequently,
\[
S_q(z_j)
=
\sum_{i=0}^k
S_q\bigl(\mathsf U_{y_j^i}\widehat Z_j^i\bigr)
+S_q(w_j)+o(1).
\]

We now pass to the limit in each term. Since $w_j\to0$ strongly
in $\mathcal X$, we have $S_q(w_j)\to0$. For $i=0$, recalling
that $y_j^0=0$ and $\widehat Z_j^0\to z^0$ strongly in
$\mathcal X$, the continuity of $S_q$ yields
\[
S_q(\widehat Z_j^0)\longrightarrow S_q(z^0).
\]

For $i\geq1$, we show that the corresponding limiting energy
is $S_{\bar q}(Z^i)$. Indeed, the translated truncated profile
is supported in $B_{2R_j}(y_j^i)$. Since
\[
|y_j^i|-2R_j\longrightarrow\infty
\]
and $q(x)\to\bar q$ as $|x|\to\infty$, we have
\[
\|q-\bar q\|_{L^\infty(B_{2R_j}(y_j^i))}
\longrightarrow0.
\]
Since $S_q$ and $S_{\bar q}$ differ only in their potential
terms, it follows that
\[
\begin{aligned}
&\left|
S_q(\mathsf U_{y_j^i}\widehat Z_j^i)
-S_{\bar q}(\mathsf U_{y_j^i}\widehat Z_j^i)
\right|\\
&\qquad\leq
C\|q-\bar q\|_{L^\infty(B_{2R_j}(y_j^i))}
\|\widehat Z_j^i\|_{L^p(\mathbb R^N)^2}^p
\longrightarrow0.
\end{aligned}
\]
Moreover, by translation invariance of $S_{\bar q}$ and the
strong convergence $\widehat Z_j^i\to Z^i$ in $\mathcal X$,
we obtain
\[
\begin{aligned}
S_q(\mathsf U_{y_j^i}\widehat Z_j^i)
&=S_{\bar q}(\mathsf U_{y_j^i}\widehat Z_j^i)+o(1)\\
&=S_{\bar q}(\widehat Z_j^i)+o(1)\\
&\longrightarrow S_{\bar q}(Z^i).
\end{aligned}
\]
Combining these limits with the preceding approximate energy
splitting, and recalling that the number of profiles is finite,
we conclude that
\[
S_q(z_j)
=
S_q(z^0)+
\sum_{i=1}^k S_{\bar q}(Z^i)+o(1),
\]
which proves (7.16). The definition of the energy then gives
\eqref{eq:singular-energy-splitting}.

Finally, the parameter-independent power-flux estimates of
Lemma~\ref{lem:fixed-nonlinear-splitting}, or the same
partition argument in $L^{p'}$, give
\begin{equation}
 \left\|S_q'(z_j)
         -\sum_{i=0}^k\mathsf U_{y_j^i}S_{\bar{q}}'(Z^i)
       \right\|_{\mathcal X^*}\longrightarrow0.
 \label{eq:singular-principal-splitting}
\end{equation}
Indeed, \[z_j=\sum_{i=0}^k \mathsf U_{y_j^i}\widehat{Z_j}^i+w_j,\quad \|w_j\|_{\mathcal{X}}\rightarrow 0.\]
On the other hand, Since $p\ge2$, the standard $p$-power inequality and
H\"older's inequality yield
\[
\|S_q'(f)-S_q'(g)\|_{\mathcal X^*}
\le C(\|f\|_{\mathcal X}+\|g\|_{\mathcal X})^{p-2}
\|f-g\|_{\mathcal X}
\le C\|f-g\|_{\mathcal X},
\]
where the last inequality follows from the uniform boundedness
of $f$ and $g$ in $\mathcal X$. And is thus, 
\[\|S_q'(z_j)-S_q'(\sum_{i=0}^k \mathsf U_{y_j^i}\widehat{Z_j}^i)\|_{\mathcal{X}^*}\le C\|w_j\|_{\mathcal{X}}\to 0.\]
Since the compact sets of $\mathsf U_{y_j^i}\widehat{Z_j}^i$ are mutually disjoint,
\begin{equation}
\begin{aligned}
S_q'(z_j)
&=\sum_{i=0}^k
S_q'(\mathsf U_{y_j^i}\widehat Z_j^i)
+o_{\mathcal X^*}(1)\\
&=\sum_{i=0}^k
\mathsf U_{y_j^i}S_{\bar{q}}'(\widehat Z_j^i)
+o_{\mathcal X^*}(1)\\
&=\sum_{i=0}^k
\mathsf U_{y_j^i}S_{\bar{q}}'(Z^i)
+o_{\mathcal X^*}(1).
\end{aligned}
\end{equation}
By \eqref{eq:singular-residual-definition},
\[
 D_\xi\ell_{\varepsilon_j}(z_j)
       =S_q'(z_j)-p\mathfrak R_{q,\varepsilon_j}(z_j),
 \qquad
 D_\xi\ell_0(Z^i)=S_{\bar{q}}'(Z^i)
\]
as elements of $\mathcal X^*$. The extensions here are
unique because compactly supported smooth test pairs are
dense. Equation
\eqref{eq:singular-source-splitting} follows.
Finally, the radial identities and
Corollary~\ref{cor:uniform-levels} give Nehari membership
and \eqref{eq:singular-profile-count}. We complete the proof.
\end{proof}
\begin{remark}[Meaning of derivative splitting]
\label{rem:singular-derivative-meaning}
Each limiting source in \eqref{eq:singular-source-splitting}
admits the canonical extension
$D_\xi\ell_0(Z^i)=S_{\bar{q}}'(Z^i)\in\mathcal X^*$.
Thus, the splitting holds in the dual norm without requiring
global Lebesgue integrability. The localised profiles also have
vanishing residuals at the original parameters. These conclusions
follow from the uniform modulation estimate
\eqref{eq:singular-modulation-source}, without requiring uniformity
in \eqref{eq:fixed-source-splitting}.
\end{remark}

\subsection{Recovery of finite configurations}
\label{subsec:singular-reverse-construction}
\begin{theorem}[Realisation of prescribed finite-energy profiles]
\label{thm:singular-profile-realisation}
Let $z^0$ be a finite-energy solution for $q$, and let
$Z^1,\ldots,Z^k$ be nonzero finite-energy solutions for
$\bar q$, where $k\geq0$. Fix centres satisfying
\eqref{eq:fixed-profile-separation} and a sequence
$0\leq\varepsilon_j\leq1/2$ with $\varepsilon_j\to0$.
Then there exist radii $R_j\to\infty$ such that
\begin{equation}
z_j=\chi_{R_j}z^0
+\sum_{i=1}^k\mathsf U_{y_j^i}(\chi_{R_j}Z^i)
\label{eq:singular-recovery-sequence}
\end{equation}
is admissible and satisfies \eqref{eq:singular-ps-hypothesis}
at the level
\begin{equation}
c=J_q(z^0)+\sum_{i=1}^kJ_{\bar q}(Z^i).
\label{eq:singular-recovery-level}
\end{equation}
All splittings in Theorem~\ref{thm:singular-global-decomposition}
hold with these prescribed profiles. Moreover,
$z_j\rightharpoonup z^0$ in $\mathcal X$, with strong convergence
if and only if $k=0$. The result also holds for
$\varepsilon_j\equiv0$.
\end{theorem}
\begin{proof}
Set $Z^0=z^0$. For each fixed profile,
Proposition~\ref{prop:singular-residual-commutator}
at zero regularisation yields
\begin{equation}
\begin{aligned}
\|\mathfrak R_{q,0}(\chi_R z^0)\|_{\mathcal X^*}
&\longrightarrow0,\\
\|\mathfrak R_{\bar q,0}(\chi_R Z^i)\|_{\mathcal X^*}
&\longrightarrow0
\qquad(1\le i\le k),
\end{aligned}
\qquad(R\to\infty),
\label{eq:singular-truncated-solution}
\end{equation}
since the profile's $L^p$ mass on
$B_{2R}\setminus B_R$ vanishes and
$\|\nabla\chi_R\|_{L^\infty}\to0$.
Moreover, $\chi_R Z^i\to Z^i$ strongly in $\mathcal X$.
The estimates \eqref{eq:singular-modulation-positive} and
\eqref{eq:singular-modulation-negative} provide an integrable
majorant for the logarithmic densities. Thus, dominated
convergence gives convergence in $L^1$, and hence
\begin{equation}
\begin{aligned}
J_q(\chi_R z^0)&\longrightarrow J_q(z^0),\\
J_{\bar q}(\chi_R Z^i)&\longrightarrow J_{\bar q}(Z^i)
\qquad(1\le i\le k).
\end{aligned}
\label{eq:singular-truncated-energy}
\end{equation}
For fixed $R$, all truncated profiles are supported in
$B_{2R}$. By Lemma~\ref{lem:local-nonlinear-convergence},
whose proof also covers $\varepsilon=0$, we have
\begin{align*}
\|D_\xi\ell_{\varepsilon_j}(\chi_R Z^i)
-D_\xi\ell_0(\chi_R Z^i)\|_{L^{p'}(B_{2R})}
&\longrightarrow0,\\
\|\ell_{\varepsilon_j}^{\pm}(\chi_R Z^i)
-\ell_0^{\pm}(\chi_R Z^i)\|_{L^1}
&\longrightarrow0.
\end{align*}
Choose $R_j\to\infty$ slowly enough for these errors to
tend to zero for all profiles, while also imposing
\eqref{eq:singular-cutoff-separation}.
Combining this choice
with \eqref{eq:singular-truncated-solution} and
\eqref{eq:singular-truncated-energy} yields
\[
\begin{aligned}
 \|\mathfrak R_{q,\varepsilon_j}(\chi_{R_j}z^0)\|_{\mathcal X^*}
 &\to0,
 &J_{q,\varepsilon_j}(\chi_{R_j}z^0)&\to J_q(z^0),\\
 \|\mathfrak R_{\bar q,\varepsilon_j}(\chi_{R_j}Z^i)\|_{\mathcal X^*}
 &\to0,
 &J_{\bar q,\varepsilon_j}(\chi_{R_j}Z^i)&\to J_{\bar q}(Z^i)
 \quad(1\le i\le k).
\end{aligned}
\]
Indeed, By \eqref{eq:singular-truncated-solution} and
\eqref{eq:singular-truncated-energy}, we can choose an increasing
sequence $R_n\to\infty$ such that
\[
\begin{aligned}
&\|\mathfrak R_{q,0}(\chi_{R_n}z^0)\|_{\mathcal X^*}
+
|J_q(\chi_{R_n}z^0)-J_q(z^0)|
\leq \frac1n,\\
&\|\mathfrak R_{\bar q,0}(\chi_{R_n}Z^i)\|_{\mathcal X^*}
+
|J_{\bar q}(\chi_{R_n}Z^i)-J_{\bar q}(Z^i)|
\leq \frac1n
\qquad(1\le i\le k).
\end{aligned}
\]
For each fixed $R_n$, the local convergence established above
implies that, for sufficiently large $j$,
\[
\begin{aligned}
&\|\mathfrak R_{q,\varepsilon_j}(\chi_{R_n}z^0)
-\mathfrak R_{q,0}(\chi_{R_n}z^0)\|_{\mathcal X^*}\\
&\quad+
|J_{q,\varepsilon_j}(\chi_{R_n}z^0)
-J_q(\chi_{R_n}z^0)|
\leq \frac1n,
\end{aligned}
\]
and
\[
\begin{aligned}
&\|\mathfrak R_{\bar q,\varepsilon_j}(\chi_{R_n}Z^i)
-\mathfrak R_{\bar q,0}(\chi_{R_n}Z^i)\|_{\mathcal X^*}\\
&\quad+
|J_{\bar q,\varepsilon_j}(\chi_{R_n}Z^i)
-J_{\bar q}(\chi_{R_n}Z^i)|
\leq \frac1n
\qquad(1\le i\le k)
\end{aligned}
\]
simultaneously for all profiles.

A diagonal choice of indices $j_n$ therefore allows us to set
\[
R_j=R_n,\qquad j_n\leq j<j_{n+1}.
\]
Then $R_j\to\infty$, and the triangle inequality gives
\[
\begin{aligned}
\|\mathfrak R_{q,\varepsilon_j}
(\chi_{R_j}z^0)\|_{\mathcal X^*}
&\leq
\|\mathfrak R_{q,\varepsilon_j}(\chi_{R_j}z^0)
-\mathfrak R_{q,0}(\chi_{R_j}z^0)\|_{\mathcal X^*}\\
&\quad+
\|\mathfrak R_{q,0}(\chi_{R_j}z^0)\|_{\mathcal X^*}
\longrightarrow0,
\end{aligned}
\]
and
\[
\begin{aligned}
\|\mathfrak R_{\bar q,\varepsilon_j}
(\chi_{R_j}Z^i)\|_{\mathcal X^*}
&\leq
\|\mathfrak R_{\bar q,\varepsilon_j}(\chi_{R_j}Z^i)
-\mathfrak R_{\bar q,0}(\chi_{R_j}Z^i)\|_{\mathcal X^*}\\
&\quad+
\|\mathfrak R_{\bar q,0}(\chi_{R_j}Z^i)\|_{\mathcal X^*}
\longrightarrow0
\qquad(1\le i\le k).
\end{aligned}
\]
Similarly,
\[
\begin{aligned}
|J_{q,\varepsilon_j}(\chi_{R_j}z^0)-J_q(z^0)|
&\leq
|J_{q,\varepsilon_j}(\chi_{R_j}z^0)
-J_q(\chi_{R_j}z^0)|\\
&\quad+
|J_q(\chi_{R_j}z^0)-J_q(z^0)|
\longrightarrow0,
\end{aligned}
\]
and
\[
\begin{aligned}
|J_{\bar q,\varepsilon_j}(\chi_{R_j}Z^i)
-J_{\bar q}(Z^i)|
&\leq
|J_{\bar q,\varepsilon_j}(\chi_{R_j}Z^i)
-J_{\bar q}(\chi_{R_j}Z^i)|\\
&\quad+
|J_{\bar q}(\chi_{R_j}Z^i)-J_{\bar q}(Z^i)|
\longrightarrow0
\qquad(1\le i\le k).
\end{aligned}
\]
This proves the desired assertions.

Since the translated supports are pairwise disjoint, the
locality of the energy and residual yields
\[
\begin{aligned}
J_{q,\varepsilon_j}(z_j)
&=J_{q,\varepsilon_j}(\chi_{R_j}z^0)
+\sum_{i=1}^k
J_{q,\varepsilon_j}
\bigl(\mathsf U_{y_j^i}(\chi_{R_j}Z^i)\bigr),\\
\mathfrak R_{q,\varepsilon_j}(z_j)
&=\mathfrak R_{q,\varepsilon_j}(\chi_{R_j}z^0)
+\sum_{i=1}^k
\mathfrak R_{q,\varepsilon_j}
\bigl(\mathsf U_{y_j^i}(\chi_{R_j}Z^i)\bigr).
\end{aligned}
\]
For each $i\geq1$, since
\[
|y_j^i|-2R_j\longrightarrow\infty,
\]
condition \eqref{hyp:continuous} gives
\[
\|q-\bar q\|_{L^\infty(B_{2R_j}(y_j^i))}
\longrightarrow0.
\]
Consequently,
\[
\begin{aligned}
&\left|
J_{q,\varepsilon_j}
\bigl(\mathsf U_{y_j^i}(\chi_{R_j}Z^i)\bigr)
-J_{\bar q,\varepsilon_j}(\chi_{R_j}Z^i)
\right|\\
&\qquad\leq
C\|q-\bar q\|_{L^\infty(B_{2R_j}(y_j^i))}
\|\chi_{R_j}Z^i\|_{L^p}^p
=o(1),
\end{aligned}
\]
and similarly,
\[
\begin{aligned}
&\left\|
\mathfrak R_{q,\varepsilon_j}
\bigl(\mathsf U_{y_j^i}(\chi_{R_j}Z^i)\bigr)
-\mathsf U_{y_j^i}
\mathfrak R_{\bar q,\varepsilon_j}(\chi_{R_j}Z^i)
\right\|_{\mathcal X^*}\\
&\qquad\leq
C\|q-\bar q\|_{L^\infty(B_{2R_j}(y_j^i))}
\|\chi_{R_j}Z^i\|_{L^p}^{p-1}
=o(1).
\end{aligned}
\]
Combining these estimates with the preceding truncation
limits and translation invariance, we obtain
\[
\begin{aligned}
\|\mathfrak R_{q,\varepsilon_j}(z_j)\|_{\mathcal X^*}
&\leq
\|\mathfrak R_{q,\varepsilon_j}(\chi_{R_j}z^0)\|_{\mathcal X^*}\\
&\quad+
\sum_{i=1}^k
\|\mathfrak R_{\bar q,\varepsilon_j}(\chi_{R_j}Z^i)\|_{\mathcal X^*}
+o(1)
\longrightarrow0,
\end{aligned}
\]
and
\[
\begin{aligned}
J_{q,\varepsilon_j}(z_j)
&=J_{q,\varepsilon_j}(\chi_{R_j}z^0)
+\sum_{i=1}^kJ_{\bar q,\varepsilon_j}(\chi_{R_j}Z^i)
+o(1)\\
&\longrightarrow
J_q(z^0)+\sum_{i=1}^kJ_{\bar q}(Z^i).
\end{aligned}
\]
Thus, $(z_j)$ is admissible and satisfies
\eqref{eq:singular-ps-hypothesis} at the level
\eqref{eq:singular-recovery-level}. The argument also applies
when $\varepsilon_j\equiv0$.
The approximation $\chi_{R_j}Z^i\to Z^i$ proves the stated
strong remainder estimate. The logarithmic density splittings
follow from the choice of radii, disjoint supports, and
dominated convergence. The radial and mixed-mass
splittings follow from their $p$-growth bounds. The source
splitting follows from the principal flux splitting and the
residual identities, exactly as in
\eqref{eq:singular-principal-splitting}.
Escaping translates converge weakly to zero, so the weak
limit is $z^0$. Since $\chi_{R_j}z^0\to z^0$ strongly in $\mathcal X$, we have
\[
S_q(z_j-z^0)
=
S_q\left(\sum_{i=1}^k
\mathsf U_{y_j^i}(\chi_{R_j}Z^i)\right)+o(1).
\]
By the disjointness of the supports,
\[
S_q(z_j-z^0)
=
\sum_{i=1}^k
S_q\bigl(\mathsf U_{y_j^i}(\chi_{R_j}Z^i)\bigr)+o(1).
\]
Since $|y_j^i|-2R_j\to\infty$, assumption~(C) allows us to
replace $q$ by $\bar q$ on each translated support, with an
error tending to zero. By translation invariance and
$\chi_{R_j}Z^i\to Z^i$ strongly in $\mathcal X$, we obtain
\[
S_q(z_j-z^0)\longrightarrow
\sum_{i=1}^k S_{\bar q}(Z^i).
\]
By the definition of the norm,
$\|z\|_q^p=S_q(z)$, and its equivalence to the
$\mathcal X$-norm, if $k=0$, then
\[
z_j\longrightarrow z^0
\qquad\text{strongly in }\mathcal X.
\]
If $k\ge1$, the nontriviality of each $Z^i$ gives
\[
\lim_{j\to\infty}\|z_j-z^0\|_q^p
=\sum_{i=1}^k S_{\bar q}(Z^i)>0.
\]
Since $z_j\rightharpoonup z^0$ in $\mathcal X$, no
subsequence can converge strongly, and we complete the proof.
\end{proof}

\begin{corollary}[Autonomous ground states as regularisation vanishes]
\label{cor:singular-autonomous-ground-states}
The autonomous Nehari levels satisfy
\begin{equation}
d_{\bar q,\varepsilon}\uparrow d_{\bar q,0}
\qquad\text{as }\varepsilon\downarrow0.
\label{eq:singular-ground-level-limit}
\end{equation}
Moreover, $d_{\bar q,0}$ is attained by a nonzero,
componentwise nonnegative finite-energy solution of the
unregularised autonomous system. Furthermore, any sequence
of ground states of $J_{\bar q,\varepsilon_j}$, with
$\varepsilon_j\downarrow0$, is precompact in $\mathcal X$
modulo translations, and every limit is an unregularised
autonomous ground state.
\end{corollary}

\begin{proof}
By Proposition~2.2 and the monotonicity in (2.17),
$d_{\bar q,\varepsilon}$ is nondecreasing as $\varepsilon$
decreases. Comparison with every unregularised finite-energy
Nehari ray gives $d_{\bar q,\varepsilon}\leq d_{\bar q,0}$.
Hence, by Corollary~2.1,
\[
0<\delta
\leq \lim_{\varepsilon\downarrow0}d_{\bar q,\varepsilon}
\leq d_{\bar q,0}.
\]
Set
\[
c:=\lim_{\varepsilon\downarrow0}d_{\bar q,\varepsilon}.
\]

Let $\varepsilon_j\downarrow0$. By
Theorem~\ref{thm:fixed-autonomous-ground-state}, there exist
nonzero, componentwise nonnegative ground states $Z_j$ such
that
\[
J_{\bar q,\varepsilon_j}(Z_j)=\bar d_{\varepsilon_j},
\qquad
J'_{\bar q,\varepsilon_j}(Z_j)=0.
\]
Thus $(Z_j)$ satisfies the hypotheses of
Theorem~\ref{thm:singular-global-decomposition} at level $c$.
Applying that theorem with $q=\bar q$, we obtain a weak limit
$Z^0$ and finitely many nonzero profiles
$Z^1,\ldots,Z^k$. Every nonzero limiting state belongs to
$\mathcal N_{\bar q,0}$, and hence has energy at least
$d_{\bar q,0}$. The exact energy identity gives
\[
c=J_{\bar q,0}(Z^0)
+\sum_{i=1}^k J_{\bar q,0}(Z^i).
\]
Consequently, if $m$ denotes the number of nonzero limiting
states, then
\[
0<c,\qquad
md_{\bar q,0}\leq c\leq d_{\bar q,0}.
\]
Therefore,
\[
m=1,\qquad c=d_{\bar q,0}.
\]
In particular, the unique nonzero limiting state $Z$
satisfies
\[
Z\in\mathcal N_{\bar q,0},
\qquad
J_{\bar q,0}(Z)=d_{\bar q,0}.
\]
Thus $d_{\bar q,0}$ is attained, and
\[
d_{\bar q,\varepsilon}\uparrow d_{\bar q,0}
\qquad\text{as }\varepsilon\downarrow0.
\]

It remains to establish compactness modulo translations.
Since there is exactly one nonzero limiting state, either
$Z^0=Z$ and $k=0$, or $Z^0=0$ and there is exactly one
nonzero escaping profile, say $Z^i$.
In the first case, the remainder estimate in
Theorem~\ref{thm:singular-global-decomposition} gives
\[
Z_j\longrightarrow Z
\qquad\text{strongly in }\mathcal X.
\]
In the second case, it gives
\[
\mathsf U_{-y_j^i}Z_j
\longrightarrow Z^i
\qquad\text{strongly in }\mathcal X.
\]
Hence every sequence of regularised ground states is
precompact in $\mathcal X$ modulo translations, with every
limit an unregularised ground state.

Finally, since the regularised ground states can be chosen
componentwise nonnegative and translations preserve
nonnegativity, their strong limit is also componentwise
nonnegative. This completes the proof.
\end{proof}
\label{subsec:singular-independent-rates}

The following examples explain why residual control and scalar
localisation are needed: strong Sobolev convergence alone does
not ensure convergence of the logarithmic energy, and the
finite-energy domain need not be closed under subtraction.

\begin{example}[The effective domain and the two remainders]
\label{rem:singular-effective-domain}
The entropy estimate for the localised remainder $w_j$
does not extend automatically to the remainder $e_j$
in \eqref{eq:singular-profile-decomposition}, which
converges strongly in $\mathcal X$. Assigning a separate
logarithmic energy to $e_j$ would require an additional
domain argument. Indeed, $\mathcal D_q$ need not be a
vector space. To see this, choose a nonnegative smooth
function $f$ such that
\[
f(x)=|x|^{-1}(\log|x|)^{-3/(2p)}
\qquad (|x|\geq e^2).
\]
Since $p=N$, we have $f\in W^{1,p}$ but
\[
\int_{\mathbb R^N}f^p|\log f|\,dx=\infty.
\]
Thus $(f,0)$ and $(0,f)$ belong to $\mathcal D_q$,
whereas their sum does not. Scalar cutoffs preserve
the domain by \eqref{eq:singular-modulation-negative},
which motivates their use in the exact decomposition
and the reverse construction.
\end{example}

\begin{example}[A diffuse entropy defect without residual control]
Choose $\phi\in C_c^\infty(\mathbb R^N)$ with $0\le\phi\le1$
and
\[
\int_{\mathbb R^N}\phi^p\,dx=1.
\]
Set
\[
t_j=e^{-j},\qquad
\varrho_j=e^j j^{-1/p},\qquad
\varepsilon_j=e^{-4j},
\]
and define
\[
\widetilde z_j(x)
=
\bigl(t_j\phi(x/\varrho_j),t_j\phi(x/\varrho_j)\bigr).
\]
Using $N=p$, we obtain
\[
S_{\bar q}(\widetilde z_j)
=
2e^{-pj}\int_{\mathbb R^N}|\nabla\phi|^p\,dx
+\frac{\bar a+\bar b}{j}
\longrightarrow0,
\qquad
M(\widetilde z_j)=\frac2j.
\]
Thus $\widetilde z_j\to0$ in $\mathcal X$ and
$0\le R_{\varepsilon_j}(\widetilde z_j)
\le M(\widetilde z_j)\to0$.
The logarithms are nonpositive for $j\ge1$, so the positive
part vanishes. A change of variables gives
\begin{equation}
L_{\varepsilon_j}^-(\widetilde z_j)
=
4-\frac2j\int_{\mathbb R^N}
\phi^p\log(\phi^2+e^{-2j})\,dx
\longrightarrow4.
\label{eq:singular-diffuse-entropy}
\end{equation}
Indeed, the integral is uniformly bounded because
$\phi^p|\log\phi^2|$ is integrable, with the convention
$0^p|\log0|=0$. Consequently,
\[
J_{\bar q,\varepsilon_j}(\widetilde z_j)
\longrightarrow\frac4p.
\]

The radial identity gives
\[
\langle J_{\bar q,\varepsilon_j}'(\widetilde z_j),
        \widetilde z_j\rangle
=
S_{\bar q}(\widetilde z_j)
-L_{\varepsilon_j}(\widetilde z_j)
-\frac2pR_{\varepsilon_j}(\widetilde z_j)
\longrightarrow4.
\]
Since $\|\widetilde z_j\|_{\bar q}\to0$, it follows that
\[
\|J_{\bar q,\varepsilon_j}'(\widetilde z_j)\|_{\mathcal X^*}
\ge
\frac{
\bigl|\langle J_{\bar q,\varepsilon_j}'(\widetilde z_j),
              \widetilde z_j\rangle\bigr|
}{
\|\widetilde z_j\|_{\bar q}
}
\longrightarrow\infty.
\]
The same energy defect occurs at $\varepsilon_j=0$, since
\[
L_0^-(\widetilde z_j)
=
4-\frac2j\int_{\mathbb R^N}\phi^p\log\phi^2\,dx
\longrightarrow4,
\qquad
J_{\bar q}(\widetilde z_j)\longrightarrow\frac4p.
\]
\end{example}

\section{Binding, escape, and the geometry of the compactness threshold}
\label{sec:binding-and-escape}

Throughout this section, we assume \eqref{hyp:continuous},
with $p=N>4$. Binding and compactness are linked by
concentration--compactness theory \cite{Lio84a,Lio84b}.
Related advances concern nonlinear Schr\"odinger ground states
\cite{GLN21} and compactness of Lieb--Thirring optimizing
sequences \cite{FGL25}. Here we treat the singular logarithmic
endpoint, where the energy is extended-valued. Profile
decomposition and uniform variational recovery yield the exact
compactness failure levels, convergence of Nehari levels, and
a binding--escape alternative for arbitrary Nehari minimising
sequences. We also establish a Danskin-type directional formula
\cite{Dan67} that incorporates escape through a zero density
functional. This identifies the superdifferential and
distinguishes G\^ateaux from Fr\'echet differentiability at the
threshold, including when the ground-state level is not attained.

\subsection{The exact singular failure levels}
\label{subsec:singular-failure-levels}

For the residual formulation introduced in
Section~\ref{subsec:singular-residual-framework}, define
\[
\mathfrak C_{q,0}
:=
\left\{
J_q(z):z\in \mathcal D_q,
\ \mathfrak R_{q,0}(z)=0
\right\}.
\]
Thus $0\in\mathfrak C_{q,0}$, while every nonzero
member is positive by
\eqref{eq:singular-profile-radial-identity}.

Fix any sequence $\varepsilon_j\in[0,1/2]$ with
$\varepsilon_j\to0$. Compactness at level $c$ means
that every sequence satisfying
\eqref{eq:singular-ps-hypothesis} has a strongly
convergent subsequence in $\mathcal X$. We now proceed to prove Theorem~\ref{thm:binding-failure-spectrum}.

\begin{proposition}[Exact singular compactness threshold]
\label{thm:binding-failure-spectrum}
Along every prescribed parameter sequence, compactness fails
precisely at the levels
\begin{equation}\label{eq:binding-failure-spectrum}
\left\{
c_0+\sum_{i=1}^k c_i:
k\ge1,\quad
c_0\in\mathfrak C_{q,0},\quad
c_i\in\mathfrak C_{\bar q,0}\setminus\{0\}
\right\}.
\end{equation}
In particular, compactness holds below \(d_{\bar q,0}\) and fails
at \(d_{\bar q,0}\). Every noncompact residual Palais--Smale
sequence at this threshold has, after passage to a subsequence,
the form
\begin{equation}\label{eq:binding-first-escape}
z_j=\mathsf U_{y_j}Z+o_{\mathcal X}(1),
\qquad |y_j|\to\infty,
\qquad J_{\bar q}(Z)=d_{\bar q,0},
\end{equation}
where \(Z\) is an autonomous ground state. In the translated
frame, the logarithmic densities converge in \(L^1\) and the
sources in \(\mathcal X^*\), as in
\eqref{eq:singular-entropy-splitting} and
\eqref{eq:singular-source-splitting}.
\end{proposition}
    \begin{proof}
We first prove the necessity. Let $(z_j)$ be a sequence
at level $c$ with vanishing residuals but no strongly
convergent subsequence. By
Theorem~\ref{thm:singular-global-decomposition}, after
passing to a subsequence, there exist a possibly zero
finite-energy solution $z^0$ for $q$ and nonzero
finite-energy solutions $Z^1,\ldots,Z^k$ for $\bar q$,
with $k\geq1$. The exact energy splitting gives
\[
c=J_q(z^0)+\sum_{i=1}^k J_{\bar q}(Z^i).
\]
By the definition of $\mathfrak C_{q,0}$ and
$\mathfrak C_{\bar q,0}$, this yields
\[
c\in
\left\{
c_0+\sum_{i=1}^k c_i:
k\geq1,\quad
c_0\in\mathfrak C_{q,0},\quad
c_i\in\mathfrak C_{\bar q,0}\setminus\{0\}
\right\}.
\]

Conversely, let $c$ belong to the above set and choose
finite-energy solutions realizing the corresponding
levels. Theorem~\ref{thm:singular-profile-realisation}
constructs a noncompact sequence at level $c$ for the
prescribed parameter sequence, including the case
$\varepsilon_j\equiv0$. This proves the exact
characterization of the failure levels.

By the radial identity, every nonzero escaping profile
belongs to $\mathcal N_{\bar q,0}$ and therefore has
energy at least $d_{\bar q,0}$. Moreover,
$J_q(z^0)\geq0$. Hence every failure level satisfies
\[
c=J_q(z^0)+\sum_{i=1}^kJ_{\bar q}(Z^i)
\geq k d_{\bar q,0}\geq d_{\bar q,0}.
\]
Thus compactness holds below $d_{\bar q,0}$. On the other
hand, Corollary~\ref{cor:singular-autonomous-ground-states}
provides an autonomous ground state $Z$ with
$J_{\bar q}(Z)=d_{\bar q,0}$. Taking $z^0=0$ and $k=1$
in the realisation theorem shows that compactness fails
at $d_{\bar q,0}$.

Finally, consider a sequence at $c=d_{\bar q,0}$ without
a strongly convergent subsequence. Its energy identity
gives
\[
d_{\bar q,0}
=J_q(z^0)+\sum_{i=1}^kJ_{\bar q}(Z^i)
\geq J_q(z^0)+kd_{\bar q,0}.
\]
Since $k\geq1$ and $J_q(z^0)\geq0$, necessarily
\[
k=1,\qquad J_q(z^0)=0,\qquad
J_{\bar q}(Z^1)=d_{\bar q,0}.
\]
The radial identity further implies $z^0=0$.
Consequently, writing $Z=Z^1$ and $y_j=y_j^1$,
the profile decomposition yields
\[
z_j=\mathsf U_{y_j}Z+o_{\mathcal X}(1),
\qquad |y_j|\to\infty,
\qquad J_{\bar q}(Z)=d_{\bar q,0}.
\]
Thus $Z$ is an autonomous ground state. The density
and source convergence follow from
\eqref{eq:singular-energy-splitting} and the
corresponding splittings in
Theorem~\ref{thm:singular-global-decomposition}.
\end{proof}

\subsection{Normalised energies and minimising sequences}
\label{subsec:continuous-normalisation}

\begin{definition}[The normalised variational problem]
\label{def:binding-normalised-problem}
Since the effective domain is independent of the choice
of bounded potential with a positive lower bound, denote
it by $\mathcal D$. We define
\[
\begin{gathered}
\Sigma:=\{w\in\mathcal D:w\geq0,\ M(w)=1\},\\
\Phi_q(w):=S_q(w)-L(w),\qquad
F(q):=\inf_{w\in\Sigma}\Phi_q(w),\\
\mathcal W(q):=\operatorname*{argmin}_{w\in\Sigma}\Phi_q(w).
\end{gathered}
\]
The set $\mathcal W(q)$ is allowed to be empty.
\end{definition}

\begin{proposition}[Representation, coercivity, and comparison]
\label{prop:binding-normalised-energy}
The following statements hold:
\begin{enumerate}[label=\textup{(\roman*)},itemsep=4pt]
\item One has
\begin{equation}
d_{q,0}
=\frac{2}{p^2}\exp\left(\frac p2 F(q)-1\right),
\qquad
0<d_{q,\varepsilon}\leq d_{\bar q,\varepsilon}
\quad (0\leq\varepsilon\leq1/2).
\label{eq:binding-exponential-comparison}
\end{equation}
Moreover, on any family of potentials with fixed positive
lower and upper bounds, every sublevel
$\Phi_q(w)\leq E$ in $\Sigma$ satisfies uniform bounds
for $\|w\|_{\bar q}$ and $L_0^\pm(w)$.

\item The functional $F$ is finite, concave, and locally
Lipschitz with respect to bounded potential perturbations
preserving the asymptotic potential $\bar q$.
\end{enumerate}
\end{proposition}

\begin{proof}
    \textup{(i)} By the ray formula
\eqref{eq:exact-ray-energy} and the restriction to
nonnegative rays in Corollary~\ref{cor:uniform-levels},
we obtain
\[
d_{q,0}
=\frac{2}{p^2}\exp\left(\frac p2F(q)-1\right).
\]
Since $0<d_{q,0}<\infty$, $F(q)$ is finite.

To prove the uniform sublevel bounds, choose
$0<\theta<1$ and $0<\sigma<p\theta$. On $\Sigma$,
the calibration estimate \eqref{eq:calibration}, with
$R_0=M=1$, gives
\[
\Phi_q(w)\geq S_q(w)-CS_q(w)^\beta,
\qquad
\beta=1-\theta+\frac{\sigma}{p}<1.
\]
Since the right-hand side tends to infinity as
$S_q(w)\to\infty$, every sublevel $\Phi_q(w)\leq E$
has a uniform bound for $S_q(w)$, and hence for
$\|w\|_{\bar q}$, under the stated potential bounds.
By \eqref{eq:positive-growth}, $L_0^+(w)$ is also
uniformly bounded. Finally,
\[
L_0^-(w)=\Phi_q(w)-S_q(w)+L_0^+(w)
\]
gives the corresponding bound for $L_0^-(w)$.

The positivity of $d_{q,\varepsilon}$ follows from
Corollary~\ref{cor:uniform-levels}. It remains to prove
$d_{q,\varepsilon}\leq d_{\bar q,\varepsilon}$.
Fix $\varepsilon\in[0,1/2]$ and let $Z_\varepsilon$
be an autonomous ground state, whose existence follows
from Theorem~\ref{thm:fixed-autonomous-ground-state}
if $\varepsilon>0$ and
Corollary~\ref{cor:singular-autonomous-ground-states}
if $\varepsilon=0$.

By \eqref{eq:compact-potential-energy} and translation
invariance of the remaining terms,
\[
S_q(\mathsf U_yZ_\varepsilon)
\longrightarrow S_{\bar q}(Z_\varepsilon)
\qquad (|y|\to\infty).
\]
Consequently, the radial Nehari equations for
$\mathsf U_yZ_\varepsilon$ converge, at every fixed
positive scale, to the autonomous radial equation for
$Z_\varepsilon$. Since the latter has a unique root at
$t=1$, strict radial monotonicity in
Proposition~\ref{prop:radial-projection} yields
\[
t_{q,\varepsilon}(\mathsf U_yZ_\varepsilon)
\longrightarrow1.
\]
Indeed, Since $Z_\varepsilon\in\mathcal N_{\bar q,\varepsilon}$,
strict radial monotonicity gives, for every $0<\delta<1$,
\[
G_{\bar q,\varepsilon}((1-\delta)Z_\varepsilon)>0,
\qquad
G_{\bar q,\varepsilon}((1+\delta)Z_\varepsilon)<0.
\]
For each fixed $t>0$, translation invariance of the nonlinear
terms yields
\[
\begin{aligned}
&G_{q,\varepsilon}(t\mathsf U_yZ_\varepsilon)
-G_{\bar q,\varepsilon}(tZ_\varepsilon)\\
&\qquad
=t^p\bigl(S_q(\mathsf U_yZ_\varepsilon)
-S_{\bar q}(Z_\varepsilon)\bigr)
\longrightarrow0
\qquad(|y|\to\infty).
\end{aligned}
\]
Thus, for sufficiently large $|y|$,
\[
G_{q,\varepsilon}((1-\delta)\mathsf U_yZ_\varepsilon)>0,
\qquad
G_{q,\varepsilon}((1+\delta)\mathsf U_yZ_\varepsilon)<0.
\]
By Proposition~\ref{prop:radial-projection},
\[
1-\delta<t_{q,\varepsilon}(\mathsf U_yZ_\varepsilon)<1+\delta.
\]
Since $\delta\in(0,1)$ is arbitrary,
\[
t_{q,\varepsilon}(\mathsf U_yZ_\varepsilon)\longrightarrow1
\qquad(|y|\to\infty).
\]
Evaluating the energy at this Nehari projection, we obtain
\[
\begin{aligned}
d_{q,\varepsilon}
&\leq
J_{q,\varepsilon}\bigl(
t_{q,\varepsilon}(\mathsf U_yZ_\varepsilon)
\mathsf U_yZ_\varepsilon\bigr)\\
&\longrightarrow
J_{\bar q,\varepsilon}(Z_\varepsilon)
=d_{\bar q,\varepsilon}.
\end{aligned}
\]
This proves the comparison.

\medskip
\noindent
\textup{(ii)} For each fixed $w\in\Sigma$,
$\Phi_q(w)$ is affine in $q$. Hence $F(q)$, being
the infimum of affine functions, is concave.

It remains to prove local Lipschitz continuity. Fix a
potential $q$ and consider a sufficiently small
$L^\infty$-neighbourhood in which the positive potential
bounds remain uniform and the asymptotic potential
$\bar q$ remains unchanged. A fixed test function in $\Sigma$
provides a uniform upper bound for $F$ in this
neighbourhood. By part~(i), all corresponding almost
minimisers are therefore uniformly bounded in
$\|\cdot\|_{\bar q}$.

Let $q_1,q_2$ belong to this neighbourhood, and choose
$w\in\Sigma$ such that
\[
\Phi_{q_1}(w)\leq F(q_1)+\eta,
\qquad \eta>0.
\]
The uniform sublevel bound gives
\[
\begin{aligned}
|\Phi_{q_2}(w)-\Phi_{q_1}(w)|
&\leq
\|q_2-q_1\|_{L^\infty}
\bigl(\|w_1\|_{L^p}^p+\|w_2\|_{L^p}^p\bigr)\\
&\leq C\|q_2-q_1\|_{L^\infty}.
\end{aligned}
\]
Consequently,
\[
F(q_2)\leq F(q_1)+\eta
+C\|q_2-q_1\|_{L^\infty}.
\]
Letting $\eta\downarrow0$ and interchanging $q_1$
and $q_2$, we conclude that
\[
|F(q_2)-F(q_1)|
\leq C\|q_2-q_1\|_{L^\infty}.
\]
Thus $F$ is locally Lipschitz, which completes the proof.
\end{proof}

We now proceed to the proof of Theorem~\ref{thm:binding-minimisers}, which applies to arbitrary Nehari minimising sequences. The required passage to small equation residuals is supplied by Lemma~\ref{lem:appendix-uniform-lifting}; its constants do not deteriorate as $\varepsilon\downarrow0$.
\paragraph{Proof of Theorem~\ref{thm:binding-minimisers}}
    \emph{Step 1.} First of all ,we show the following fact: 
    \begin{equation}
 d_{q,\varepsilon}\uparrow d_{q,0}\qquad(\varepsilon\downarrow0).
 \label{eq:binding-level-convergence}
\end{equation}
Indeed, assume that \(\varepsilon_1<\varepsilon_2\), and take \(z\in \mathcal{N}_{q,\varepsilon_1}\), there exists \(t>0\), such that \(tz\in \mathcal{N}_{q,\varepsilon_2}\). Then, 
\[d_{q,\varepsilon_2}\le J_{q,\varepsilon_2}(tz)=\max\limits_{s>0} J_{q,\varepsilon_2}(sz).\]
By monotonicity of the energy, for any \(s>0\), there holds \(J_{q,\varepsilon_2}(sz)\le J_{q,\varepsilon_1}(sz)\). And is thus, \[\max\limits_{s>0} J_{q,\varepsilon_2}(sz)\le \max\limits_{s>0} J_{q,\varepsilon_1}(sz)=J_{q,\varepsilon_1}(z).\]
Take the infimum on the right side gives \[d_{q,\varepsilon_2}\le d_{q,\varepsilon_1}.\]
Its limit $c$ is positive by Corollary~\ref{cor:uniform-levels}.
For any positive $\varepsilon_j\to0$, choose near minimisers on
$\mathcal N_{q,\varepsilon_j}$ and apply
Lemma~\ref{lem:appendix-uniform-lifting}. This gives
\[
 J_{q,\varepsilon_j}(\zeta_j)\to c,\qquad
 \|J_{q,\varepsilon_j}'(\zeta_j)\|_{\mathcal X^*}\to0.
\]
Applying Theorem~\ref{thm:singular-global-decomposition}, \[c=J_{q,0}(z^0)+\sum_{i=1}^{k} J_{\bar{q},0}(Z^i).\]
If \(z^0\neq 0\), since \(z^0\) is a nontrivial critical point of the energy with potential \(q\), there holds \(J_{q,0}(z^0)\ge d_{q,0}\). If \(Z^i\neq 0\), then \(J_{\bar{q},0}(Z^i)\ge \bar{d_0}\).
Proposition~\ref{prop:binding-normalised-energy} gives \[d_{q,0}\le d_{\bar{q},0},\] 
which immediately implies that \(c\ge d_{q,0}\). On the other hand, by the monotonicity comparison, it follows that \(c\le d_{q,0}\), which implies \(c=d_{q,0}\). Thus, \eqref{eq:binding-level-convergence} holds.

\emph{Step 2.}
By the minimising hypothesis and \eqref{eq:binding-level-convergence},
\[
 J_{q,\varepsilon_j}(\widehat z_j)\longrightarrow d_{q,0}.
\]
We approximate \(\widehat z_j\) by points for which the equation
residual tends to zero.

If \(\varepsilon_j>0\), apply
Lemma~\ref{lem:appendix-uniform-lifting} directly to
\(\widehat z_j\in\mathcal N_{q,\varepsilon_j}\), since
\(J_{q,\varepsilon_j}(\widehat z_j)-d_{q,\varepsilon_j}\to0\).

If \(\varepsilon_j=0\), fix \(j\). Since
\(\widehat z_j\in\mathcal N_{q,0}\), its radial projection at
parameter zero is \(\widehat z_j\) itself. By
Lemma~\ref{lem:appendix-ray-recovery} and
\eqref{eq:binding-level-convergence}, we may choose a positive
parameter \(\varepsilon<1/j\) such that
\[
 \begin{aligned}
 \bigl\|t_{q,\varepsilon}(\widehat z_j)\widehat z_j
                 -\widehat z_j\bigr\|_{\bar q}
 &<\frac1j,\\
 \bigl|J_{q,\varepsilon}
       \bigl(t_{q,\varepsilon}(\widehat z_j)\widehat z_j\bigr)
                 -J_q(\widehat z_j)\bigr|
 &<\frac1j,
 \qquad
 |d_{q,\varepsilon}-d_{q,0}|<\frac1j.
 \end{aligned}
\]
The projected point belongs to \(\mathcal N_{q,\varepsilon}\),
and its energy exceeds the corresponding Nehari level by at most
\[
 J_{q,\varepsilon}
 \bigl(t_{q,\varepsilon}(\widehat z_j)\widehat z_j\bigr)
 -d_{q,\varepsilon}
 \le J_q(\widehat z_j)-d_{q,0}+\frac2j
 \longrightarrow0.
\]
We may therefore apply Lemma~\ref{lem:appendix-uniform-lifting}
to the projected point.

For each \(j\), let \(\varepsilon\) denote the positive parameter
used above: \(\varepsilon=\varepsilon_j\) if \(\varepsilon_j>0\),
and \(0<\varepsilon<1/j\) otherwise. In either case,
\(\varepsilon\to0\) as \(j\to\infty\), and the projected point
\(t_{q,\varepsilon}(\widehat z_j)\widehat z_j\) has energy
\(d_{q,\varepsilon}+o(1)\). Applying
Lemma~\ref{lem:appendix-uniform-lifting} to this point gives
\(\zeta_j\in\mathcal N_{q,\varepsilon}\) within distance \(o(1)\)
of it, with \(\|J_{q,\varepsilon}'(\zeta_j)\|_{\mathcal X^*}=o(1)\).

The projection leaves \(\widehat z_j\) unchanged if
\(\varepsilon_j>0\) and moves it by less than \(1/j\) otherwise.
Thus
\[
 \|\zeta_j-\widehat z_j\|_{\bar q}
 \le \|\zeta_j-t_{q,\varepsilon}(\widehat z_j)\widehat z_j\|_{\bar q}
      +\frac1j
 \longrightarrow0.
\]
Moreover, the energy bounds in the lifting lemma give
\[
 d_{q,\varepsilon}
 \le J_{q,\varepsilon}(\zeta_j)
 \le J_{q,\varepsilon}
       \bigl(t_{q,\varepsilon}(\widehat z_j)\widehat z_j\bigr)
 =d_{q,\varepsilon}+o(1).
\]
Since \(d_{q,\varepsilon}\to d_{q,0}\) by
\eqref{eq:binding-level-convergence}, we obtain
\begin{equation}
 \begin{gathered}
 \|\zeta_j-\widehat z_j\|_{\bar q}\longrightarrow0,
 \qquad
 J_{q,\varepsilon}(\zeta_j)\longrightarrow d_{q,0},\\
 \|J_{q,\varepsilon}'(\zeta_j)\|_{\mathcal X^*}
 \longrightarrow0.
 \end{gathered}
 \label{eq:binding-minimiser-lifting}
\end{equation}
Applying Theorem~\ref{thm:singular-global-decomposition} to
\(\zeta_j\), we obtain, after passage to a subsequence,
\[
 d_{q,0}=J_q(z^0)+\sum_{i=1}^k J_{\bar q}(Z^i).
\]
Every nonzero state contributes at least \(d_{q,0}>0\), since
each escaping profile has energy at least \(d_{\bar q,0}\ge d_{q,0}\).
Hence exactly one state is nonzero.

If \(z^0\ne0\), then \(k=0\), \(J_q(z^0)=d_{q,0}\), and
\(\zeta_j\to z^0\) in \(\mathcal X\).
If \(z^0=0\), then \(k=1\) and
\[
 d_{q,0}=J_{\bar q}(Z^1)\ge d_{\bar q,0}\ge d_{q,0}.
\]
Thus \(Z^1\) is an autonomous ground state and
\(\|\zeta_j-\mathsf U_{y_j^1}Z^1\|_{\bar q}\to0\),
with \(|y_j^1|\to\infty\).

Finally, \eqref{eq:binding-minimiser-lifting} gives
\(\|\widehat z_j-\zeta_j\|_{\bar q}\to0\).
In the first case, the triangle inequality yields
\[
 \|\widehat z_j-z^0\|_{\bar q}
 \le \|\widehat z_j-\zeta_j\|_{\bar q}
    +\|\zeta_j-z^0\|_{\bar q}\longrightarrow0.
\]
In the second case, the same estimate with
\(\mathsf U_{y_j^1}Z^1\) in place of \(z^0\) gives
\(\|\widehat z_j-\mathsf U_{y_j^1}Z^1\|_{\bar q}\to0\).

\emph{Step 3.}
We first show that every minimizer is a weak solution.
Suppose that \(z_*\in\mathcal N_{q,0}\) satisfies
\(J_{q,0}(z_*)=d_{q,0}\). Taking the constant minimizing
sequence \(\widehat z_j=z_*\), the binding--escape alternative
applies. The escape alternative is impossible, since
\(\mathsf U_{y_j}Z\rightharpoonup0\) in \(X\) as
\(|y_j|\to\infty\), whereas
\(z_*-\mathsf U_{y_j}Z\to0\) strongly in \(X\) would imply
\(z_*=0\), a contradiction. Hence
\(\mathcal R_{q,0}(z_*)=0\), so \(z_*\) is a weak solution.

Next, assume that \(d_{q,0}<d_{\bar q,0}\). By
Corollary~2.1, there exists a nonnegative minimizing sequence
\(\{\widehat z_j\}\subset\mathcal N_{q,0}\).
The escape alternative would imply
\(d_{q,0}=d_{\bar q,0}\), which is impossible. Thus, up to a
subsequence,
\[
\widehat z_j\longrightarrow z\quad\text{strongly in }X,
\qquad
\mathcal R_{q,0}(z)=0,
\qquad
J_{q,0}(z)=d_{q,0}.
\]
Since \(\widehat z_j\ge0\), we have \(z\ge0\).
Consequently, \(z\) is a nonnegative ground state.
Moreover, the same argument applies to every minimizing
sequence, showing that each has a strongly convergent
subsequence in \(X\).

Finally, suppose that \(d_{q,0}=d_{\bar q,0}\), and let \(Z\)
be a ground state of the autonomous problem. Choose
\(|y_j|\to\infty\) and set
\[
\widehat z_j=t_j\mathsf U_{y_j}Z\in\mathcal N_{q,0},
\qquad
t_j=t_{q,0}(\mathsf U_{y_j}Z)>0.
\]
Subtracting the Nehari identities for \(\widehat z_j\)
and \(Z\), we obtain
\[
2M(Z)\log t_j
=S_q(\mathsf U_{y_j}Z)-S_{\bar q}(Z).
\]
Since \(q(x)\to\bar q\) and the remaining terms are
translation invariant,
\[
S_q(\mathsf U_{y_j}Z)\longrightarrow S_{\bar q}(Z),
\qquad t_j\longrightarrow1.
\]
Consequently, the Nehari energy identity yields
\[
J_{q,0}(\widehat z_j)
=\frac{2}{p^2}t_j^pM(Z)
\longrightarrow
\frac{2}{p^2}M(Z)
=d_{\bar q,0}=d_{q,0}.
\]
Thus \(\{\widehat z_j\}\) is a minimizing sequence. However,
\[
\widehat z_j\rightharpoonup0\quad\text{in }X,
\qquad
\|\widehat z_j\|_{\bar q}
=t_j\|Z\|_{\bar q}
\longrightarrow\|Z\|_{\bar q}>0.
\]
Hence no subsequence converges strongly in \(X\).

Combining the above conclusions with
\(d_{q,0}\le d_{\bar q,0}\), we obtain
\(d_{q,0}<d_{\bar q,0}\) is equivalent to every minimizing sequence in $\mathcal N_{q,0}$ has a strongly convergent subsequence in $\mathcal{X}$.
The proof is completed.

\begin{corollary}[Stability inside the binding region]
\label{cor:binding-stability}
Assume \(d_{q,0}<d_{\bar q,0}\). Then the following hold.
\begin{enumerate}
\item For all sufficiently small \(\varepsilon>0\),
\(d_{q,\varepsilon}<d_{\bar q,\varepsilon}\), and
\(J_{q,\varepsilon}\) admits a componentwise nonnegative ground
state at level \(d_{q,\varepsilon}\). If \(\varepsilon_j\to0\)
and \(z_j\) are such ground states, then, after passage to a
subsequence, \(z_j\to z\) strongly in \(\mathcal X\) for some
ground state \(z\) of \(J_q\).

\item Suppose \(q_j\to q\) in \(L^\infty(\mathbb R^N)^2\)
and each \(q_j\) has the same limiting pair \(\bar q\) at infinity.
Then \(\mathcal W(q_j)\ne\varnothing\) for all sufficiently large
\(j\), and every sequence \(w_j\in\mathcal W(q_j)\) has a
subsequence converging strongly in \(\mathcal X\) to some
\(w\in\mathcal W(q)\).
\end{enumerate}
\end{corollary}

\begin{proof}
\emph{(i)}
By \eqref{eq:binding-level-convergence}, applied to \(q\) and
\(\bar q\),
\[
 d_{\bar q,\varepsilon}-d_{q,\varepsilon}
 \longrightarrow d_{\bar q,0}-d_{q,0}>0.
\]
Fix \(\varepsilon>0\) sufficiently small that
\(d_{q,\varepsilon}<d_{\bar q,\varepsilon}\).
By Corollary~\ref{cor:uniform-levels}, we may choose a
componentwise nonnegative minimising sequence
\(z_n\in\mathcal N_{q,\varepsilon}\).
Lemma~\ref{lem:appendix-uniform-lifting} gives a Palais--Smale
sequence \((\zeta_n)\) at level \(d_{q,\varepsilon}\) with
\(\|\zeta_n-z_n\|_{\bar q}\to0\).
Since this level lies below \(d_{\bar q,\varepsilon}\),
Corollary~\ref{cor:fixed-ps-failure-levels} gives, after passage
to a subsequence, \(\zeta_n\to z_\varepsilon\) in \(\mathcal X\).
Then \(z_n\to z_\varepsilon\) as well, so \(z_\varepsilon\ge0\).
Since \(J_{q,\varepsilon}\) is \(C^1\),
\[
 J_{q,\varepsilon}'(z_\varepsilon)=0,
 \qquad J_{q,\varepsilon}(z_\varepsilon)=d_{q,\varepsilon}>0.
\]
Thus \(z_\varepsilon\) is a nonnegative ground state.

Now let \(\varepsilon_j\to0\) and choose any such ground states
\(z_j\). Since \(J_{q,\varepsilon_j}(z_j)=d_{q,\varepsilon_j}\),
Theorem~\ref{thm:binding-minimisers} applies. Strict binding
excludes escape, so a subsequence converges strongly in
\(\mathcal X\) to a ground state \(z\) of \(J_q\).
Moreover, \(S_q(z_j)\to S_q(z)\) and
\(J_{q,\varepsilon_j}(z_j)\to d_{q,0}=J_q(z)\).
The asserted density convergence follows from
Lemma~\ref{lem:appendix-entropy-recovery}.

\emph{(ii)}
By local Lipschitz continuity of \(F\) and
\eqref{eq:binding-exponential-comparison},
\[
 d_{q_j,0}
 =\frac{2}{p^2}\exp\!\left(\frac p2F(q_j)-1\right)
 \longrightarrow d_{q,0}<d_{\bar q,0}.
\]
Since the limiting pair \(\bar q\) is common to all \(q_j\),
strict binding holds for \(q_j\) for all sufficiently large \(j\).
Theorem~\ref{thm:binding-minimisers} then gives a nonnegative
ground state for \(q_j\); normalising it to have mixed mass one
gives an element of \(\mathcal W(q_j)\).

Take any \(w_j=(u_j,v_j)\in\mathcal W(q_j)\). Since
\(\Phi_{q_j}(w_j)=F(q_j)\to F(q)\), uniform coercivity gives
boundedness of \((w_j)\) in \(\mathcal X\). Only the potential
term changes with \(q\), so
\[
 |\Phi_q(w_j)-F(q_j)|
 \le \|q-q_j\|_\infty
      \bigl(\|u_j\|_{L^p}^p+\|v_j\|_{L^p}^p\bigr)
 \longrightarrow0.
\]
Thus \(\Phi_q(w_j)\to F(q)\). Set \(t_j=t_{q,0}(w_j)\).
As \(M(w_j)=1\), the radial formula \eqref{eq:exact-ray-scale}
gives
\[
 t_j=\exp\!\left(\frac{\Phi_q(w_j)}2-\frac1p\right)
 \longrightarrow
 t_*:=\exp\!\left(\frac{F(q)}2-\frac1p\right)>0.
\]
Hence \(t_jw_j\in\mathcal N_{q,0}\) and
\[
 J_q(t_jw_j)=\frac{2}{p^2}t_j^p
 \longrightarrow \frac{2}{p^2}t_*^p=d_{q,0}.
\]
Applying Theorem~\ref{thm:binding-minimisers} at the fixed
potential \(q\), strict binding gives \(t_jw_j\to z\) in
\(\mathcal X\) along a subsequence, where \(z\) is a nonnegative
ground state. Since \(M(z)=\lim_j M(t_jw_j)=t_*^p\), we conclude
that
\[
 w_j\longrightarrow \frac{z}{t_*}
 =\frac{z}{M(z)^{1/p}}\in\mathcal W(q)
 \qquad\text{in }\mathcal X.
\]
\end{proof}

\subsection{Potential derivatives with an escape state}
\label{subsec:binding-density-geometry}
Fix the limiting potential \(\bar q=(\bar a,\bar b)\) and set
\begin{equation}
 \begin{aligned}
 \mathcal P_0
 &=\left\{
 h\in L^\infty(\mathbb R^N)^2:
 \lim_{R\to\infty}\|h\|_{L^\infty(\{|x|>R\})}=0
 \right\},\\
 \Omega_\infty
 &=\left\{
 q=(a,b)\in\bar q+\mathcal P_0:
 \operatorname*{ess\,inf}a>0,\ 
 \operatorname*{ess\,inf}b>0
 \right\}.
 \end{aligned}
 \label{eq:binding-potential-space}
\end{equation}
We equip \(\mathcal P_0\) with the maximum of the two component
\(L^\infty\) norms. Then \(\Omega_\infty\) is open in
\(\bar q+\mathcal P_0\).

For \(w=(w_u,w_v)\in\Sigma\), the pair
\((|w_u|^p,|w_v|^p)\in L^1(\mathbb R^N)^2\) defines a
functional \(\rho(w)\in\mathcal P_0^*\) by
\begin{equation}
 \rho(w)[h]
 =\int_{\mathbb R^N}
   \bigl(h_a|w_u|^p+h_b|w_v|^p\bigr)\,dx,
 \qquad h=(h_a,h_b)\in\mathcal P_0.
 \label{eq:binding-density-functional}
\end{equation}
Define
\begin{equation}
 \mathcal A(q)=
 \begin{cases}
 \{\rho(w):w\in\mathcal W(q)\},
     & F(q)<F(\bar q),\\
 \{\rho(w):w\in\mathcal W(q)\}\cup\{0\},
     & F(q)=F(\bar q).
 \end{cases}
 \label{eq:binding-active-densities}
\end{equation}
At the threshold, the zero functional represents escape:
for every fixed \(h\in\mathcal P_0\),
\(\rho(\mathsf U_y w)[h]\to0\) as \(|y|\to\infty\), since
\(h\) vanishes at infinity and the density pair is integrable.
No \(w\in\Sigma\) itself has zero density, since \(M(w)=1\).

Finally, the superdifferential of the concave function \(F\) at
\(q\in\Omega_\infty\) is
\[
 \partial^+F(q)
 =\left\{\lambda\in\mathcal P_0^*:
 F(q')\le F(q)+\lambda[q'-q]
 \ \text{for every }q'\in\Omega_\infty\right\}.
\]
We now proceed to prove Theorem~\ref{thm:binding-density-response}.
\paragraph{Proof of Theorem~\ref{thm:binding-density-response}}
\emph{Step 1.}
Let \(w_j\in\Sigma\) satisfy \(\Phi_q(w_j)\to F(q)\).
Since \(M(w_j)=1\), the radial projection formula gives
\(s_jw_j\in\mathcal N_{q,0}\), where
\[
 s_j=\exp\!\left(\frac{\Phi_q(w_j)}{2}-\frac1p\right)
 \longrightarrow
 \exp\!\left(\frac{F(q)}{2}-\frac1p\right)>0.
\]
By the Nehari identity and \eqref{eq:binding-exponential-comparison},
\[
 J_q(s_jw_j)=\frac{2}{p^2}M(s_jw_j)
            =\frac{2}{p^2}s_j^p
 \longrightarrow d_{q,0}.
\]
Thus Theorem~\ref{thm:binding-minimisers} applies to \(s_jw_j\).
Dividing by \(s_j\) normalises the limiting ground state to have
mixed mass one. Hence, after passing to a subsequence, either
\(w_j\to w\) strongly in \(\mathcal X\) for some
\(w\in\mathcal W(q)\), or \(d_{q,0}=d_{\bar q,0}\) and
\begin{equation}
 w_j-\mathsf U_{y_j}W\longrightarrow0
 \quad\text{in }\mathcal X,
 \qquad |y_j|\longrightarrow\infty,
 \qquad W\in\mathcal W(\bar q).
 \label{eq:binding-normalised-alternative}
\end{equation}
By \eqref{eq:binding-exponential-comparison} and the strict
monotonicity of the exponential function, this is equivalent
to \(F(q)=F(\bar q)\).

In the first case, \(w_j\to w\) strongly in \(L^p(\mathbb R^N)^2\).
Hence the \(p\)th powers of the absolute values of the components
converge in \(L^1(\mathbb R^N)\), and
\[
 \rho(w_j)[h]\longrightarrow\rho(w)[h]
 \qquad\text{for every }h\in\mathcal P_0.
\]

In the escape case, \(w_j-\mathsf U_{y_j}W\to0\) in
\(L^p(\mathbb R^N)^2\), so the corresponding densities differ
by \(o(1)\) in \(L^1(\mathbb R^N)^2\). Thus, for each fixed
\(h\in\mathcal P_0\), a change of variables gives
\[
 \rho(w_j)[h]=\rho(W)[h(\cdot+y_j)]+o(1).
\]
Since \(|y_j|\to\infty\) and \(h\) vanishes at infinity,
\(\|h(\cdot+y_j)\|_{L^\infty(B_R)}\to0\) for each fixed \(R>0\).
Outside \(B_R\), the absolute value of the integral is bounded by
\[
 \|h\|_\infty\int_{\mathbb R^N\setminus B_R}
       \bigl(|W_u|^p+|W_v|^p\bigr)\,dx,
\]
which tends to zero as \(R\to\infty\). It follows that
\(\rho(w_j)[h]\to0\).

Thus every minimising sequence for \(F(q)\) has a subsequence
whose density functionals converge weakly-* in \(\mathcal P_0^*\)
to an element of \(\mathcal A(q)\): the limit is \(\rho(w)\)
in the first case and \(0\) in the escape case.

\emph{Step 2.} For every \(w\in \mathcal{W}(q)\), there is
\[F(q+th)\le F(q)+t\rho(w)[h].\]
Indeed, by the definition of \(F\), \[F(q+th)\le \Phi_{q+th}(w)=\Phi_q(w)+t\rho(w)[h]=F(q)+t\rho(w)[h],\]
where we use the fact that \(w\in \mathcal{W}(q)\) in the last equality.
At the threshold, \(F(q+th)\le \Phi_{q+th}(\mathsf U_y W)\), \(\lim\limits_{|y|\to\infty} \Phi_{q+th}(\mathsf U_y W)=\Phi_{\bar{q}}(W)=F(q)\), that is \(F(q+th)\le F(q)\). These inequalities give the upper bound in \eqref{eq:binding-directional-response}.

For the reverse bound, let $t_j\downarrow0$ and choose
$w_j\in\Sigma$ so that
\[
 \Phi_{\mathrm q+t_jh}(w_j)
          \le F(q+t_jh)+t_j^2,
\]
it follows that 
\[
   F(q+t_jh)\ge \Phi_q(w_j)+t_j \rho(w_j)[h]-t_j^2\ge F(q)+t_j \rho(w_j)[h]-t_j^2,
\]
and thus 
\begin{equation}
  \frac{F(q+t_jh)-F(q)}{t_j}\ge \rho(w_j)[h]-t_j
\end{equation}
Choose \(t_j\downarrow0\) so that the difference quotients converge
to their lower limit. Since \(\Phi_q(w_j)\to F(q)\), Step~1 gives,
after passage to a further subsequence,
\(\rho(w_j)[h]\to\lambda_*[h]\) for some
\(\lambda_*\in\mathcal A(q)\). Passing to the limit in the lower
bound above and combining it with the upper bound yields
\[
\begin{aligned}
 \inf_{\lambda\in\mathcal A(q)}\lambda[h]
 &\le \lambda_*[h]
 \le \liminf_{t\downarrow0}\frac{F(q+th)-F(q)}{t}\\
 &\le \limsup_{t\downarrow0}\frac{F(q+th)-F(q)}{t}
 \le \inf_{\lambda\in\mathcal A(q)}\lambda[h].
\end{aligned}
\]
Thus the right directional derivative exists, and the infimum is
attained at \(\lambda_*\). Hence
\[
 F'(q;h)=\min_{\lambda\in\mathcal A(q)}\lambda[h],
\]
and \eqref{eq:binding-directional-response} is established.

To prove the Hadamard assertion, let \(t_j\downarrow0\) and
\(h_j\to h\) in \(\mathcal P_0\). By local Lipschitz continuity
of \(F\) at \(q\), for all sufficiently large \(j\),
\[
 \left|
 \frac{F(q+t_jh_j)-F(q)}{t_j}
 -
 \frac{F(q+t_jh)-F(q)}{t_j}
 \right|
 \le C\|h_j-h\|_\infty\longrightarrow0.
\]
Since the quotient in the fixed direction \(h\) converges to
\(F'(q;h)\), it follows that
\[
 \frac{F(q+t_jh_j)-F(q)}{t_j}
 \longrightarrow F'(q;h)
 =\min_{\lambda\in\mathcal A(q)}\lambda[h].
\]

For every \(w\in\mathcal W(q)\) and \(q'\in\Omega_\infty\),
\[
 F(q')\le \Phi_{q'}(w)
 =\Phi_q(w)+\rho(w)[q'-q]
 =F(q)+\rho(w)[q'-q].
\]
Thus \(\rho(w)\in\partial^+F(q)\). If \(F(q)=F(\bar q)\), then
\(F(q')\le F(\bar q)=F(q)\) for every \(q'\in\Omega_\infty\), so
\(0\in\partial^+F(q)\) as well. Since \(\partial^+F(q)\) is an
intersection of weak-star closed half-spaces, it follows that
\[
 \overline{\operatorname{co}}^{\,w^*}\mathcal A(q)
 \subseteq \partial^+F(q).
\]

Conversely, let \(\lambda\in\partial^+F(q)\). For every
\(h\in\mathcal P_0\) and all sufficiently small \(t>0\),
the supergradient inequality gives
\[
 \frac{F(q+th)-F(q)}{t}\le\lambda[h].
\]
Letting \(t\downarrow0\) and using
\eqref{eq:binding-directional-response}, we obtain
\[
 \min_{\eta\in\mathcal A(q)}\eta[h]
 =F'(q;h)\le\lambda[h].
\]
Suppose, to the contrary, that
\(\lambda\notin\overline{\operatorname{co}}^{\,w^*}\mathcal A(q)\).
By weak-star separation, there exists \(h\in\mathcal P_0\) such that
\[
 \lambda[h]
 <
 \inf_{\eta\in\overline{\operatorname{co}}^{\,w^*}\mathcal A(q)}
      \eta[h]
 = \inf_{\eta\in\mathcal A(q)}\eta[h]
 = F'(q;h)
 \le \lambda[h],
\]
a contradiction, where the first equality follows from the weak-star continuity
and linearity of \(\eta\mapsto\eta[h]\). Hence
\[
 \partial^+F(q)
 \subseteq \overline{\operatorname{co}}^{\,w^*}\mathcal A(q).
\]
Together with the reverse inclusion proved above, this gives
\eqref{eq:binding-superdifferential}.

\emph{Step 3.} 
Suppose that \(F(q)<F(\bar q)\) and \(\mathcal W(q)=\{w\}\).
For \(\|h\|_\infty\) sufficiently small, choose
\(w_h\in\mathcal W(q+h)\). By
Corollary~\ref{cor:binding-stability} and uniqueness,
\(w_h\to w\) strongly in \(\mathcal X\) as \(\|h\|_\infty\to0\).
The component densities therefore converge in \(L^1\), so
\(\|\rho(w_h)-\rho(w)\|_{\mathcal P_0^*}\to0\).

Testing \(F(q+h)\) with \(w\) and \(F(q)\) with \(w_h\) gives
\[
 \rho(w_h)[h]\le F(q+h)-F(q)\le \rho(w)[h].
\]
Hence
\[
 \bigl|F(q+h)-F(q)-\rho(w)[h]\bigr|
 \le \|\rho(w_h)-\rho(w)\|_{\mathcal P_0^*}\,\|h\|_\infty
 =o(\|h\|_\infty).
\]
Thus \(F\) is Fr\'echet differentiable at \(q\), 
\(DF(q)=\rho(w)\).

The density functional distinguishes distinct nonnegative functions. In fact, assume \(\rho(w_1)=\rho(w_2)\), take \(h=(\phi,0)\), where \(\phi\) is a compactly support function, and thus belongs to \(\mathcal{P}_0\). It follows that
\[\int_{\mathbb{R}^N}\phi (|w_{1,u}|^p-|w_{2,u}|^p) dx=0.\]
Since \(\phi\) is arbitrary, we have \(w_{1,u}=w_{2,u}, a.e. x\in \mathbb{R}^N\). Similarly, we obtain \(w_{1,v}=w_{2,v}\), and is thus, \(w_1=w_2\).

We now claim that if the minimiser is not unique, then \(F\) is not G\^ateaux differentiable. Indeed, we denote two different minimisers as \(w_1\) and \(w_2\), then from the analysis above, there exists a direction \(h\) such that \(\rho(w_1)[h]<\rho(w_2)[h]\). The directional response formula gives \[F'(q;h)\le \rho(w_1)[h],\quad F'(q;-h)\le -\rho(w_2)[h].\]
To sum up, we have \[
F'(q;h)
\le \rho(w_1)[h]
<\rho(w_2)[h]
\le -F'(q;-h).
\] 
a contradiction.

\emph{Step 4.}
At the threshold, If \(w\in\mathcal W(q)\), then
\(0,\rho(w)\in\mathcal A(q)\) by definition. Since \(M(w)=1\),
\(\rho(w)\ne0\), so we may choose \(h\in\mathcal P_0\) with
\(\rho(w)[h]>0\). By \eqref{eq:binding-directional-response},
\[
 F'(q;h)\le 0,
 \qquad
 F'(q;-h)\le -\rho(w)[h]<0.
\]
Thus \(F'(q;h)+F'(q;-h)<0\), which rules out G\^ateaux
differentiability at \(q\).

If \(\mathcal W(q)=\varnothing\), then \(\mathcal A(q)=\{0\}\).
The same formula gives \(F'(q;h)=F'(q;-h)=0\) for every
\(h\in\mathcal P_0\). Hence the two-sided derivative in every
direction is zero, so \(F\) is G\^ateaux differentiable at \(q\)
with \(DF(q)=0\).

To rule out Fr\'echet differentiability in the latter case,
fix $W\in\mathcal W(\bar q)$ and choose
$\chi\in C_c^\infty$, $0\le\chi\le1$, $\|\chi\|_\infty=1$,
so that
\[
T=\int\chi(|W_u|^p+|W_v|^p) dx>0.
\]
By \eqref{eq:compact-potential-energy},
\(\Phi_q(\mathsf U_yW)\to\Phi_{\bar q}(W)=F(q)\) as
\(|y|\to\infty\). Thus, given \(s_j\downarrow0\), we may choose
\(|y_j|\to\infty\) so that
\[
 \Phi_q(\mathsf U_{y_j}W)\le F(q)+\tfrac12s_jT.
\]
Set
\begin{equation}
 h_j=-s_j(\mathsf U_{y_j}\chi,\mathsf U_{y_j}\chi).
 \label{eq:binding-moving-wells}
\end{equation}
Then \(h_j\in\mathcal P_0\), \(\|h_j\|_\infty=s_j\), and
\(q+h_j\in\Omega_\infty\) for all sufficiently large \(j\).
Testing with \(\mathsf U_{y_j}W\) and changing variables gives
\begin{equation}
 \begin{aligned}
 F(q+h_j)-F(q)
 &\le \Phi_{q+h_j}(\mathsf U_{y_j}W)-F(q)\\
 &=\Phi_q(\mathsf U_{y_j}W)-F(q)-s_jT
 \le -\tfrac12s_jT.
 \end{aligned}
 \label{eq:binding-nonfrechet-defect}
\end{equation}
Consequently,
\[
 \frac{|F(q+h_j)-F(q)|}{\|h_j\|_\infty}\ge \frac{T}{2}>0.
\]
Since the G\^ateaux derivative at \(q\) is zero, this rules out
Fr\'echet differentiability. We complete the proof.

\begin{remark}
The binding region \(\{q\in\Omega_\infty:F(q)<F(\bar q)\}\) is open
by continuity of \(F\). It is dense because, at every threshold
potential \(q\), \eqref{eq:binding-nonfrechet-defect} gives
\(h_j\to0\) in \(L^\infty\) with \(F(q+h_j)<F(\bar q)\).
\end{remark}

\begin{corollary}[Attainment and selection at the threshold]
\label{cor:binding-open-dense}
Suppose \(F(q)=F(\bar q)\).
\begin{enumerate}
\item[(i)] The set \(\mathcal W(q)\) is nonempty if and only if
there exists a compactly supported \(h\in\mathcal P_0\) such that
\(F'(q;h)<0\).

\item[(ii)] If \(h\in\mathcal P_0\) satisfies \(F'(q;h)<0\), then
\(F(q+th)<F(\bar q)\) for all sufficiently small \(t>0\).
For every sequence \(t_j\downarrow0\) and every choice of
\(w_j\in\mathcal W(q+t_jh)\), a subsequence converges strongly
in \(\mathcal X\) to some \(w\) satisfying
\begin{equation}
 w\in\mathcal W(q),
 \qquad \rho(w)[h]=F'(q;h).
 \label{eq:binding-threshold-selection}
\end{equation}
\end{enumerate}
\end{corollary}

\begin{proof}
\emph{(i)} Let \(w=(u,v)\in\mathcal W(q)\). Since \(M(w)=1\),
there exists \(R>0\) such that
\(\int_{B_R}(|u|^p+|v|^p)\,dx>0\).
Taking \(h=-(\mathbf 1_{B_R},\mathbf 1_{B_R})\) in
\eqref{eq:binding-directional-response}, we obtain
\[
 F'(q;h)\le \rho(w)[h]
 =-\int_{B_R}(|u|^p+|v|^p)\,dx<0.
\]
Conversely, if \(\mathcal W(q)=\varnothing\), then
\(\mathcal A(q)=\{0\}\), so the same formula gives
\(F'(q;h)=0\) for every \(h\in\mathcal P_0\).

\emph{(ii)} Since \(F(q)=F(\bar q)\) and \(F'(q;h)<0\),
\[
 F(q+th)=F(\bar q)+tF'(q;h)+o(t)<F(\bar q)
 \qquad (t\downarrow0).
\]
Thus strict binding holds at \(q+th\) for all sufficiently small
\(t>0\). By Theorem~\ref{thm:binding-minimisers}, the Nehari level
\(d_{q+th,0}\) is attained by a nonnegative ground state, whose
normalisation belongs to \(\mathcal W(q+th)\).

Let \(t_j\downarrow0\) and \(w_j\in\mathcal W(q+t_jh)\).
Uniform coercivity gives boundedness of \((w_j)\) in \(\mathcal X\),
hence boundedness of \(\rho(w_j)[h]\). By minimality at \(q+t_jh\),
\[
 F(q)\le \Phi_q(w_j)
 =F(q+t_jh)-t_j\rho(w_j)[h]\longrightarrow F(q).
\]
Thus \((w_j)\) is a minimising sequence at \(q\). Rearranging
the same inequality gives
\[
 \rho(w_j)[h]\le
 \frac{F(q+t_jh)-F(q)}{t_j}\longrightarrow F'(q;h)<0.
\]
By Step~1 of Theorem~\ref{thm:binding-density-response}, an escaping
subsequence would satisfy \(\rho(w_j)[h]\to0\), contradicting this
bound. Hence, after passage to a subsequence, \(w_j\to w\) strongly
in \(\mathcal X\) for some \(w\in\mathcal W(q)\).
The component densities then converge in \(L^1\), so
\[
 \rho(w)[h]=\lim_{j\to\infty}\rho(w_j)[h]
 \le F'(q;h)\le \rho(w)[h],
\]
where the last inequality follows from
\eqref{eq:binding-directional-response}. This proves
\eqref{eq:binding-threshold-selection}.
\end{proof}

\begin{appendices}
\renewcommand{\theequation}{\Alph{section}.\arabic{equation}}

\section{Uniform passage from Nehari minimization to small residuals}
\label{app:uniform-variational-recovery}

Throughout this appendix, we assume \eqref{hyp:continuous} and
use the reference norm \(\|\cdot\|_{\bar q}\) on \(\mathcal X\).
The following lemma gives the variational estimate used in the
binding--escape alternative without requiring a uniform bound
on \(DG_{q,\varepsilon}\).

\begin{lemma}[Uniform Ekeland lifting]
\label{lem:appendix-uniform-lifting}
Fix \(q\) and \(\sup\limits_{0<\varepsilon\le 1/2}d_{q,\varepsilon}<E\).
There exist \(C>0\) and \(\delta_0>0\), independent of
\(0<\varepsilon\le 1/2\), with the following property.
If \(0<\delta<\delta_0\) and
\(\widehat z\in\mathcal N_{q,\varepsilon}\) satisfies
\[
 J_{q,\varepsilon}(\widehat z)
 \le d_{q,\varepsilon}+\delta^2,
\]
then there exists \(\zeta\in\mathcal N_{q,\varepsilon}\) such that
\begin{equation}
 \begin{gathered}
 \|\zeta-\widehat z\|_{\bar q}\le\delta,
 \qquad
 d_{q,\varepsilon}\le J_{q,\varepsilon}(\zeta)
                    \le J_{q,\varepsilon}(\widehat z),\\
 \|J_{q,\varepsilon}'(\zeta)\|_{\mathcal X^*}\le C\delta.
 \end{gathered}
 \label{eq:appendix-uniform-lifting}
\end{equation}
\end{lemma}

\begin{proof}
For fixed \(\varepsilon>0\), the Nehari set
\(\mathcal N_{q,\varepsilon}\) is complete in the reference norm.
Indeed, if a sequence in this set converges strongly to \(z\),
then \(G_{q,\varepsilon}(z)=0\) by continuity, and the uniform
lower bound on the Nehari norm excludes \(z=0\).
If \(M(z)=0\), then \(L_\varepsilon(z)=R_\varepsilon(z)=0\),
so \(G_{q,\varepsilon}(z)=S_q(z)>0\), a contradiction.
Thus \(z\in\mathcal N_{q,\varepsilon}\).

Choose \(\delta_0>0\) so that
\[
 \delta_0^2<E-\sup_{0<\varepsilon\le 1/2}d_{q,\varepsilon}.
\]
Ekeland's variational principle gives
\(\zeta\in\mathcal N_{q,\varepsilon}\) satisfying the distance
and energy estimates in \eqref{eq:appendix-uniform-lifting} and
\[
 J_{q,\varepsilon}(w)\ge J_{q,\varepsilon}(\zeta)
                  -\delta\|w-\zeta\|_{\bar q}
 \qquad (w\in\mathcal N_{q,\varepsilon}).
\]
Applying this inequality along \(C^1\) curves in the Nehari
hypersurface through \(\zeta\), in both tangent directions, gives
\begin{equation}
 \bigl|\langle J_{q,\varepsilon}'(\zeta),k\rangle\bigr|
 \le\delta\|k\|_{\bar q}
 \qquad\text{whenever }DG_{q,\varepsilon}(\zeta)[k]=0.
 \label{eq:appendix-tangential-ekeland}
\end{equation}

Indeed, if \(DG_{q,\varepsilon}(\zeta)[k]=0\), then \(k\) is
tangent to the \(C^1\) hypersurface \(\mathcal N_{q,\varepsilon}\)
at \(\zeta\). Hence there is a \(C^1\) curve \(\gamma\) in
\(\mathcal N_{q,\varepsilon}\) with \(\gamma(0)=\zeta\) and
\(\gamma'(0)=k\). Taking \(w=\gamma(t)\) in Ekeland's inequality
and dividing by \(t>0\), we obtain
\[
 \frac{J_{q,\varepsilon}(\gamma(t))-J_{q,\varepsilon}(\zeta)}{t}
 \ge -\delta
 \left\|\frac{\gamma(t)-\gamma(0)}{t}\right\|_{\bar q}.
\]
As \(t\downarrow0\), the left-hand side tends by the chain rule to
\(\langle J_{q,\varepsilon}'(\zeta),k\rangle\), while the norm on
the right tends to \(\|k\|_{\bar q}\). Thus
\[
 \langle J_{q,\varepsilon}'(\zeta),k\rangle
 \ge -\delta\|k\|_{\bar q}.
\]
Applying the same argument to the reversed curve \(t\mapsto\gamma(-t)\),
whose velocity at zero is \(-k\), gives the upper bound
\(\langle J_{q,\varepsilon}'(\zeta),k\rangle\le\delta\|k\|_{\bar q}\).
This proves \eqref{eq:appendix-tangential-ekeland}.

Since \(J_{q,\varepsilon}(\zeta)\le E\), the uniform bound on
Nehari sublevels gives \(\|\zeta\|_{\bar q}\le B\), independently
of \(\varepsilon\). The pointwise estimate
\(\lvert D_\xi r_\varepsilon(\xi)\rvert\le C|\xi|^{p-1}\)
and H\"older's inequality therefore give
\begin{equation}
 \|DR_\varepsilon(\zeta)\|_{\mathcal X^*}\le C_B.
 \label{eq:appendix-radial-density-derivative}
\end{equation}
Differentiating
\(G_{q,\varepsilon}=pJ_{q,\varepsilon}-(2/p)R_\varepsilon\)
yields
\begin{equation}
 DG_{q,\varepsilon}
 =pJ_{q,\varepsilon}'-\frac{2}{p}DR_\varepsilon.
 \label{eq:appendix-constraint-derivative}
\end{equation}
Moreover, radial transversality and the uniform positive lower
bound on the Nehari levels give
\begin{equation}
 -DG_{q,\varepsilon}(\zeta)[\zeta]
 \ge 2R_\varepsilon(\zeta)
 =p^2J_{q,\varepsilon}(\zeta)
 \ge p^2d_{q,\varepsilon}\ge b_*>0,
 \label{eq:appendix-uniform-transversality}
\end{equation}
where \(b_*\) is independent of \(\varepsilon\).

For \(h\in\mathcal X\), set
\[
 Ph=h-\frac{DG_{q,\varepsilon}(\zeta)[h]}
              {DG_{q,\varepsilon}(\zeta)[\zeta]}\,\zeta.
\]
Then \(DG_{q,\varepsilon}(\zeta)[Ph]=0\). Since
\(\langle J_{q,\varepsilon}'(\zeta),\zeta\rangle=0\), we also have
\[
 \langle J_{q,\varepsilon}'(\zeta),h\rangle
 =\langle J_{q,\varepsilon}'(\zeta),Ph\rangle.
\]
Writing \(K=\|J_{q,\varepsilon}'(\zeta)\|_{\mathcal X^*}\)
and using
\eqref{eq:appendix-radial-density-derivative}--%
\eqref{eq:appendix-uniform-transversality}, we obtain
\[
 \|Ph\|_{\bar q}
 \le \left(1+\frac{B}{b_*}
                   \left(pK+\frac{2}{p}C_B\right)\right)
        \|h\|_{\bar q}
 \le C_0(1+K)\|h\|_{\bar q}.
\]
Thus \eqref{eq:appendix-tangential-ekeland}, applied to \(Ph\),
gives \(K\le C_0\delta(1+K)\). Reducing \(\delta_0\) so that
\(C_0\delta_0\le 1/2\), we conclude that \(K\le 2C_0\delta\).
\end{proof}

\section{Radial recovery and logarithmic densities}
\label{app:entropy-recovery}

\begin{lemma}[Recovery on a fixed finite-energy ray]
\label{lem:appendix-ray-recovery}
For every \(z\in\mathcal D_q\) with \(M(z)>0\),
\begin{equation}
 \begin{gathered}
 t_{q,\varepsilon}(z)\longrightarrow t_{q,0}(z),\\
 J_{q,\varepsilon}\bigl(t_{q,\varepsilon}(z)z\bigr)
 \longrightarrow J_q\bigl(t_{q,0}(z)z\bigr)
 \qquad (\varepsilon\downarrow0).
 \end{gathered}
 \label{eq:appendix-ray-recovery}
\end{equation}
\end{lemma}

\begin{proof}
For each fixed \(t>0\), we have \(tz\in\mathcal D_q\).
Monotone regularization gives \(L_\varepsilon(tz)\to L(tz)\),
while \(0\le r_\varepsilon(tz)\le m(tz)\) gives
\(R_\varepsilon(tz)\to M(tz)\) by dominated convergence.
Consequently,
\[
 t^{-p}G_{q,\varepsilon}(tz)\longrightarrow t^{-p}G_{q,0}(tz).
\]
Given \(0<a<t_{q,0}(z)<b\), the limiting expression is positive
at \(a\) and negative at \(b\). The same signs therefore hold
for all sufficiently small \(\varepsilon>0\). By uniqueness
of the radial zero,
\[
 a<t_{q,\varepsilon}(z)<b.
\]
Letting \(a\) and \(b\) approach \(t_{q,0}(z)\) proves the
convergence of the scales.

Fix such \(a,b\). For all sufficiently small \(\varepsilon\),
the scale \(t_{q,\varepsilon}(z)\) lies in \([a,b]\).
For \(a\le t\le b\) and any fixed \(0<\sigma<1\), the
positive-growth estimate and the termwise bound
\[
 [\log(t^2s^2+\varepsilon)]_-
 \le [\log s^2]_-+2|\log t|
 \qquad (s\ne0)
\]
give
\[
 \ell_\varepsilon^+(tz)+\ell_\varepsilon^-(tz)
 \le C_{a,b,\sigma}
       \bigl(|z|^{p+\sigma}+\ell_0^-(z)+m(z)\bigr).
\]
The right-hand side is integrable since \(z\in\mathcal D_q\).
Dominated convergence at \(t=t_{q,\varepsilon}(z)\) now gives
convergence of the logarithmic integrals. Together with
\(S_q(tz)=t^pS_q(z)\), this proves the energy convergence
in \eqref{eq:appendix-ray-recovery}.
\end{proof}

\begin{lemma}[Recovery of the logarithmic and mass densities]
\label{lem:appendix-entropy-recovery}
Let \(0\le\varepsilon_j\le1/2\), \(\varepsilon_j\to0\), and
\(w_j\to W\) strongly in \(\mathcal X\), with \(W\) of finite energy.
Then
\begin{equation}
 \ell_{\varepsilon_j}^+(w_j)\to\ell_0^+(W),\qquad
 r_{\varepsilon_j}(w_j)\to m(W),\qquad
 m(w_j)\to m(W)
 \quad\text{in }L^1(\mathbb R^N).
 \label{eq:appendix-positive-mass-recovery}
\end{equation}
If, in addition,
\begin{equation}
 S_{q_j}(w_j)\to S_{q_*}(W),\qquad
 J_{q_j,\varepsilon_j}(w_j)\to J_{q_*}(W)<\infty,
 \label{eq:appendix-energy-recovery-hypotheses}
\end{equation}
then
\begin{equation}
 \ell_{\varepsilon_j}^-(w_j)\to\ell_0^-(W)
 \quad\text{in }L^1(\mathbb R^N).
 \label{eq:appendix-negative-recovery}
\end{equation}
\end{lemma}

\begin{proof}
Fix \(0<\sigma<1\). Strong convergence in \(\mathcal X\) gives
strong convergence in both \(L^p\) and \(L^{p+\sigma}\).
Thus every subsequence admits a further subsequence on which
\(w_j\to W\) almost everywhere and \(|w_j|\le g\) for some
\(g\in L^p\cap L^{p+\sigma}\). Along this subsequence the densities
converge pointwise, and
\[
 0\le\ell_{\varepsilon_j}^+(w_j)\le C_\sigma g^{p+\sigma},
 \qquad
 0\le r_{\varepsilon_j}(w_j)\le m(w_j)\le Cg^p.
\]
Dominated convergence and the subsequence criterion give
\eqref{eq:appendix-positive-mass-recovery}.

Under \eqref{eq:appendix-energy-recovery-hypotheses}, the energy
identity now gives
\[
 L_{\varepsilon_j}^-(w_j)
 =pJ_{q_j,\varepsilon_j}(w_j)-S_{q_j}(w_j)
   +L_{\varepsilon_j}^+(w_j)
 \longrightarrow L_0^-(W).
\]
The negative densities are nonnegative and converge almost everywhere
along the above subsequences. Since their integrals also converge,
Scheff\'e's lemma and the same subsequence criterion yield
\eqref{eq:appendix-negative-recovery}.
\end{proof}

\section*{Open problems}
\addcontentsline{toc}{section}{Open problems}

\paragraph{1. Convergence rates for regularized ground states and their energies.}
Assume \(d_{q,0}<d_{\bar q,0}\), and let \(\mathcal K(q)\) denote the set
of ground states of the unregularized problem. As \(\varepsilon\downarrow0\),
how fast does the energy gap \(d_{q,0}-d_{q,\varepsilon}\) tend to zero?
For a ground state \(z_\varepsilon\) of the regularized problem, can one
also determine the decay rate of
$\operatorname{dist}_{\mathcal X}(z_\varepsilon,\mathcal K(q))$?
If \(\mathcal K(q)=\{z_0\}\), this is precisely the rate of norm
convergence of \(z_\varepsilon\) to \(z_0\).

\paragraph{2. Onset of binding for a fixed local perturbation.}
Let \(q_\lambda=q_0-\lambda(\chi,\chi)\), where
\(q_0\in\Omega_\infty\) and \(0\le\chi\in C_c^\infty(\mathbb R^N)\)
is nonzero. Suppose the family crosses from the threshold into
strict binding while both potentials remain uniformly positive, and set
\[
 \lambda_{\mathrm b}
 =\inf\{\lambda\ge0:d_{q_\lambda,0}<d_{\bar q,0}\}.
\]
What conditions determine whether \(d_{q_{\lambda_{\mathrm b}},0}\)
is attained? What is the leading-order behaviour of
\(d_{\bar q,0}-d_{q_\lambda,0}\) as \(\lambda\downarrow\lambda_{\mathrm b}\)
from the binding side?

\paragraph{3. Classification of autonomous ground states.}
Can the nonnegative ground states for \(\bar q=(\bar a,\bar b)\)
be classified up to translations and the symmetries of the system?
In particular, under what conditions must their components satisfy
\(u=\kappa v\) for some constant \(\kappa>0\), and is this ratio
unique? Such a classification would identify which profiles can
escape at the first loss-of-compactness level \(d_{\bar q,0}\).

\paragraph{4. Several limiting potentials at infinity.}
Suppose \(q\) has no single limit at infinity, but every sequence
of translates \(q(\cdot+y_j)\), with \(|y_j|\to\infty\), has a
subsequence converging locally uniformly to a positive constant pair.
Does the singular profile decomposition still hold, with each escaping
profile solving the system for its own limiting pair? In particular,
is the first loss-of-compactness level the infimum of the ground-state
energies over all such limiting pairs?

\end{appendices}

\end{document}